\documentclass[journal]{IEEEtran}

\IEEEoverridecommandlockouts                              % This command is only
\usepackage{centernot}
\usepackage[noadjust]{cite}
\usepackage{longtable}
\usepackage{amsmath}
\usepackage{graphicx, amssymb}
\usepackage[dvips]{epsfig}
\usepackage{color}
\usepackage[dvipsnames]{xcolor}
\usepackage{comment}
\usepackage{balance}
\usepackage{tikz}

\newtheorem{thm}{Theorem}
\newtheorem{corr}{Corollary}
\newtheorem{pro}{Proposition}
\newtheorem{deff}{Definition}
\newtheorem{lem}{Lemma}
\newtheorem{ex}{Example}

\def\tt{\tilde \theta}

\def\L2{{\cal L}_2}

\newlength{\defbaselineskip}
\newcommand{\setlinespacing}[1]%
           {\setlength{\baselineskip}{#1 \defbaselineskip}}

\newcommand*{\Scale}[2][4]{\scalebox{#1}{$#2$}}%

\newcommand{\beq}{\begin{equation}}
\newcommand{\eeq}{\end{equation}}

\newcommand{\carre} {\hfill $\blacksquare$}

\makeatletter
\newcommand{\xleftrightarrow}[2][]{\ext@arrow 3359\leftrightarrowfill@{#1}{#2}}
\makeatother

\newcommand{\xdasharrow}[2][-->]{
\tikz[baseline=-\the\dimexpr\fontdimen22\textfont2\relax]{
\node[anchor=south,font=\scriptsize, inner ysep=1.5pt,outer xsep=8pt](x){#2};
\draw[shorten <=3.4pt,shorten >=3.4pt,dashed,#1](x.south west)--(x.south east);
}
}

\def\BibTeX{{\rm B\kern-.05em{\sc i\kern-.025em b}\kern-.08em
    T\kern-.1667em\lower.7ex\hbox{E}\kern-.125emX}}

\usepackage{longtable}
\usepackage{amsmath}
\usepackage{graphicx}
\usepackage{color}
\usepackage[dvips]{epsfig}
\usepackage{times}

\usepackage{graphicx, amssymb}
\usepackage{amsmath}

\smallskip

\usepackage{graphicx,amsmath}
\newcommand{\colvec}[2][.8]{%
  \scalebox{#1}{%
    \renewcommand{\arraystretch}{.8}%
    $\begin{bmatrix}#2\end{bmatrix}$%
  }
}

\allowdisplaybreaks

\def\BibTeX{{\rm B\kern-.05em{\sc i\kern-.025em b}\kern-.08em
    T\kern-.1667em\lower.7ex\hbox{E}\kern-.125emX}}

\begin{document}
\title{Modal Strong Structural Controllability for Networks with Dynamical Nodes}
\author{Shima Sadat Mousavi, \IEEEmembership{Member, IEEE}, Anastasios Kouvelas, \IEEEmembership{Senior Member, IEEE},  Karl H. Johansson, \IEEEmembership{Fellow, IEEE}
\thanks{This work was supported in part by the Swiss National Science Foundation (SNSF) under the project RECCE ``Real-time traffic estimation and control in a connected environment'', contract No. 200021-188622. The work by K. H. J. is supported by the Swedish Research Council and the  Knut and Alice Wallenberg Foundation.}
\thanks{S.~S.~Mousavi and A.~Kouvelas are with the Institute for Transport Planning and Systems, Department of Civil, Environmental and Geomatic Engineering, ETH Zurich, 8093, Switzerland (e-mails: mousavis@ethz.ch, kouvelas@ethz.ch). }
\thanks{K.~H.~Johansson is with the School of Electrical Engineering and Computer Science, KTH Royal Institute of Technology, 10044, Stockholm, Sweden. He is also affiliated with Digital Futures, Sweden (e-mail: kallej@kth.se).}}

\maketitle

\begin{abstract}
In this article, a new notion of modal strong structural controllability is introduced and examined for a family of linear time-invariant (LTI) networks. These  networks include structured  LTI subsystems (or dynamical nodes), whose system matrices have the same zero/nonzero/arbitrary pattern. % in the sense that some entries are fixed to be zero, some take any nonzero value, and the others are arbitrary.  
An eigenvalue (or the corresponding mode) associated with a system matrix is controllable if it can be directly influenced by the control inputs. We consider an arbitrary set $\Delta\subseteq \mathbb{C}$, and we refer to  a network as modal strongly  structurally controllable with respect to $\Delta$ if, for all systems in a specific family of LTI networks, every $\lambda\in \Delta$ is a controllable eigenvalue.   %The family of structured systems studied in this work includes networks of dynamical subsystems or the so-called dynamical nodes with the same  structure. However, 
For this family of LTI networks, not only is the zero/nonzero/arbitrary pattern of system matrices available, but also for a given $\Delta$, there might be extra information about the intersection of the spectrum associated with  some subsystems and $\Delta$.
%Besides the same zero/nonzero/arbitrary pattern of system matrices, for a given $\Delta$, there might be extra information about the intersection of the spectrum associated with  some subsystems and $\Delta$ for all systems in this family of LTI networks. 
For instance, for a $\Delta$ defined as the set of all complex numbers in the closed right-half plane, we may know that some structured subsystems are stable and have no eigenvalue in $\Delta$. In this case, modal strong structural controllability  is equivalent to  strong structural stabilizability.  Given a set $\Delta$, we first define a $\Delta$-network graph, and by introducing a coloring process of this graph, we establish a correspondence between the set of control subsystems and the so-called zero forcing sets. For  networks of one-dimensional subsystems, it is shown that the graph-theoretic condition is necessary and sufficient. %As a special case, by considering $\Delta$ as the set of all complex numbers in the right-half plane, the strong structural stabilizability can be examined.  
 We also demonstrate how with $\Delta=\{0\}$ or $\Delta=\mathbb{C}\setminus \{0\}$,  existing results on strong structural controllability can be derived through our approach. In fact,   compared to  existing works on strong structural controllability, a more restricted family of LTI networks is considered in this work, and for this family, the derived controllability conditions are less conservative.
\end{abstract}


\section{Introduction}

%In the two last decades, a surge of interest in studying networks of dynamical systems has arisen in both the control and network communities \cite{mesbahi2010graph}. Many features of networks of dynamical systems can be interpreted in terms of  properties of the corresponding graph that captures the network structure.  One of the key features of large-scale networks is  controllability \cite{kalman1963mathematical}. % One of the key features  that emerges when interacting with large-scale networks is their controllability, that is a fundamental system property \cite{kalman1963mathematical}.  
%In fact, controllability analysis of networks can shed light on  the control and understanding of  large-scale networks from a topological (structural) point of view. 
%In fact, controllability analysis of networks can shed light on  their understanding  from a topological (structural) point of view. There are two lines of research on controllability analysis of networks with linear time-invariant (LTI) dynamics. In the first approach, a network with predetermined weights of interconnections has been assumed \cite{yuan2013exact}. In this setting, the system matrix can be considered as the adjacency or the Laplacian matrix associated with the network graph \cite{godsil2010control, tanner2004controllability, aguilar2017almost, mousavi2018controllability,  mousavi2020laplacian}.  On the other hand, in the second approach, the controllability of a  network with unspecified link weights is studied \cite{lin1974structural, liu2011controllability,  mousavi2020structural}. 

In the two last decades, a surge of interest in studying networks of dynamical systems has arisen in both the control and network communities \cite{mesbahi2010graph}. Many features of networks of dynamical systems can be interpreted in terms of  properties of the corresponding graph that captures the network structure.  One of the key features of large-scale networks is  controllability \cite{kalman1963mathematical}. % One of the key features  that emerges when interacting with large-scale networks is their controllability, that is a fundamental system property \cite{kalman1963mathematical}.  
%In fact, controllability analysis of networks can shed light on  the control and understanding of  large-scale networks from a topological (structural) point of view. 
In fact, controllability analysis of networks can shed light on  their understanding  from a topological (structural) point of view. There are two lines of research on controllability  of networks with linear time-invariant (LTI) dynamics. In the first approach, a network with predetermined weights of interconnections is assumed \cite{yuan2013exact}. In this setting, the system matrix can be considered as the adjacency or the Laplacian matrix of the network  \cite{godsil2010control, tanner2004controllability, aguilar2017almost, mousavi2018controllability,  mousavi2020laplacian}.  In the second approach, the controllability of a  network with unspecified link weights is studied \cite{lin1974structural, liu2011controllability,  mousavi2020structural}.  

When the interaction strengths along the edges of a network are unknown, classical controllability tests cannot be utilised. However, in many practical cases, the underlying pattern of interconnections between the nodes of a network is available, so  controllability can be examined in a structural framework. In this framework, the two notions of \emph{weak} and \emph{strong structural controllability} have been introduced in the literature. Given a zero/nonzero pattern for the system matrices, an LTI system is called weakly structurally controllable if for \emph{almost all} values of the nonzero parameters, the system is controllable \cite{lin1974structural, shields1976structural, liu2011controllability,mousavi2018structural}. Strong structural controllability has been introduced to ensure that for \emph{all} nonzero weights of interconnections, the network remains controllable. Accordingly, in the structural framework, the network is viewed in terms of the pattern of zero and nonzero values of system matrices, and controllability is examined through combinatorial and graph-theoretic
conditions \cite{menara2017number}.

This paper deals with a new notion of modal strong structural controllability of LTI networks. Almost all  literature on strong structural controllability deals with the controllability of all \emph{systems}  with the same zero/nonzero pattern. However, in this work, we take a  different approach and  focus on the controllability of specific \emph{eigenvalues} (corresponding modes) for a specific family of LTI networks with the same  pattern. An eigenvalue of a system matrix associated with an LTI system is called controllable if it can be directly influenced by the control input, and by applying an appropriate control signal, it can be moved throughout the complex plane. The controllability of an eigenvalue can be checked through the Popov-Belevitch-Hautus (PBH) test \cite{sontag2013mathematical}. Now, consider an arbitrary subset of complex numbers $\Delta\subseteq \mathbb{C}$, that can be either discrete or continuous. We call a network modal strongly structurally controllable with respect to $\Delta$ if for all systems in a specific family of LTI networks, whose system matrices have some eigenvalue $\lambda\in \Delta$, $\lambda$ is a controllable eigenvalue.

Modal strong structural controllability, where strong structural controllability is examined
with respect to particular subspaces associated with the network,
is of importance in numerous applications. For instance, for $\Delta=\{y\in \mathbb{C}\mid\Re(y)\geq 0\}$, 
the modal strong structural controllability leads to strong structural stabilizability. 
 As an applied example, consider a platoon of heterogeneous vehicles whose dynamics can be represented by the optimal velocity model  \cite{bando1995dynamical}. Although the nonzero parameters in every vehicle's dynamics are uncertain,  due to the physical constraints, one can find that the dynamics of every vehicle can be modeled as a stable LTI system \cite{cui2017stabilizing}.  Accordingly, by considering the platoon of vehicles as a network of stable subsystems and analysing  the modal strong  structural controllability with respect to $\Delta=\{y\in \mathbb{C}\mid\Re(y)\geq 0\}$, one can examine the strong structural stabilizability of the whole network. 

Another application of modal strong structural controllability is in the case where the weights associated with self-loops in a network of one-dimensional subsystems (N1DS) cannot be any nonzero real numbers; rather, they can only take real values from a given (discrete or continuous) interval. In this case, the existing results on strong structural controllability may be very conservative. We show that the structured systems studied in the relevant works on strong structural controllability (e.g., \cite{monshizadeh2014zero, trefois2015zero, jia2020unifying}) are special cases of the ones we analyze in the current work. In fact, we show how for $\Delta=\{0\}$ and $\Delta=\mathbb{C}\setminus\{ 0\}$, the existing results can be reproduced through our approach, and thus, our results can establish a unifying framework for strong structural controllability of LTI networks.

In what follows, we first present a brief description of the relevant  literature on strong structural controllability, and then, we discuss  the main contributions of this work.

\subsection{Literature Review}

Strong structural controllability was first introduced by Mayeda and Yamada in \cite{mayeda1979strong}.  
Subsequently, this notion of controllability has been explored by a number of works using spanning cycles~\cite{Reinschke1992, jarczyk2011strong,pequito2013framework,reissig2014strong} and 
constrained $t$-matchings~\cite{chapman2013strong}.
Moreover, in \cite{monshizadeh2014zero}, a one-to-one correspondence
between zero forcing sets of a {loop-free graph} and strong structural controllability has been stated;   
the latter result has been then extended in \cite{shima,mousavi2018structural,monshizadeh2015strong, van2017distance,mousavi2019strong}. {Additionally in \cite{trefois2015zero}, a link between constrained matchings and zero forcing sets has been established, and conditions for strong structural controllability of {loop directed networks} have been obtained in terms of zero forcing sets.} Later, in \cite{jia2019strong}, the results of \cite{mousavi2018structural} have been extended to the case where some nonzero entries of system matrices are restricted to have identical values. 
 
%More recently, in \cite{mousavi2017robust,mousavi2020strong,popli2019selective,jia2020unifying}, the problem of strong structural controllability has been studied, in the cases where not only the nonzero weights of connections are uncertain, but  the knowledge about the structure of the associated graph is  incomplete as well. This means that in these networks, examples of which can be seen among social and biological systems \cite{kossinets2006effects, clauset2008hierarchical}, it is uncertain whether some specific edges in the $\Delta$-network graph  exist or not. Thus, in these cases, rather than a zero/nonzero structure, we need to consider a zero/nonzero/arbitrary pattern for the system matrices. In such a pattern, some entries are fixed to zero, some are arbitrary nonzero real numbers, and others can be any real numbers (including zero and nonzero). In \cite{jia2020unifying}, considering a zero/nonzero/arbitrary pattern, a strong structural controllability condition has been provided in terms of zero forcing sets. 
%

In \cite{mousavi2017robust,mousavi2020strong,popli2019selective,jia2020unifying}, the problem of strong structural controllability has been studied in the cases where not only the nonzero weights of connections are uncertain, but also the structural characteristics of the associated graph are not completely specified. In such networks, examples of which are found in social and biological systems \cite{kossinets2006effects, clauset2008hierarchical}, the existence of some of the edges in the network graph is unknown. In these cases, we need to consider a zero/nonzero/arbitrary pattern for the system matrices instead of a zero/nonzero structure. In such a pattern, the entries are of three types: fixed to zero, arbitrary nonzero real numbers, and arbitrary real numbers (including zero and nonzero). In the corresponding graph, the links associated with nonzero and arbitrary entries are, respectively, represented by solid and dotted edges.   Networks owing a zero/nonzero/arbitrary patterns are considered in \cite{jia2020unifying} where a strong structural controllability condition is determined in terms of zero forcing sets.

The notion of zero forcing set, {related to a particular coloring of graph nodes}, was first proposed in \cite{work2008zero} to study the minimum rank problem. A zero forcing set is a subset of black nodes in a graph that can force its other nodes to be black through applying a coloring rule as many times as possible. The minimum rank problem was motivated by the inverse eigenvalue problem of a graph, whose goal is to obtain information about the possible eigenvalues of a family of patterned matrices; the first step towards this aim is to find the maximum multiplicity of an arbitrary eigenvalue of all matrices in this family.  
{With this in mind, using the notion of zero forcing sets,  \cite{work2008zero} and \cite{barioli2009minimum} have presented upper bounds on the maximum multiplicity of the \emph{zero eigenvalue} for loop-free undirected and loop directed graphs, respectively. 

%In this direction, for a given set of eigenvalues $\Delta$, by introducing the new notion of \emph{$\Delta$-network graphs}, we provide an upper bound on the maximum multiplicity of \emph{all eigenvalues in $\Delta$} associated with a {loop} directed network of dynamical nodes with a zero/nonzero/arbitrary structure.} 

Most of the works on weak and strong structural controllability have been dedicated to the networks, whose any node state is scalar. However, recently, the problem of weak structural controllability has been studied for networks with dynamical nodes, where every node is made of a linear dynamic system \cite{commault2019structural}. Moreover, in \cite{jia2020scalable}, the strong structural controllability of a network of single-input single output structured subsystems has been investigated, and it has been shown that an LTI network is strongly structurally controllable if and only if an associated  network of an order at most twice the number of the included subsystems is strongly structurally controllable.

\subsection{Main Contributions}

In this paper, we introduce a new notion of modal strong structural controllability and analyze it for a family of LTI networks, which include structured LTI subsystems. %Almost all the works in the literature on strong structural controllability deal with the controllability of all \emph{systems}  with the same zero/nonzero or zero/nonzero/arbitrary pattern system matrices. However, in this work, we take a  different approach and  focus on the controllability of specific \emph{eigenvalues} (or the corresponding modes) for a specified family of LTI networks with the same structure. An eigenvalue of a system matrix associated with an LTI system is called controllable if it can be directly influenced by the control input, and by applying an appropriate control signal, it can be moved throughout the complex plane. The controllability of an eigenvalue can be tested through the Popov-Belevitch-Hautus (PBH) test \cite{sontag2013mathematical}. Now, consider an arbitrary subset of complex numbers $\Delta\subseteq \mathbb{C}$, that can be either discrete or continuous. We call a network modal strongly structurally controllable with respect to $\Delta$ if for all systems in a specific family of LTI networks, whose system matrices have some eigenvalue $\lambda\in \Delta$, $\lambda$ is a controllable eigenvalue. 
%
%
%In this paper, we study an LTI network including LTI subsystems. However, our approach defers from \cite{jia2020scalable}, mainly because other than the zero/nonzero/arbitrary structure of the network, we assume that other restrictive information about the subsystems may be available. In fact, based on our assumption, there can be some information about the intersection of a given subset of complex numbers $\Delta$ and the spectra of some of the subsystems included in the network. For instance, let $\Delta$ be the set of all eigenvalues in the closed right-half complex plane, and assume that some of the subsystems are stable, in the sense that all of their associated eigenvalues are in the open left-half plane. Therefore, the family of networks includes all systems of the same zero/nonzero/arbitrary structure, with this extra property that some of the subsystems are stable.% that have extra features as well.
%
We assume that, other than the zero/nonzero/arbitrary pattern of system matrices, some restrictive information about the subsystems included in the LTI network may be available. In fact, based on our assumption, we might have information about the intersection of a given subset of complex numbers $\Delta$ and the spectra of some of the subsystems. In this direction, a $\Delta$-characteristic vector %associated with a network of dynamical subsystems
is defined, which captures the information about the intersection of the spectrum associated with any subsystem and $\Delta$. Moreover, a $\Delta$-specified pattern class associated with an LTI network includes  all system matrices of the same zero/nonzero/arbitrary pattern, which have the extra properties described by the $\Delta$-characteristic vector.  For instance, let $\Delta$ be the set of all eigenvalues in the closed right-half complex plane, and assume that some of the subsystems are stable, in the sense that all of their associated eigenvalues are in the open left-half plane. Therefore, the specific family of LTI networks includes all systems of the same zero/nonzero/arbitrary pattern, with this extra property that some given subsystems are stable.
 We also introduce a notion of \emph{$\Delta$-network graph} associated with a $\Delta$-characteristic vector and  propose a coloring process applied to this graph. The main contributions of this work are:

}
 
 %For a given subset of complex numbers $\Delta$, in Section II, we define the problem of modal strong structural controllability for LTI networks of dynamical subsystems and present new definitions, e.g, $\Delta$-characteristic vectors and $\Delta$-specified pattern classes for networks. A $\Delta$-characteristic vector associated with a network of dynamical subsystems captures the information about the intersection of the spectrum associated with any subsystem and $\Delta$. A $\Delta$-specified pattern class associated with a network is the set of all system matrices of the same structure, which have the extra properties described by the $\Delta$-characteristic vector.   We also introduce a notion of \emph{$\Delta$-network graph} associated with a $\Delta$-characteristic vector and  propose a coloring process applied to this graph in Section III.  

\begin{enumerate}
%   \item 
   % \item For a given subset of complex numbers $\Delta$, in Section II, we define the problem of modal strong structural controllability for LTI networks of dynamical subsystems and present new definitions, e.g, $\Delta$-characteristic vectors and $\Delta$-specified pattern classes for networks. A $\Delta$-characteristic vector associated with a network of dynamical subsystems captures the information about the intersection of the spectrum associated with any subsystem and $\Delta$. A $\Delta$-specified pattern class associated with a network is the set of all system matrices of the same structure, which have the extra properties described by the $\Delta$-characteristic vector.   We also introduce a notion of \emph{$\Delta$-network graph} associated with a $\Delta$-characteristic vector and  propose a coloring process applied to this graph in Section III. 
    \item We develop a correspondence between zero forcing sets of a $\Delta$-network graph and sets of control subsystems, rendering it modal strongly structurally controllable (Theorem 1). Since the family of LTI networks in this work is more specified than the family of networks with only the same structure, the set of control subsystems obtained through our approach can be of a smaller cardinality. Moreover,   given a $\Delta$ that includes at least one real number,  we establish a necessary and sufficient condition for modal strong structural controllability of N1DSs (Theorem 2). As a particular case, we show how the derived graph-theoretic conditions can be utilised for the analysis of the strong structural stabilizability of LTI networks. 
    \item 
    Finding the maximum geometric multiplicity of an arbitrary eigenvalue for all matrices of the same pattern is of interest as one goal of the minimum rank problem. In this direction, we take a step forward, and for any given set $\Delta$, we provide an upper bound on the maximum geometric multiplicity of \emph{all eigenvalues in $\Delta$} associated with the system matrices in a $\Delta$-specified pattern class of an LTI network (Proposition 2). %{loop} directed network of dynamical nodes with the same zero/nonzero/arbitrary structure. 
    We also derive combinatorial conditions, under which, no matrix in a $\Delta$-specified pattern class associated with an N1DS has any eigenvalue in $\Delta$ (Theorem 5). For example, through a graph-theoretic condition, one can see whether any system matrix associated with an N1DS is stable or not. 
    \item Although there exist results in the literature on the full rankness of a pattern matrix with a zero/nonzero \cite{work2008zero,barioli2009minimum} or zero/nonzero/arbitrary structure \cite{jia2020unifying}, we provide an equivalent combinatorial condition, which can be tested by considering the corresponding bipartite graph (Proposition 3). This can facilitate forming of $\Delta$-network graph associated with a $\Delta$-specified pattern class from the corresponding global graph. 
    %\item We show how our results can extend the existing knowledge on strong structural controllability of LTI networks. We provide numerous examples to better illustrate definitions and the obtained results.
\end{enumerate}

Finally, we note that the strong structural controllability of LTI networks including LTI subsystems has been  investigated in \cite{jia2020scalable} as well; however,  our approach defers from \cite{jia2020scalable}, mainly because other than the zero/nonzero/arbitrary pattern of the system matrices, we assume that other restrictive information about the subsystems may be available. In fact, a set of system matrices represented by a $
\Delta$-specified pattern class is smaller and more restrictive than a set of system matrices with only the same zero/nonzero/arbitrary pattern. Thus, the controllability conditions that we derive are less conservative than the existing results on strong structural controllability.  Additionally, we show how our results can extend the existing knowledge on strong structural controllability of LTI networks. We also provide numerous examples to better illustrate definitions and the obtained results.

\subsection{Outline}

The paper is organized as follows. Preliminaries, including all necessary definitions and problem formulation, are presented in Section~\ref{prelim}. In Section~\ref{Delta_color}, we introduce a $\Delta$-network graph associated with a given set $\Delta$. Moreover, the definition of the coloring process and zero forcing sets are presented in this section. In Section~\ref{results}, the main results  are established.  Finally, Section~\ref{concl} concludes the paper.

%%%%%%%%%%%%%%%%%%%%%%%%%%%%%%%%%%%%%%%%%%%%%%%%%%%%%%%
%%%%%%%%%%%%%%%%%%%%%%%%%%%%%%%%%%%%%%%%%%%%%%%%%%%%%%%%%%%

\section{Preliminaries}\label{prelim}

The set of real and complex numbers are denoted by $\mathbb{R}$ and $\mathbb{C}$, respectively.  
The $i$-th element of the vector $v$ is designated by $v(i)$, and $M({i,j})$ is the entry in row $i$ and column $j$ of $M$. Moreover, for $i_2\geq i_1$ and $j_2\geq j_1$, $M(i_1:i_2,j_1:j_2)$ is a submatrix of $M$ formed from the successive rows $i_1, \ldots, i_2$ and successive columns $j_1,\ldots, j_2$. A subvector $v(X)$ is comprised of $v(i)$, for $i\in X$, ordered lexiographically. $\mathbf{1}_n$ denotes the vector of all ones in $\mathbb{R}^n$.
We denote the $n\times n$ identity matrix by $I_n$ and represent its $j$-th column by $e_j$.  The cardinality of a set $S$ is designated by $\mid S\mid$.

\subsection{Definitions}

\textit{Pattern Matrices:} A pattern matrix $\mathcal{A}$ is a matrix whose entries are chosen from the set of symbols $\{0,*,?\}$. A pattern class of a $q\times p$ pattern matrix $\mathcal{A}$, denoted by $\mathcal{P}(\mathcal{A})$, is defined as $\mathcal{P}(\mathcal{A})=\{A\in \mathbb{R}^{q\times p} \mid A({i,j})=0 \:\:\:\:\:\mbox{if}\:\:\:\:\: \mathcal{A}({i,j})=0, \:\:\:\mbox{and}\:\: A({i,j})\neq 0 \:\:\:\mbox{if} \:\:\mathcal{A}({i,j})=* \}$. By this definition, if $\mathcal{A}({i,j})=?$,  $A({i,j})$ can be any arbitrary real number, including zero and nonzero. 
We say that a $q\times p$  pattern matrix $\mathcal{A}$ %with $q\leq p$,  
has full row rank if every $A\in\mathcal{P}(\mathcal{A})$ has full row rank.

\emph{Graphs:} Let $\mathcal{A}\in \{0,*,?\}^{q\times p}$ be a pattern matrix, where $l=\max (q,p)$.  
Then, we define an associated graph $G(\mathcal{A})=(V,E)$ with the node set $V=\{v_1,\ldots,v_l\}$ and edge set $E\subseteq V\times V$. We have $(v_j, v_i) \in E$ if and only if $\mathcal{A}(i,j)\neq 0$. %$\mathcal{A}(i,j)=*$  or $\mathcal{A}({i,j})=?$.
Then, there is an edge from node $v_j$ to node $v_i$. In this case, node $v_i$ (respectively, node $v_j$) is said to be an out-neighbor (respectively, in-neighbor) of node $v_j$ (respectively, node $v_i$). We denote by $N_{\rm out}(v_j)$ %
the set of out-neighbors
 of node $v_j$. For an \emph{undirected} graph, $(v_i,v_j)\in E$ if and only if $(v_j,v_i)\in E$, and thus the associated pattern matrix $\mathcal{A}$ should be symmetric in this case.  Note that  a graph ${G({\mathcal{A}})}$ can contain loops as $(v_i,v_i)$ (i.e., self-loops) for some $v_i\in V$.  Let $E=E^*\cup E^?$, where $(v_{j},v_i)\in E^*$ if and only if $\mathcal{A}(i,j)=*$, and $(v_{j},v_i)\in E^?$ if and only if $\mathcal{A}(i,j)=?$. If $(v_{j},v_i)\in E^*$ (respectively, $(v_{j},v_i)\in E^?$), $v_i$ is called a \emph{strong} (respectively, \emph{weak}) out-neighbor of $v_j$. To distinguish the edges in $E^*$ and $E^?$, we show them by solid and dotted arrows, respectively.

 \emph{Bipartite graphs:} Given a pattern matrix $\mathcal{A}\in \{0,*,?\}^{q\times p}$, one can associate a  bipartite graph  as $G_b(\mathcal{A})=(V_r,V_c,E_b)$ with $V_r=\{v^r_1,\ldots,v^r_q\}$ and $V_c=\{v^c_1,\ldots,v^c_p\}$. The set $E_b$ is a subset of  edges from nodes of $V_c$ to nodes of  $V_r$, where $(v^c_i, v^r_j) \in E_b$ if and only if $\mathcal{A}(j,i)\neq 0$. We let $E_b=E^*_b\cup E^?_b$, where $(v^c_i, v^r_j) \in E_b^*$ if and only if $\mathcal{A}(j,i)=*$, and $(v^c_i, v^r_j) \in E_b^?$ if and only if $\mathcal{A}(j,i)=?$. %
 The edges in $E^*$ and $E^?$ are shown by solid and dotted arrows, respectively.

\emph{Network graphs:} Consider a network $\mathcal{N}$ with $n$ nodes. An associated network graph ${G_{\mathcal{N}}}$ is denoted by $G_{{\mathcal{N}}}=(V_{{{\mathcal{N}}}},E_{{{\mathcal{N}}}})$, where $V_{\mathcal{N}}=\{1,2,\ldots,n\}$ is the node set, and $E_{\mathcal{N}}\subseteq V_{\mathcal{N}}\times V_{\mathcal{N}}$ is the edge set of the graph. In the next section of the paper, we provide a detailed description of the edge set of a network graph regarding our problem. Essentially, we assume that a set of numbers $\Delta$ is given, and we demonstrate how one can form an associated \emph{$\Delta$-network graph $G_{\mathcal{N}}^{\Delta}$}. 

\emph{Node graphs:} Every node $i$ in a network graph ${G_{\mathcal{N}}}$, $i=1,2,\ldots, n$, is indeed a ``super node", in the sense that it can represent a system itself. Thus, it has $l_i$ internal vertices and is illustrated by a node graph $G_{i}=(V_{i},E_{i})$, where $V_{i}=\{i^1,\ldots,i^{l_i}\}$. To make a distinction, we refer to $V_{i}$ as the set of ``vertices'' of $G_{i}$, while $V_{\mathcal{N}}$ is the set of ``nodes'' of ${G_{\mathcal{N}}}$.  For a subset of nodes $Z\subseteq V_{\mathcal{N}}$, we define the \emph{set of vertices} of $Z$ as $\mbox{Ver}(Z)=\bigcup_{i\in Z}V_i.$

\emph{Eigenvalues:} Let $\Lambda(A)$ denote the spectrum or the set of eigenvalues of matrix $A$.  
The geometric multiplicity of eigenvalue $\lambda\in \Lambda(A)$, which is denoted by $\psi_{A}(\lambda)$, is the dimension of the subspace $\mathcal{V}_{A}(\lambda)=\{\nu\in \mathbb{R}^n \mid \nu^TA=\lambda \nu^T \}$. For a subset $\mathcal{M}\subseteq \Lambda(A)$, the maximum geometric multiplicty of the eigenvalues of $A$ belonging to $\mathcal{M}$ is defined as $\Psi_{\mathcal{M}}(A)=\mathrm{max}\{\psi_{A}(\lambda)\mid \lambda\in \mathcal{M}\}$.

 \emph{System matrices:} For $i=1,2,\ldots, n$, let $\mathcal{A}_{ii}\in\{0,*,?\}^{l_i\times l_i}$ be a pattern matrix. Now, 
consider a network $\mathcal{N}$, including some LTI  subsystem $i$ with system matrix $A_{ii}\in \mathcal{P}(\mathcal{A}_{ii})$, $i=1,2,\ldots, n$. 
Note that for some $1\leq i\leq n$, we may have $l_i=1$. %If $l_i=1$, for all $1\leq i\leq n$, we have a \emph{network of one-dimensional subsystems} (\emph{N1DS}).
Now, let $N=\sum_{i=1}^n l_i$. The global system matrix associated with this network is
\begin{equation}
A=\begin{bmatrix}A_{11} & A_{12} & \ldots& A_{1n}\\ A_{21} & A_{22} & \ldots& A_{2n}\\
\ddots & \ddots & \ddots & \ddots\\
A_{n1} & A_{n2} & \ldots & A_{nn}\end{bmatrix}\in \mathbb{R}^{N\times N}, 
\label{sm}
\end{equation}
and we let $A\in \mathcal{P}(\mathcal{{A}})$, for a pattern matrix $\mathcal{A}\in\{0,*,?\}^{N\times N}$. Note that for all $1\leq i,j\leq n$, we have $A_{ij}\in\mathbb{R}^{l_i\times l_j}$, and  $A_{ij}\in \mathcal{P}(\mathcal{A}_{ij})$, where $\mathcal{A}_{ij}\in\{0,*,?\}^{l_i\times l_j}$.
For $i=1,2,\ldots, n$, node $i$ is a super node, representing the subsystem $i$. Node graph $G_{i}=G(\mathcal{A}_{ii})=(V_{i}, E_{i})$ illustrates the internal dynamics of node $i$. On the other hand, 
graph $\mathcal{G}=G(\mathcal{A})=(V_{\mathcal{G}}, E_{\mathcal{G}})$, with $V_{\mathcal{G}}=\bigcup_{i=1}^n V_{i}$, is called the \emph{global graph}, representing the structure of the entire system.

\begin{figure}[t]
\centering

\tikzset{every picture/.style={line width=0.75pt}} %set default line width to 0.75pt        

\begin{tikzpicture}[x=0.75pt,y=0.75pt,yscale=-1,xscale=1]
%uncomment if require: \path (0,11616); %set diagram left start at 0, and has height of 11616

%Shape: Circle [id:dp004541612384300953] 
\draw  [fill={rgb, 255:red, 255; green, 255; blue, 255 }  ,fill opacity=1 ] (113.44,8160.35) .. controls (118,8160.34) and (121.7,8164.02) .. (121.7,8168.58) .. controls (121.71,8173.14) and (118.02,8176.84) .. (113.47,8176.85) .. controls (108.91,8176.85) and (105.21,8173.17) .. (105.2,8168.61) .. controls (105.2,8164.05) and (108.88,8160.35) .. (113.44,8160.35) -- cycle ;
%Shape: Circle [id:dp17737262064014714] 
\draw   (173.44,8160.35) .. controls (178,8160.34) and (181.7,8164.02) .. (181.7,8168.58) .. controls (181.71,8173.14) and (178.02,8176.84) .. (173.47,8176.85) .. controls (168.91,8176.85) and (165.21,8173.17) .. (165.2,8168.61) .. controls (165.2,8164.05) and (168.88,8160.35) .. (173.44,8160.35) -- cycle ;
%Shape: Circle [id:dp8853780994720417] 
\draw   (144.44,8110.35) .. controls (149,8110.34) and (152.7,8114.02) .. (152.7,8118.58) .. controls (152.71,8123.14) and (149.02,8126.84) .. (144.47,8126.85) .. controls (139.91,8126.85) and (136.21,8123.17) .. (136.2,8118.61) .. controls (136.2,8114.05) and (139.88,8110.35) .. (144.44,8110.35) -- cycle ;
%Straight Lines [id:da65864651248012] 
\draw    (113.44,8160.35) -- (136.8,8126.47) ;
\draw [shift={(138.5,8124)}, rotate = 484.59] [fill={rgb, 255:red, 0; green, 0; blue, 0 }  ][line width=0.08]  [draw opacity=0] (8.93,-4.29) -- (0,0) -- (8.93,4.29) -- cycle    ;
%Straight Lines [id:da7458409019061885] 
\draw  [dash pattern={on 4.5pt off 4.5pt}]  (121.7,8168.58) -- (162.2,8168.61) ;
\draw [shift={(165.2,8168.61)}, rotate = 180.04] [fill={rgb, 255:red, 0; green, 0; blue, 0 }  ][line width=0.08]  [draw opacity=0] (8.93,-4.29) -- (0,0) -- (8.93,4.29) -- cycle    ;
%Curve Lines [id:da010630378868456303] 
\draw    (152.48,8125.79) .. controls (164.15,8134.67) and (172.19,8147.88) .. (173.44,8160.35) ;
\draw [shift={(150,8124)}, rotate = 34.2] [fill={rgb, 255:red, 0; green, 0; blue, 0 }  ][line width=0.08]  [draw opacity=0] (8.93,-4.29) -- (0,0) -- (8.93,4.29) -- cycle    ;
%Curve Lines [id:da16516597109113285] 
\draw  [dash pattern={on 4.5pt off 4.5pt}]  (150,8124) .. controls (150.31,8132.88) and (159.46,8151.1) .. (171,8158.88) ;
\draw [shift={(173.44,8160.35)}, rotate = 207.76] [fill={rgb, 255:red, 0; green, 0; blue, 0 }  ][line width=0.08]  [draw opacity=0] (8.93,-4.29) -- (0,0) -- (8.93,4.29) -- cycle    ;
%Curve Lines [id:da030150673947198614] 
\draw [color={rgb, 255:red, 0; green, 0; blue, 0 }  ,draw opacity=1 ]   (113.47,8176.85) .. controls (107.87,8197.42) and (90.29,8195.66) .. (104.04,8174.54) ;
\draw [shift={(105.67,8172.17)}, rotate = 485.54] [fill={rgb, 255:red, 0; green, 0; blue, 0 }  ,fill opacity=1 ][line width=0.08]  [draw opacity=0] (8.93,-4.29) -- (0,0) -- (8.93,4.29) -- cycle    ;
%Curve Lines [id:da7879789938584612] 
\draw [color={rgb, 255:red, 0; green, 0; blue, 0 }  ,draw opacity=1 ] [dash pattern={on 4.5pt off 4.5pt}]  (140.33,8111.17) .. controls (131.37,8078.21) and (151.9,8085.49) .. (149.43,8107.95) ;
\draw [shift={(149,8110.83)}, rotate = 280.71] [fill={rgb, 255:red, 0; green, 0; blue, 0 }  ,fill opacity=1 ][line width=0.08]  [draw opacity=0] (8.93,-4.29) -- (0,0) -- (8.93,4.29) -- cycle    ;
%Shape: Circle [id:dp5574728551908104] 
\draw  [fill={rgb, 255:red, 255; green, 255; blue, 255 }  ,fill opacity=1 ] (271.11,8160.35) .. controls (275.66,8160.34) and (279.36,8164.02) .. (279.37,8168.58) .. controls (279.38,8173.14) and (275.69,8176.84) .. (271.13,8176.85) .. controls (266.58,8176.85) and (262.88,8173.17) .. (262.87,8168.61) .. controls (262.86,8164.05) and (266.55,8160.35) .. (271.11,8160.35) -- cycle ;
%Shape: Circle [id:dp022042543937917047] 
\draw  [fill={rgb, 255:red, 255; green, 255; blue, 255 }  ,fill opacity=1 ] (271.11,8110.01) .. controls (275.66,8110) and (279.36,8113.69) .. (279.37,8118.25) .. controls (279.38,8122.8) and (275.69,8126.5) .. (271.13,8126.51) .. controls (266.58,8126.52) and (262.88,8122.83) .. (262.87,8118.28) .. controls (262.86,8113.72) and (266.55,8110.02) .. (271.11,8110.01) -- cycle ;
%Shape: Circle [id:dp6733635847701134] 
\draw  [fill={rgb, 255:red, 255; green, 255; blue, 255 }  ,fill opacity=1 ] (330.77,8160.01) .. controls (335.33,8160) and (339.03,8163.69) .. (339.04,8168.25) .. controls (339.04,8172.8) and (335.36,8176.5) .. (330.8,8176.51) .. controls (326.24,8176.52) and (322.54,8172.83) .. (322.54,8168.28) .. controls (322.53,8163.72) and (326.22,8160.02) .. (330.77,8160.01) -- cycle ;
%Shape: Circle [id:dp8469492723614058] 
\draw  [fill={rgb, 255:red, 255; green, 255; blue, 255 }  ,fill opacity=1 ] (330.77,8109.68) .. controls (335.33,8109.67) and (339.03,8113.36) .. (339.04,8117.91) .. controls (339.04,8122.47) and (335.36,8126.17) .. (330.8,8126.18) .. controls (326.24,8126.19) and (322.54,8122.5) .. (322.54,8117.94) .. controls (322.53,8113.39) and (326.22,8109.69) .. (330.77,8109.68) -- cycle ;
%Shape: Circle [id:dp45345938235906225] 
\draw  [fill={rgb, 255:red, 255; green, 255; blue, 255 }  ,fill opacity=1 ] (143.11,8322.35) .. controls (147.66,8322.34) and (151.36,8326.02) .. (151.37,8330.58) .. controls (151.38,8335.14) and (147.69,8338.84) .. (143.13,8338.85) .. controls (138.58,8338.85) and (134.88,8335.17) .. (134.87,8330.61) .. controls (134.86,8326.05) and (138.55,8322.35) .. (143.11,8322.35) -- cycle ;
%Shape: Circle [id:dp8697679575653472] 
\draw  [fill={rgb, 255:red, 255; green, 255; blue, 255 }  ,fill opacity=1 ] (143.11,8262.01) .. controls (147.66,8262) and (151.36,8265.69) .. (151.37,8270.25) .. controls (151.38,8274.8) and (147.69,8278.5) .. (143.13,8278.51) .. controls (138.58,8278.52) and (134.88,8274.83) .. (134.87,8270.28) .. controls (134.86,8265.72) and (138.55,8262.02) .. (143.11,8262.01) -- cycle ;
%Shape: Circle [id:dp4455545957671221] 
\draw  [fill={rgb, 255:red, 255; green, 255; blue, 255 }  ,fill opacity=1 ] (270.77,8262.01) .. controls (275.33,8262) and (279.03,8265.69) .. (279.04,8270.25) .. controls (279.04,8274.8) and (275.36,8278.5) .. (270.8,8278.51) .. controls (266.24,8278.52) and (262.54,8274.83) .. (262.54,8270.28) .. controls (262.53,8265.72) and (266.22,8262.02) .. (270.77,8262.01) -- cycle ;
%Shape: Circle [id:dp5889529261676885] 
\draw  [fill={rgb, 255:red, 255; green, 255; blue, 255 }  ,fill opacity=1 ] (333.44,8261.68) .. controls (338,8261.67) and (341.7,8265.36) .. (341.7,8269.91) .. controls (341.71,8274.47) and (338.02,8278.17) .. (333.47,8278.18) .. controls (328.91,8278.19) and (325.21,8274.5) .. (325.2,8269.94) .. controls (325.2,8265.39) and (328.88,8261.69) .. (333.44,8261.68) -- cycle ;
%Curve Lines [id:da7888882370873904] 
\draw [color={rgb, 255:red, 0; green, 0; blue, 0 }  ,draw opacity=1 ]   (266.67,8111.17) .. controls (262.48,8082.54) and (279.1,8083.4) .. (275.78,8108.06) ;
\draw [shift={(275.33,8110.83)}, rotate = 280.37] [fill={rgb, 255:red, 0; green, 0; blue, 0 }  ,fill opacity=1 ][line width=0.08]  [draw opacity=0] (8.93,-4.29) -- (0,0) -- (8.93,4.29) -- cycle    ;
%Curve Lines [id:da042117754916256134] 
\draw [color={rgb, 255:red, 0; green, 0; blue, 0 }  ,draw opacity=1 ]   (326,8111.5) .. controls (321.8,8082.72) and (337.67,8081.55) .. (334.97,8108.58) ;
\draw [shift={(334.67,8111.17)}, rotate = 277.68] [fill={rgb, 255:red, 0; green, 0; blue, 0 }  ,fill opacity=1 ][line width=0.08]  [draw opacity=0] (8.93,-4.29) -- (0,0) -- (8.93,4.29) -- cycle    ;
%Straight Lines [id:da047339059009849915] 
\draw  [dash pattern={on 4.5pt off 4.5pt}]  (330.77,8160.01) -- (330.8,8129.18) ;
\draw [shift={(330.8,8126.18)}, rotate = 450.05] [fill={rgb, 255:red, 0; green, 0; blue, 0 }  ][line width=0.08]  [draw opacity=0] (8.93,-4.29) -- (0,0) -- (8.93,4.29) -- cycle    ;
%Curve Lines [id:da9247162241597631] 
\draw    (279.37,8168.58) .. controls (286.69,8162.42) and (300.5,8157.88) .. (320.06,8167.06) ;
\draw [shift={(322.54,8168.28)}, rotate = 207.31] [fill={rgb, 255:red, 0; green, 0; blue, 0 }  ][line width=0.08]  [draw opacity=0] (8.93,-4.29) -- (0,0) -- (8.93,4.29) -- cycle    ;
%Straight Lines [id:da02282153323821423] 
\draw  [dash pattern={on 4.5pt off 4.5pt}]  (282.37,8118.23) -- (322.54,8117.94) ;
\draw [shift={(279.37,8118.25)}, rotate = 359.6] [fill={rgb, 255:red, 0; green, 0; blue, 0 }  ][line width=0.08]  [draw opacity=0] (8.93,-4.29) -- (0,0) -- (8.93,4.29) -- cycle    ;
%Curve Lines [id:da6492014300892364] 
\draw [color={rgb, 255:red, 0; green, 0; blue, 0 }  ,draw opacity=1 ] [dash pattern={on 4.5pt off 4.5pt}]  (335.33,8175.5) .. controls (340.48,8194.8) and (325.14,8204.17) .. (326.46,8178.44) ;
\draw [shift={(326.67,8175.5)}, rotate = 455.19] [fill={rgb, 255:red, 0; green, 0; blue, 0 }  ,fill opacity=1 ][line width=0.08]  [draw opacity=0] (8.93,-4.29) -- (0,0) -- (8.93,4.29) -- cycle    ;
%Curve Lines [id:da8993070450908078] 
\draw  [dash pattern={on 4.5pt off 4.5pt}]  (279.04,8270.25) .. controls (288.6,8261.85) and (304.34,8257.83) .. (322.87,8268.53) ;
\draw [shift={(325.2,8269.94)}, rotate = 212.49] [fill={rgb, 255:red, 0; green, 0; blue, 0 }  ][line width=0.08]  [draw opacity=0] (8.93,-4.29) -- (0,0) -- (8.93,4.29) -- cycle    ;
%Curve Lines [id:da7195465528620832] 
\draw  [dash pattern={on 4.5pt off 4.5pt}]  (143.13,8278.51) .. controls (151.87,8285.15) and (151.71,8302.93) .. (144.33,8319.71) ;
\draw [shift={(143.11,8322.35)}, rotate = 296.05] [fill={rgb, 255:red, 0; green, 0; blue, 0 }  ][line width=0.08]  [draw opacity=0] (8.93,-4.29) -- (0,0) -- (8.93,4.29) -- cycle    ;
%Curve Lines [id:da8166056345252264] 
\draw  [dash pattern={on 4.5pt off 4.5pt}]  (143.11,8322.35) .. controls (137.27,8317.08) and (133.97,8298.6) .. (141.93,8280.99) ;
\draw [shift={(143.13,8278.51)}, rotate = 477.33] [fill={rgb, 255:red, 0; green, 0; blue, 0 }  ][line width=0.08]  [draw opacity=0] (8.93,-4.29) -- (0,0) -- (8.93,4.29) -- cycle    ;
%Curve Lines [id:da7787158109149253] 
\draw [color={rgb, 255:red, 0; green, 0; blue, 0 }  ,draw opacity=1 ] [dash pattern={on 4.5pt off 4.5pt}]  (139,8262.5) .. controls (133.77,8248.45) and (146.47,8215.52) .. (147.61,8260.01) ;
\draw [shift={(147.67,8262.83)}, rotate = 269.22] [fill={rgb, 255:red, 0; green, 0; blue, 0 }  ,fill opacity=1 ][line width=0.08]  [draw opacity=0] (8.93,-4.29) -- (0,0) -- (8.93,4.29) -- cycle    ;
%Curve Lines [id:da03637827887160494] 
\draw [color={rgb, 255:red, 0; green, 0; blue, 0 }  ,draw opacity=1 ] [dash pattern={on 4.5pt off 4.5pt}]  (147.67,8337.17) .. controls (150.88,8365.47) and (135.48,8364.6) .. (139.19,8340.54) ;
\draw [shift={(139.67,8337.83)}, rotate = 461.31] [fill={rgb, 255:red, 0; green, 0; blue, 0 }  ,fill opacity=1 ][line width=0.08]  [draw opacity=0] (8.93,-4.29) -- (0,0) -- (8.93,4.29) -- cycle    ;
%Curve Lines [id:da16956494999068394] 
\draw  [dash pattern={on 4.5pt off 4.5pt}]  (325.2,8269.94) .. controls (311.57,8277.2) and (299.96,8282.4) .. (281.39,8271.67) ;
\draw [shift={(279.04,8270.25)}, rotate = 392.23] [fill={rgb, 255:red, 0; green, 0; blue, 0 }  ][line width=0.08]  [draw opacity=0] (8.93,-4.29) -- (0,0) -- (8.93,4.29) -- cycle    ;
%Curve Lines [id:da7527323824097971] 
\draw    (322.54,8168.28) .. controls (309.54,8176.49) and (300.11,8180.81) .. (281.72,8170.01) ;
\draw [shift={(279.37,8168.58)}, rotate = 392.23] [fill={rgb, 255:red, 0; green, 0; blue, 0 }  ][line width=0.08]  [draw opacity=0] (8.93,-4.29) -- (0,0) -- (8.93,4.29) -- cycle    ;
%Straight Lines [id:da10404812613178782] 
\draw  [dash pattern={on 4.5pt off 4.5pt}]  (152.7,8118.58) -- (259.87,8118.28) ;
\draw [shift={(262.87,8118.28)}, rotate = 539.8399999999999] [fill={rgb, 255:red, 0; green, 0; blue, 0 }  ][line width=0.08]  [draw opacity=0] (8.93,-4.29) -- (0,0) -- (8.93,4.29) -- cycle    ;
%Straight Lines [id:da667970233163045] 
\draw    (178.33,8162.17) -- (262.91,8126.01) ;
\draw [shift={(265.67,8124.83)}, rotate = 516.85] [fill={rgb, 255:red, 0; green, 0; blue, 0 }  ][line width=0.08]  [draw opacity=0] (8.93,-4.29) -- (0,0) -- (8.93,4.29) -- cycle    ;
%Straight Lines [id:da34356330629343423] 
\draw    (181.7,8168.58) -- (259.87,8168.61) ;
\draw [shift={(262.87,8168.61)}, rotate = 180.02] [fill={rgb, 255:red, 0; green, 0; blue, 0 }  ][line width=0.08]  [draw opacity=0] (8.93,-4.29) -- (0,0) -- (8.93,4.29) -- cycle    ;
%Straight Lines [id:da7919336810292574] 
\draw    (181,8172.83) -- (327.1,8261.28) ;
\draw [shift={(329.67,8262.83)}, rotate = 211.19] [fill={rgb, 255:red, 0; green, 0; blue, 0 }  ][line width=0.08]  [draw opacity=0] (8.93,-4.29) -- (0,0) -- (8.93,4.29) -- cycle    ;
%Straight Lines [id:da5230449552943881] 
\draw    (119,8174.17) -- (262.45,8262.59) ;
\draw [shift={(265,8264.17)}, rotate = 211.65] [fill={rgb, 255:red, 0; green, 0; blue, 0 }  ][line width=0.08]  [draw opacity=0] (8.93,-4.29) -- (0,0) -- (8.93,4.29) -- cycle    ;
%Straight Lines [id:da8887115329722699] 
\draw    (262.54,8270.28) -- (154.37,8270.25) ;
\draw [shift={(151.37,8270.25)}, rotate = 360.01] [fill={rgb, 255:red, 0; green, 0; blue, 0 }  ][line width=0.08]  [draw opacity=0] (8.93,-4.29) -- (0,0) -- (8.93,4.29) -- cycle    ;
%Straight Lines [id:da2900746074278133] 
\draw    (329,8276.83) -- (154.24,8329.71) ;
\draw [shift={(151.37,8330.58)}, rotate = 343.16999999999996] [fill={rgb, 255:red, 0; green, 0; blue, 0 }  ][line width=0.08]  [draw opacity=0] (8.93,-4.29) -- (0,0) -- (8.93,4.29) -- cycle    ;
%Straight Lines [id:da15941268255094299] 
\draw  [dash pattern={on 4.5pt off 4.5pt}]  (265,8275.5) -- (152.42,8323.65) ;
\draw [shift={(149.67,8324.83)}, rotate = 336.84000000000003] [fill={rgb, 255:red, 0; green, 0; blue, 0 }  ][line width=0.08]  [draw opacity=0] (8.93,-4.29) -- (0,0) -- (8.93,4.29) -- cycle    ;
%Curve Lines [id:da46094881888788763] 
\draw  [dash pattern={on 4.5pt off 4.5pt}]  (150.33,8266.17) .. controls (166.22,8259.26) and (223.38,8255.68) .. (257.74,8266.63) ;
\draw [shift={(260.33,8267.5)}, rotate = 199.44] [fill={rgb, 255:red, 0; green, 0; blue, 0 }  ][line width=0.08]  [draw opacity=0] (8.93,-4.29) -- (0,0) -- (8.93,4.29) -- cycle    ;
%Curve Lines [id:da8303721476833326] 
\draw  [dash pattern={on 4.5pt off 4.5pt}]  (330.27,8280.12) .. controls (277.92,8311.5) and (223.45,8327.37) .. (147.67,8337.17) ;
\draw [shift={(333.47,8278.18)}, rotate = 148.51] [fill={rgb, 255:red, 0; green, 0; blue, 0 }  ][line width=0.08]  [draw opacity=0] (8.93,-4.29) -- (0,0) -- (8.93,4.29) -- cycle    ;
%Curve Lines [id:da2689992069764009] 
\draw    (118.52,8177.43) .. controls (115.59,8198.46) and (115.43,8226.9) .. (136.33,8266.17) ;
\draw [shift={(119,8174.17)}, rotate = 98.88] [fill={rgb, 255:red, 0; green, 0; blue, 0 }  ][line width=0.08]  [draw opacity=0] (8.93,-4.29) -- (0,0) -- (8.93,4.29) -- cycle    ;
%Shape: Circle [id:dp4151777001229171] 
\draw  [color={rgb, 255:red, 155; green, 155; blue, 155 }  ,draw opacity=0.6 ] (82.67,8149.67) .. controls (82.67,8114.6) and (111.1,8086.17) .. (146.17,8086.17) .. controls (181.24,8086.17) and (209.67,8114.6) .. (209.67,8149.67) .. controls (209.67,8184.74) and (181.24,8213.17) .. (146.17,8213.17) .. controls (111.1,8213.17) and (82.67,8184.74) .. (82.67,8149.67) -- cycle ;
%Shape: Circle [id:dp41993966935304083] 
\draw  [color={rgb, 255:red, 155; green, 155; blue, 155 }  ,draw opacity=0.6 ] (237.33,8140.33) .. controls (237.33,8105.26) and (265.76,8076.83) .. (300.83,8076.83) .. controls (335.9,8076.83) and (364.33,8105.26) .. (364.33,8140.33) .. controls (364.33,8175.4) and (335.9,8203.83) .. (300.83,8203.83) .. controls (265.76,8203.83) and (237.33,8175.4) .. (237.33,8140.33) -- cycle ;
%Shape: Circle [id:dp956448311484482] 
\draw  [color={rgb, 255:red, 155; green, 155; blue, 155 }  ,draw opacity=0.6 ] (79.67,8298) .. controls (79.67,8262.93) and (108.1,8234.5) .. (143.17,8234.5) .. controls (178.24,8234.5) and (206.67,8262.93) .. (206.67,8298) .. controls (206.67,8333.07) and (178.24,8361.5) .. (143.17,8361.5) .. controls (108.1,8361.5) and (79.67,8333.07) .. (79.67,8298) -- cycle ;
%Shape: Circle [id:dp6882566434382558] 
\draw  [color={rgb, 255:red, 155; green, 155; blue, 155 }  ,draw opacity=0.6 ] (237.33,8271) .. controls (237.33,8235.93) and (265.76,8207.5) .. (300.83,8207.5) .. controls (335.9,8207.5) and (364.33,8235.93) .. (364.33,8271) .. controls (364.33,8306.07) and (335.9,8334.5) .. (300.83,8334.5) .. controls (265.76,8334.5) and (237.33,8306.07) .. (237.33,8271) -- cycle ;
%Curve Lines [id:da1778027310138357] 
\draw  [dash pattern={on 4.5pt off 4.5pt}]  (185.33,8169.75) .. controls (263.88,8195.23) and (291.64,8216.85) .. (333.44,8261.68) ;
\draw [shift={(181.7,8168.58)}, rotate = 17.73] [fill={rgb, 255:red, 0; green, 0; blue, 0 }  ][line width=0.08]  [draw opacity=0] (8.93,-4.29) -- (0,0) -- (8.93,4.29) -- cycle    ;
%Curve Lines [id:da3822636587604338] 
\draw    (340.1,8171.73) .. controls (353.49,8216.28) and (349.91,8231.6) .. (333.44,8261.68) ;
\draw [shift={(339.04,8168.25)}, rotate = 72.8] [fill={rgb, 255:red, 0; green, 0; blue, 0 }  ][line width=0.08]  [draw opacity=0] (8.93,-4.29) -- (0,0) -- (8.93,4.29) -- cycle    ;
%Shape: Circle [id:dp23239117448548496] 
\draw  [color={rgb, 255:red, 0; green, 0; blue, 0 }  ,draw opacity=1 ][fill={rgb, 255:red, 74; green, 74; blue, 74 }  ,fill opacity=1 ] (113.33,7924.83) .. controls (113.33,7922.26) and (115.42,7920.17) .. (118,7920.17) .. controls (120.58,7920.17) and (122.67,7922.26) .. (122.67,7924.83) .. controls (122.67,7927.41) and (120.58,7929.5) .. (118,7929.5) .. controls (115.42,7929.5) and (113.33,7927.41) .. (113.33,7924.83) -- cycle ;
%Shape: Circle [id:dp55910454201637] 
\draw  [color={rgb, 255:red, 0; green, 0; blue, 0 }  ,draw opacity=1 ][fill={rgb, 255:red, 74; green, 74; blue, 74 }  ,fill opacity=1 ] (95.33,7954.83) .. controls (95.33,7952.26) and (97.42,7950.17) .. (100,7950.17) .. controls (102.58,7950.17) and (104.67,7952.26) .. (104.67,7954.83) .. controls (104.67,7957.41) and (102.58,7959.5) .. (100,7959.5) .. controls (97.42,7959.5) and (95.33,7957.41) .. (95.33,7954.83) -- cycle ;
%Shape: Circle [id:dp6210428674830062] 
\draw  [color={rgb, 255:red, 0; green, 0; blue, 0 }  ,draw opacity=1 ][fill={rgb, 255:red, 74; green, 74; blue, 74 }  ,fill opacity=1 ] (133.33,7954.83) .. controls (133.33,7952.26) and (135.42,7950.17) .. (138,7950.17) .. controls (140.58,7950.17) and (142.67,7952.26) .. (142.67,7954.83) .. controls (142.67,7957.41) and (140.58,7959.5) .. (138,7959.5) .. controls (135.42,7959.5) and (133.33,7957.41) .. (133.33,7954.83) -- cycle ;
%Curve Lines [id:da6967902945998676] 
\draw [color={rgb, 255:red, 0; green, 0; blue, 0 }  ,draw opacity=1 ] [dash pattern={on 4.5pt off 4.5pt}]  (115.33,7922.17) .. controls (107.1,7910.13) and (124.14,7902.31) .. (121.86,7918.71) ;
\draw [shift={(121.33,7921.5)}, rotate = 283.13] [fill={rgb, 255:red, 0; green, 0; blue, 0 }  ,fill opacity=1 ][line width=0.08]  [draw opacity=0] (5.36,-2.57) -- (0,0) -- (5.36,2.57) -- cycle    ;
%Curve Lines [id:da641812626653862] 
\draw [color={rgb, 255:red, 0; green, 0; blue, 0 }  ,draw opacity=1 ]   (100,7959.5) .. controls (94.96,7973.99) and (84.56,7970.02) .. (93.56,7957.17) ;
\draw [shift={(95.33,7954.83)}, rotate = 489.29] [fill={rgb, 255:red, 0; green, 0; blue, 0 }  ,fill opacity=1 ][line width=0.08]  [draw opacity=0] (5.36,-2.57) -- (0,0) -- (5.36,2.57) -- cycle    ;
%Straight Lines [id:da8204701860219279] 
\draw    (100,7950.17) -- (113.36,7929.68) ;
\draw [shift={(115,7927.17)}, rotate = 483.11] [fill={rgb, 255:red, 0; green, 0; blue, 0 }  ][line width=0.08]  [draw opacity=0] (5.36,-2.57) -- (0,0) -- (5.36,2.57) -- cycle    ;
%Straight Lines [id:da4993124052803892] 
\draw  [dash pattern={on 4.5pt off 4.5pt}]  (104.67,7954.83) -- (130.33,7954.83) ;
\draw [shift={(133.33,7954.83)}, rotate = 180] [fill={rgb, 255:red, 0; green, 0; blue, 0 }  ][line width=0.08]  [draw opacity=0] (5.36,-2.57) -- (0,0) -- (5.36,2.57) -- cycle    ;
%Curve Lines [id:da7342251642528286] 
\draw    (125.24,7926.66) .. controls (135.98,7934.59) and (138,7940.87) .. (138,7950.17) ;
\draw [shift={(122.67,7924.83)}, rotate = 34.2] [fill={rgb, 255:red, 0; green, 0; blue, 0 }  ][line width=0.08]  [draw opacity=0] (5.36,-2.57) -- (0,0) -- (5.36,2.57) -- cycle    ;
%Curve Lines [id:da19155406447876477] 
\draw  [dash pattern={on 4.5pt off 4.5pt}]  (122.67,7924.83) .. controls (122.98,7933.67) and (125.02,7942.22) .. (135.49,7948.73) ;
\draw [shift={(138,7950.17)}, rotate = 207.76] [fill={rgb, 255:red, 0; green, 0; blue, 0 }  ][line width=0.08]  [draw opacity=0] (5.36,-2.57) -- (0,0) -- (5.36,2.57) -- cycle    ;
%Shape: Circle [id:dp8003093330750504] 
\draw  [color={rgb, 255:red, 0; green, 0; blue, 0 }  ,draw opacity=1 ][fill={rgb, 255:red, 74; green, 74; blue, 74 }  ,fill opacity=1 ] (175,7925.17) .. controls (175,7922.59) and (177.09,7920.5) .. (179.67,7920.5) .. controls (182.24,7920.5) and (184.33,7922.59) .. (184.33,7925.17) .. controls (184.33,7927.74) and (182.24,7929.83) .. (179.67,7929.83) .. controls (177.09,7929.83) and (175,7927.74) .. (175,7925.17) -- cycle ;
%Shape: Circle [id:dp7909782191660482] 
\draw  [color={rgb, 255:red, 0; green, 0; blue, 0 }  ,draw opacity=1 ][fill={rgb, 255:red, 74; green, 74; blue, 74 }  ,fill opacity=1 ] (213,7925.17) .. controls (213,7922.59) and (215.09,7920.5) .. (217.67,7920.5) .. controls (220.24,7920.5) and (222.33,7922.59) .. (222.33,7925.17) .. controls (222.33,7927.74) and (220.24,7929.83) .. (217.67,7929.83) .. controls (215.09,7929.83) and (213,7927.74) .. (213,7925.17) -- cycle ;
%Shape: Circle [id:dp049222116368238256] 
\draw  [color={rgb, 255:red, 0; green, 0; blue, 0 }  ,draw opacity=1 ][fill={rgb, 255:red, 74; green, 74; blue, 74 }  ,fill opacity=1 ] (175.33,7954.83) .. controls (175.33,7952.26) and (177.42,7950.17) .. (180,7950.17) .. controls (182.58,7950.17) and (184.67,7952.26) .. (184.67,7954.83) .. controls (184.67,7957.41) and (182.58,7959.5) .. (180,7959.5) .. controls (177.42,7959.5) and (175.33,7957.41) .. (175.33,7954.83) -- cycle ;
%Shape: Circle [id:dp01310000091092789] 
\draw  [color={rgb, 255:red, 0; green, 0; blue, 0 }  ,draw opacity=1 ][fill={rgb, 255:red, 74; green, 74; blue, 74 }  ,fill opacity=1 ] (213.33,7954.83) .. controls (213.33,7952.26) and (215.42,7950.17) .. (218,7950.17) .. controls (220.58,7950.17) and (222.67,7952.26) .. (222.67,7954.83) .. controls (222.67,7957.41) and (220.58,7959.5) .. (218,7959.5) .. controls (215.42,7959.5) and (213.33,7957.41) .. (213.33,7954.83) -- cycle ;
%Straight Lines [id:da17204143283328377] 
\draw  [dash pattern={on 4.5pt off 4.5pt}]  (187.33,7925.17) -- (213,7925.17) ;
\draw [shift={(184.33,7925.17)}, rotate = 0] [fill={rgb, 255:red, 0; green, 0; blue, 0 }  ][line width=0.08]  [draw opacity=0] (5.36,-2.57) -- (0,0) -- (5.36,2.57) -- cycle    ;
%Curve Lines [id:da22875953739189647] 
\draw    (213.33,7954.83) .. controls (203.67,7961.69) and (198.38,7963.3) .. (187.11,7956.4) ;
\draw [shift={(184.67,7954.83)}, rotate = 393.69] [fill={rgb, 255:red, 0; green, 0; blue, 0 }  ][line width=0.08]  [draw opacity=0] (5.36,-2.57) -- (0,0) -- (5.36,2.57) -- cycle    ;
%Curve Lines [id:da4795965460430833] 
\draw    (184.67,7954.83) .. controls (191.72,7948.9) and (200.27,7948.24) .. (210.75,7953.45) ;
\draw [shift={(213.33,7954.83)}, rotate = 209.74] [fill={rgb, 255:red, 0; green, 0; blue, 0 }  ][line width=0.08]  [draw opacity=0] (5.36,-2.57) -- (0,0) -- (5.36,2.57) -- cycle    ;
%Straight Lines [id:da7137783723649642] 
\draw  [dash pattern={on 4.5pt off 4.5pt}]  (218,7950.17) -- (217.72,7932.83) ;
\draw [shift={(217.67,7929.83)}, rotate = 449.06] [fill={rgb, 255:red, 0; green, 0; blue, 0 }  ][line width=0.08]  [draw opacity=0] (5.36,-2.57) -- (0,0) -- (5.36,2.57) -- cycle    ;
%Curve Lines [id:da24970244782323014] 
\draw [color={rgb, 255:red, 0; green, 0; blue, 0 }  ,draw opacity=1 ]   (176.33,7921.5) .. controls (175.08,7902.07) and (182.66,7904.44) .. (183.02,7918.64) ;
\draw [shift={(183,7921.5)}, rotate = 272.29] [fill={rgb, 255:red, 0; green, 0; blue, 0 }  ,fill opacity=1 ][line width=0.08]  [draw opacity=0] (5.36,-2.57) -- (0,0) -- (5.36,2.57) -- cycle    ;
%Curve Lines [id:da07963128257803165] 
\draw [color={rgb, 255:red, 0; green, 0; blue, 0 }  ,draw opacity=1 ]   (214.33,7921.5) .. controls (212.42,7904.31) and (222.07,7900.5) .. (221.2,7918.77) ;
\draw [shift={(221,7921.5)}, rotate = 275.36] [fill={rgb, 255:red, 0; green, 0; blue, 0 }  ,fill opacity=1 ][line width=0.08]  [draw opacity=0] (5.36,-2.57) -- (0,0) -- (5.36,2.57) -- cycle    ;
%Curve Lines [id:da6039552321054946] 
\draw [color={rgb, 255:red, 0; green, 0; blue, 0 }  ,draw opacity=1 ] [dash pattern={on 4.5pt off 4.5pt}]  (221.67,7958.83) .. controls (223.55,7978.89) and (214.82,7976.56) .. (214.32,7961.81) ;
\draw [shift={(214.33,7958.83)}, rotate = 452.2] [fill={rgb, 255:red, 0; green, 0; blue, 0 }  ,fill opacity=1 ][line width=0.08]  [draw opacity=0] (5.36,-2.57) -- (0,0) -- (5.36,2.57) -- cycle    ;
%Shape: Circle [id:dp13974700574683685] 
\draw  [color={rgb, 255:red, 0; green, 0; blue, 0 }  ,draw opacity=1 ][fill={rgb, 255:red, 74; green, 74; blue, 74 }  ,fill opacity=1 ] (332,7918.83) .. controls (332,7916.26) and (334.09,7914.17) .. (336.67,7914.17) .. controls (339.24,7914.17) and (341.33,7916.26) .. (341.33,7918.83) .. controls (341.33,7921.41) and (339.24,7923.5) .. (336.67,7923.5) .. controls (334.09,7923.5) and (332,7921.41) .. (332,7918.83) -- cycle ;
%Shape: Circle [id:dp4915521776109293] 
\draw  [color={rgb, 255:red, 0; green, 0; blue, 0 }  ,draw opacity=1 ][fill={rgb, 255:red, 74; green, 74; blue, 74 }  ,fill opacity=1 ] (332.33,7956.5) .. controls (332.33,7953.92) and (334.42,7951.83) .. (337,7951.83) .. controls (339.58,7951.83) and (341.67,7953.92) .. (341.67,7956.5) .. controls (341.67,7959.08) and (339.58,7961.17) .. (337,7961.17) .. controls (334.42,7961.17) and (332.33,7959.08) .. (332.33,7956.5) -- cycle ;
%Curve Lines [id:da6702356289688982] 
\draw [color={rgb, 255:red, 0; green, 0; blue, 0 }  ,draw opacity=1 ] [dash pattern={on 4.5pt off 4.5pt}]  (333.33,7915.17) .. controls (331.42,7897.98) and (341.07,7894.16) .. (340.2,7912.44) ;
\draw [shift={(340,7915.17)}, rotate = 275.36] [fill={rgb, 255:red, 0; green, 0; blue, 0 }  ,fill opacity=1 ][line width=0.08]  [draw opacity=0] (5.36,-2.57) -- (0,0) -- (5.36,2.57) -- cycle    ;
%Curve Lines [id:da18642341810442598] 
\draw [color={rgb, 255:red, 0; green, 0; blue, 0 }  ,draw opacity=1 ] [dash pattern={on 4.5pt off 4.5pt}]  (340.67,7960.5) .. controls (342.55,7980.55) and (333.82,7978.22) .. (333.32,7963.47) ;
\draw [shift={(333.33,7960.5)}, rotate = 452.2] [fill={rgb, 255:red, 0; green, 0; blue, 0 }  ,fill opacity=1 ][line width=0.08]  [draw opacity=0] (5.36,-2.57) -- (0,0) -- (5.36,2.57) -- cycle    ;
%Curve Lines [id:da05469626662158378] 
\draw  [dash pattern={on 4.5pt off 4.5pt}]  (337,7951.83) .. controls (331.63,7946.98) and (332.35,7933.63) .. (335.33,7926.18) ;
\draw [shift={(336.67,7923.5)}, rotate = 482.28] [fill={rgb, 255:red, 0; green, 0; blue, 0 }  ][line width=0.08]  [draw opacity=0] (5.36,-2.57) -- (0,0) -- (5.36,2.57) -- cycle    ;
%Curve Lines [id:da5661588540937952] 
\draw  [dash pattern={on 4.5pt off 4.5pt}]  (336.67,7923.5) .. controls (343.05,7932.49) and (341.35,7942.74) .. (338.39,7949.21) ;
\draw [shift={(337,7951.83)}, rotate = 301.43] [fill={rgb, 255:red, 0; green, 0; blue, 0 }  ][line width=0.08]  [draw opacity=0] (5.36,-2.57) -- (0,0) -- (5.36,2.57) -- cycle    ;
%Shape: Circle [id:dp07650268425973317] 
\draw  [color={rgb, 255:red, 0; green, 0; blue, 0 }  ,draw opacity=1 ][fill={rgb, 255:red, 74; green, 74; blue, 74 }  ,fill opacity=1 ] (254,7948.83) .. controls (254,7946.26) and (256.09,7944.17) .. (258.67,7944.17) .. controls (261.24,7944.17) and (263.33,7946.26) .. (263.33,7948.83) .. controls (263.33,7951.41) and (261.24,7953.5) .. (258.67,7953.5) .. controls (256.09,7953.5) and (254,7951.41) .. (254,7948.83) -- cycle ;
%Shape: Circle [id:dp5090252750831064] 
\draw  [color={rgb, 255:red, 0; green, 0; blue, 0 }  ,draw opacity=1 ][fill={rgb, 255:red, 74; green, 74; blue, 74 }  ,fill opacity=1 ] (292,7948.83) .. controls (292,7946.26) and (294.09,7944.17) .. (296.67,7944.17) .. controls (299.24,7944.17) and (301.33,7946.26) .. (301.33,7948.83) .. controls (301.33,7951.41) and (299.24,7953.5) .. (296.67,7953.5) .. controls (294.09,7953.5) and (292,7951.41) .. (292,7948.83) -- cycle ;
%Curve Lines [id:da010689684841189928] 
\draw  [dash pattern={on 4.5pt off 4.5pt}]  (292,7948.83) .. controls (282.34,7955.69) and (277.05,7957.3) .. (265.77,7950.4) ;
\draw [shift={(263.33,7948.83)}, rotate = 393.69] [fill={rgb, 255:red, 0; green, 0; blue, 0 }  ][line width=0.08]  [draw opacity=0] (5.36,-2.57) -- (0,0) -- (5.36,2.57) -- cycle    ;
%Curve Lines [id:da3515090761062116] 
\draw  [dash pattern={on 4.5pt off 4.5pt}]  (263.33,7948.83) .. controls (270.39,7942.9) and (278.94,7942.24) .. (289.41,7947.45) ;
\draw [shift={(292,7948.83)}, rotate = 209.74] [fill={rgb, 255:red, 0; green, 0; blue, 0 }  ][line width=0.08]  [draw opacity=0] (5.36,-2.57) -- (0,0) -- (5.36,2.57) -- cycle    ;

% Text Node
\draw (113.45,8167.6) node  [font=\tiny,rotate=-1.27] [align=left] {$\displaystyle 1^{1}$};
% Text Node
\draw (173.45,8167.6) node  [font=\tiny,rotate=-1.27] [align=left] {$\displaystyle 1^{3}$};
% Text Node
\draw (144.45,8117.6) node  [font=\tiny,rotate=-1.27] [align=left] {$\displaystyle 1^{2}$};
% Text Node
\draw (271.12,8167.6) node  [font=\tiny,rotate=-1.27] [align=left] {$\displaystyle 2^{1}$};
% Text Node
\draw (271.12,8117.26) node  [font=\tiny,rotate=-1.27] [align=left] {$\displaystyle 2^{2}$};
% Text Node
\draw (330.79,8167.26) node  [font=\tiny,rotate=-1.27] [align=left] {$\displaystyle 2^{4}$};
% Text Node
\draw (330.79,8116.93) node  [font=\tiny,rotate=-1.27] [align=left] {$\displaystyle 2^{3}$};
% Text Node
\draw (143.12,8329.6) node  [font=\tiny,rotate=-1.27] [align=left] {$\displaystyle 4^{2}$};
% Text Node
\draw (143.12,8269.26) node  [font=\tiny,rotate=-1.27] [align=left] {$\displaystyle 4^{1}$};
% Text Node
\draw (270.79,8269.26) node  [font=\tiny,rotate=-1.27] [align=left] {$\displaystyle 3^{1}$};
% Text Node
\draw (333.45,8268.93) node  [font=\tiny,rotate=-1.27] [align=left] {$\displaystyle 3^{2}$};
% Text Node
\draw (143,8073.33) node   [align=left] {\begin{minipage}[lt]{10.88000000000001pt}\setlength\topsep{0pt}
$\displaystyle 1$
\end{minipage}};
% Text Node
\draw (299.67,8064.67) node   [align=left] {\begin{minipage}[lt]{10.88000000000001pt}\setlength\topsep{0pt}
$\displaystyle 2$
\end{minipage}};
% Text Node
\draw (329.67,8341) node   [align=left] {\begin{minipage}[lt]{10.88000000000001pt}\setlength\topsep{0pt}
$\displaystyle 3$
\end{minipage}};
% Text Node
\draw (199,8348.33) node   [align=left] {\begin{minipage}[lt]{10.88000000000001pt}\setlength\topsep{0pt}
$\displaystyle 4$
\end{minipage}};
% Text Node
\draw (233,8390.33) node   [align=left] {\begin{minipage}[lt]{15.64pt}\setlength\topsep{0pt}
(b)
\end{minipage}};
% Text Node
\draw (125.33,7989.5) node   [align=left] {\begin{minipage}[lt]{24.48pt}\setlength\topsep{0pt}
$\displaystyle G_{1}$
\end{minipage}};
% Text Node
\draw (200.33,7989) node   [align=left] {\begin{minipage}[lt]{25.386666666666695pt}\setlength\topsep{0pt}
 $\displaystyle G_{2}$
\end{minipage}};
% Text Node
\draw (279,7989.33) node   [align=left] {\begin{minipage}[lt]{25.386666666666695pt}\setlength\topsep{0pt}
 $\displaystyle G_{3}$
\end{minipage}};
% Text Node
\draw (342,7989.33) node   [align=left] {\begin{minipage}[lt]{25.386666666666695pt}\setlength\topsep{0pt}
 $\displaystyle G_{4}$
\end{minipage}};
% Text Node
\draw (234.42,8025.67) node   [align=left] {\begin{minipage}[lt]{17.793333333333308pt}\setlength\topsep{0pt}
(a)
\end{minipage}};

\end{tikzpicture}

%%%%%%%%%%%%%%%%%%%%%%%%%%%%%%%%
\caption{(a) Node graphs $G_1$, $G_2$, $G_3$, and $G_4$; (b) global graph $\mathcal{G}$.}
\label{globgraph}
\end{figure}
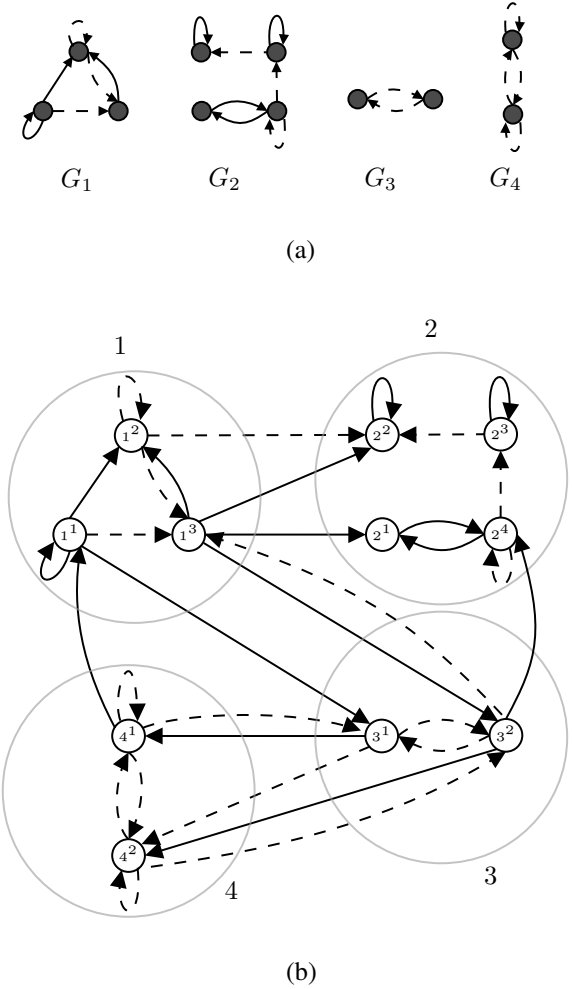

\begin{ex}
Consider an LTI network $\mathcal{N}$, including subsystems $1$, $2$, $3$, and $4$, which are depicted in Fig.~\ref{globgraph}. The node graphs $G_1$, $G_2$, $G_3$, and $G_4$ of size $l_1=3$, $l_2=4$, $l_3=2$, and $l_4=2$, respectively, are shown in Fig.~\ref{globgraph}(a). The global graph $\mathcal{G}$ of size $N=11$ is also illustrated in Fig.~\ref{globgraph}(b). For instance, the corresponding pattern matrices $\mathcal{A}_{21}$ and $\mathcal{A}_{11}$, where $G_1=G(\mathcal{A}_{11})$, are:% 

\begin{equation*}\mathcal{A}_{21}=\begin{bmatrix} \Scale[1]0  & \Scale[1] 0 & \Scale[1] *\\ \Scale[1] 0 & \Scale[1] ? & \Scale[1] * \\ \Scale[1] 0 & \Scale[1] 0 & \Scale[1] 0\\ \Scale[1] 0 & \Scale[1] 0 & \Scale[1] 0
\end{bmatrix},\:\:
\mathcal{A}_{11}=\begin{bmatrix} \Scale[1] * & \Scale[1] 0 & \Scale[1] 0\\ \Scale[1] * & \Scale[1] ? & \Scale[1] * \\ \Scale[1] ? & \Scale[1] ? & \Scale[1] 0
\end{bmatrix}.\end{equation*} 
\end{ex}

\subsection{LTI Networks}

 Consider an LTI network $\mathcal{N}$, including  LTI subsystem $i$ with the system matrix $A_{ii}\in\mathbb{R}^{l_i\times l_i}$, $i=1,2,\ldots, n$.  If all the subsystems are one-dimensional, that is,  $l_i=1$, for all $1\leq i\leq n$, we have a \emph{network of one-dimensional subsystems} (\emph{N1DS}).%Let $N=\sum_{i=1}^n{l_i}$. 

The global dynamics of an LTI network is described by:
\begin{align}
\vspace{-.5em}
\dot{x}=Ax+Bu, 
\label{e1}
\end{align}  
where $A\in \mathcal{P}(\mathcal{A})$ is described in~(\ref{sm}). Let $x_i=\begin{bmatrix}x_{i}^1 & \ldots & x_{i}^{l_i}\end{bmatrix}^T\in \mathbb{R}^{l_i}$ be the state vector of subsystem $i$. Then,  $x=\begin{bmatrix}x_1^T & \ldots & x_n^T\end{bmatrix}^T\in \mathbb{R}^N$ is the aggregated state vector of the global system. Now, assume that we choose subsystem $j_k$, $1\leq  k\leq  m$, and inject independent control signals $u_{k}^i$, $1\leq i\leq l_{j_k}$, to its any state $x_{j_k}^i$. Let the vertices of the node graph $G_{j_k}$ be indexed by $j_k^1,\ldots, j_k^{l_{j_k}}$. Now, let $M=\Sigma_{k=j_1}^{j_m} l_k $. Then, $u=\begin{bmatrix}u_1^T &\ldots u_m^T \end{bmatrix}^T\in\mathbb{R}^{M}$, where $u_k=\begin{bmatrix}u_k^1 &\ldots u_k^{l_{j_k}} \end{bmatrix}^T\in\mathbb{R}^{l_{j_k}}$, $k=1,2,\ldots, m$, is the vector of control signals injected into vertices of node graph $G_{j_k}$. 
Moreover, the input matrix  
$B=\begin{bmatrix}B_1 & \ldots & B_m\end{bmatrix}\in \mathbb{R}^{N\times M}$ is a binary matrix, where for $k=1,2,\ldots, m$, $B_k=\begin{bmatrix}e_{j_k^1}, \ldots, e_{j_k^{l_{j_k}}}\end{bmatrix}\in \mathbb{R}^{N\times l_{j_k}}$. In this case, the subsystems  $j_k$, $k=1,2,\ldots,m$, are referred to as 
\textit{control subsystems}, and  we define $V_C=\{j_1,\ldots,j_m\}$.

\begin{ex}
Consider an LTI network $\mathcal{N}$ with dynamics~(\ref{e1}), whose global graph $\mathcal{G}$ is depicted in Fig.~\ref{globgraph}(b). Let $x=\begin{bmatrix}x_1^T & x_2^T & x_3^T & x_4^T\end{bmatrix}^T$, where for instance, $x_1=\begin{bmatrix}x_1^1 & x_1^2 & x_1^3\end{bmatrix}^T$. Assume that we choose subsystem 1 as the control subsystem, i.e. $V_C=\{1\}$. Then, we have $u=u_1\in \mathbb{R}^3$, where  $u_1=\begin{bmatrix}u_1^1 & u_1^2 & u_1^3 \end{bmatrix}^T$. Moreover, $B=\begin{bmatrix}e_1 & e_2 & e_3\end{bmatrix}\in\mathbb{R}^{11\times 3}$, where $e_j$ is the $j$-th column of the identity matrix $I_{11}$. 
\end{ex}

{An LTI system with 
the pair $(A,B)$ in dynamics~(\ref{e1})
is controllable 
if with a suitable choice of the inputs, the states of all nodes can be driven from any initial state to any final state within a finite time. }
{As} controllability is preserved under the equivalent transformations,
when LTI system~(\ref{e1}) is uncontrollable, 
there exists a nonsingular matrix $T\in \mathbb{R}^{N\times N}$ such that {for some $0\leq q<N$,}
\vspace{-.5em}
{\begin{equation}
\Scale[1]{T^{-1}AT=\begin{bmatrix}\hat{A}_{11} & \hat{A}_{12}\\ 0 & \hat{A}_{22}\end{bmatrix}, \;
T^{-1}B=\begin{bmatrix}\hat{B}_1\\0\end{bmatrix}},
\label{decom}
\end{equation}}where $(\hat{A}_{11},\hat{B}_{1})$ is controllable, with $\hat{A}_{11}\in\mathbb{R}^{q\times q}$ and $\hat{B}_{1}\in\mathbb{R}^{q\times M}$. 

{When $\lambda\in \Lambda(A)$ and $\lambda\notin \Lambda(\hat{A}_{22})$, then it is referred to as a \textit{controllable eigenvalue} of the LTI system~(\ref{e1}). On the other hand, if  $\lambda\notin \Lambda(\hat{A}_{11})$, then it is called an \textit{uncontrollable eigenvalue}, which cannot be influenced by the system's input. {The} PBH test is used for checking the controllability of eigenvalues.}

{\begin{pro}[\cite{sontag2013mathematical}]
For the LTI model (\ref{e1}), the eigenvalue $\lambda$ of $A$ is controllable
if and only if for all nonzero $w$ for which $w^{T}A=\lambda w^T,$ also $w^{T}B\neq0$ holds.
\end{pro}}

\subsection{$\Delta$-Specified Pattern System Matrices}
Consider the block matrix $A$ in (\ref{sm}). In order to characterize the submatrices $A_{ii}$'s,  $i=1,2,\ldots,n$, in the block matrix $A$, it suffices to define a sequence $\mathcal{L}=(l_1,l_2,\ldots, l_n)$, including the dimension  of the  diagonal blocks. Let $\mathcal{L}^{k}=\sum_{j=1}^k\mathcal{L}(j)$. Then, given the sequence $\mathcal{L}$, one can simply write $A_{ii}=A(\mathcal{L}^{i-1}+1:\mathcal{L}^{i},\mathcal{L}^{i-1}+1:\mathcal{L}^{i})$, for $i=1,2,\ldots, n$.

Now, let $\Delta\subseteq \mathbb{C}$ be an arbitrary nonempty set, which can be either discrete (finite or infinite) or continuous. Let $A\in \mathbb{R}^{N\times N}$, defined in (\ref{sm}), be a system matrix associated with a network $\mathcal{N}$, including some subsystem $i$, $i=1,2,\ldots, n$, with system matrix $A_{ii}\in \mathbb{R}^{l_i\times l_i}$.  Let $A\in \mathcal{P}(\mathcal{A})$,  
where $\mathcal{A}\in\{0,*,?\}^{N\times N}$ 
be a pattern matrix. 
Assume that in addition to the zero/nonzero/arbitrary pattern of $A$, for some $1\leq i\leq n$, we may have some further information about subsystem $i$, in terms of the intersection of $\Delta$ and spectrum of ${A}_{ii}$. For instance, for some $1\leq i\leq n$, we may know that $\Lambda({A}_{ii})\cap \Delta=\emptyset$.
For a matrix $A\in \mathbb{R}^{N\times N}$,  a sequence $\mathcal{L}=(l_1,\ldots, l_n)$, with $\sum_{i=1}^n{l_i}=N$, is given, which can determine its diagonal blocks $A_{ii}$'s, $i=1,2,\ldots,n$. Now,  let us define a vector 
$f^{\mathcal{L}}_{\Delta}\in \{0,*,?\}^n$, representing the extra information about $\Lambda(A_{ii})\cap \Delta$, i.e., the intersection of $\Delta$ and spectrum of ${A}_{ii}$'s, $i=1,2,\ldots, n$. We refer to $f^{\mathcal{L}}_{\Delta}$ as a $\Delta$\emph{-characteristic vector}. If $\Delta$ is not a singleton (i.e., $|\Delta|>1$), we let $f^{\mathcal{L}}_{\Delta}\in \{*,?\}^n$. 

Now, given a  pattern matrix $\mathcal{A}\in\{0,*,?\}^{N\times N}$ and a $\Delta$-characteristic vector $f^{\mathcal{L}}_{\Delta}$, a \emph{$\Delta$-specified pattern class} of $f^{\mathcal{L}}_{\Delta}$ and $\mathcal{A}$, denoted by $\mathcal{S}(f^{\mathcal{L}}_{\Delta}, \mathcal{A})$, is defined as:
\begin{equation}
\begin{aligned}
 & \mathcal{S}(f^{\mathcal{L}}_{\Delta}, \mathcal{A})=\{ A\in \mathcal{P}(\mathcal{A})\mid \Lambda(A_{kk})\cap \Delta=\emptyset \:\: \text{if}\ f^{\mathcal{L}}_{\Delta}(k)=*, \\& \text{and}\  {A}_{kk}=\mu I_{l_k} \:\: \text{if for}\ \Delta=\{\mu\},\:\: f^{\mathcal{L}}_{\Delta}(k)=0\}.  
\end{aligned}
\label{sp}
\end{equation}

Accordingly, $\mathcal{S}(f^{\mathcal{L}}_{\Delta}, \mathcal{A})$ is the set of all matrices $A\in \mathbb{R}^{N\times N}$ which have the same pattern $\mathcal{A}$; moreover,  $\Delta $ and the spectrum of $A_{kk}$ should have no member in common if $f^{\mathcal{L}}_{\Delta}(k)=*$. Furthermore, $f^{\mathcal{L}}_{\Delta}(k)=0$ implies that for some $\mu\in \mathbb{R}$, $\Delta=\{\mu\}$ and ${A}_{kk}=\mu I_{l_k}$. In other words, if $f^{\mathcal{L}}_{\Delta}(k)=0$, $\Delta$ should be a singleton with only member $\mu$, and $A_{kk}$ is a diagonal matrix whose all diagonal entries are equal to $\mu$. Finally,  $f^{\mathcal{L}}_{\Delta}(k)=?$ implies that  $\Lambda(A_{kk})\cap \Delta$ can be empty or nonempty, or if we have $\Lambda(A_{kk})\cap \Delta\neq \emptyset$, then either $|\Delta|>1$, or for $\Delta=\{\mu\}$, $A_{kk}$ should not  necessarily equal $\mu I_{l_k}$.

Now, in the reverse direction, assume that  a set of system matrices $\mathcal{S}^*$ along with a nonempty set $\Delta\subseteq \mathbb{C}$ are given, and we aim to define a corresponding pattern matrix $\mathcal{A}$ and a $\Delta$-characteristic vector $f^{\mathcal{L}}_{\Delta}$, where $\mathcal{S}^*\subseteq \mathcal{P}(\mathcal{A})$ and $\mathcal{S}^*\subseteq \mathcal{S}(f^{\mathcal{L}}_{\Delta}, \mathcal{A})$. One can define $\mathcal{A}$ and $f^{\mathcal{L}}_{\Delta}$ as follows. For all $1\leq i,j\leq N$, one has:
\begin{equation}
   \mathcal{A}(i,j)=
    \begin{cases}
      *, & \text{if}\ A(i,j)\neq 0, \:\:  \forall A\in \mathcal{S}^*,  
      \\
      0, & \text{if}\ A(i,j)= 0, \:\:  \forall A\in \mathcal{S}^*,\\
      ?, & \text{otherwise}.
    \end{cases}
    \label{A}
\end{equation}
Moreover, for $k=1,2,\ldots, n$, one can define:
\begin{equation}
   f^{\mathcal{L}}_{\Delta}(k)=
    \begin{cases}
      *, & \text{if}\ \Lambda(A_{kk})\cap \Delta=\emptyset, \:\:  \forall A\in \mathcal{S}^*,  
      \\
      0, & \text{if}\ \Delta=\{\mu\}\:\: \&\:\: {A}_{kk}=\mu I_{l_k}, \:\:  \forall A\in \mathcal{S}^*,\\
      ?, & \text{otherwise}.
    \end{cases}
\label{f}
\end{equation}
\begin{ex}
Let $\Delta=\{y\in\mathbb{C} \mid \Re(y)\geq 0\}$. Now, assume that for some $J\subseteq\{1,2,\ldots,n\}$, we have $f^{\mathcal{L}}_{\Delta}(k)=*$ if and only if  $k\in J$. Moreover, since $|\Delta|>1$, we let $f^{\mathcal{L}}_{\Delta}(k)=?$ for all  $k\in \{1,2,\ldots,n\}\setminus J$. Therefore, $\mathcal{S}(f^{\mathcal{L}}_{\Delta}, \mathcal{A})$ is the set of all matrices $A\in \mathcal{P}(\mathcal{A})$, where for all $k\in J$,  $A_{kk}$ is stable, in the sense that all of its eigenvalues  are in the open left-half plane; moreover, for $k\in \{1,2,\ldots,n\}\setminus J$, either we do not have any information about the stability of $A_{kk}$ or we know that it is not stable.  Now, consider an N1DS, where $l_i=1$, for all $1\leq i\leq n$. In this case,  $\Lambda(A_{kk})\cap \Delta=\emptyset$ implies that $A_{kk}=A(k,k)<0$. Thus, $\mathcal{S}(f^{\mathcal{L}}_{\Delta}, \mathcal{A})$ is the set of all $A\in\mathcal{P}(\mathcal{A})$, where $A(k,k)<0$, for every $k\in J$. Therefore, besides the zero/nonzero/arbitrary pattern of the matrix $A$,  we have extra information about the sign of some of its diagonal entries.
\label{stab}
\end{ex}

\begin{ex}
In this example, consider an N1DS, and for some $a,b\in\mathbb{R}$, where $a<b$, assume that $\Delta=[a,b]=\{y\in \mathbb{R}\mid a\leq x\leq b\}$. Now, suppose that for some $J\subseteq\{1,2,\ldots,n\}$,  $f^{\mathcal{L}}_{\Delta}(k)=*$ if and only if  $k\in J$. Thus, one can conclude that if $k\in J$, then either $A(k,k)<a$ or $A(k,k)>b$. Otherwise, for $k\notin J$, either there is no information about the value of $A(k,k)$ or  we know that it can take values from the interval $[a,b]$. \label{int} 
\end{ex}

\begin{ex}
Let $\Delta=\{0\}$, and consider an N1DS. For a given pattern matrix $\mathcal{A}\in\{0,*,?\}^{n\times n}$, assume that  the set of system matrices $\mathcal{S}^*$  is the same as the set $\mathcal{P}(\mathcal{A})$. In this case, $\Lambda(A_{kk})\cap \Delta=\emptyset$ implies that $A(k,k)\neq 0$, for all $A\in \mathcal{S}^*$. Thus, (\ref{A}) leads to $\mathcal{A}(k,k)=*$. %
Moreover, from (\ref{f}), for $k=1,2,\ldots,n$, $f_{\Delta}^{\mathcal{L}}(k)=*$ if $\mathcal{A}(k,k)=*$, and $f_{\Delta}^{\mathcal{L}}(k)=0$ if $\mathcal{A}(k,k)=0$. Otherwise, if $\mathcal{A}(k,k)=?$, we have $f_{\Delta}^{\mathcal{L}}(k)=?$. \label{ex5}
\end{ex}

\begin{ex}
Let $\Delta=\mathbb{C}\setminus \{0\}$, and consider an N1DS. One can conclude that $\Lambda(A_{kk})\cap \Delta=\emptyset$ implies that $A(k,k)= 0$. Note that since $|\Delta|>1$, $f_{\Delta}^{\mathcal{L}}\in\{*,?\}^{n}$. Hence, considering $\mathcal{S}^*=\mathcal{P}(\mathcal{A})$, for a given pattern matrix $\mathcal{A}\in\{0,*,?\}^{n\times n}$ and for $k=1,2,\ldots, n$, we have  $f_{\Delta}^{\mathcal{L}}(k)=*$ if $\mathcal{A}(k,k)=0$, and $f_{\Delta}^{\mathcal{L}}(k)=?$ otherwise. \label{ex6}
\end{ex}

\begin{ex}
For an N1DS, let $\Delta=\mathbb{C}$. Since $|\Delta|>1$, one should have $f_{\Delta}^{\mathcal{L}}\in\{*,?\}^{n}$. Moreover, for all $1\leq k\leq n$, $\Lambda(A_{kk})\cap \Delta\neq\emptyset$. Now, for a given pattern matrix $\mathcal{A}\in\{0,*,?\}^{n\times n}$, let $\mathcal{S}^*=\mathcal{P}(\mathcal{A})$. Then, from (\ref{f}), one has   $f_{\Delta}^{\mathcal{L}}(k)=?$, for all $1\leq k\leq n$. \label{ex7}
\end{ex}

\begin{ex}
As the last example in this part, for an N1DS, assume that $\mathcal{A}\in\{0,*,?\}^{n\times n}$ is a pattern matrix, whose all diagonals are $?$. This pattern matrix is called an \emph{arbitrary-diagonal pattern matrix}. Now, consider three sets $\Delta_1=\{0\}$, $\Delta_2=\mathbb{C}\setminus\{0\}$, and $\Delta_3=\mathbb{C}$, and let  $\mathcal{S}^*=\mathcal{P}(\mathcal{A})$.  In this case, one can observe that for every $k$, $1\leq k\leq n$, since $\mathcal{A}(k,k)=?$, there is some $A\in\mathcal{P}(\mathcal{A})$ with $A_{kk}=A(k,k)=0$. Thus, for $i=1,3$, $\Lambda(A_{kk})\cap \Delta_i\neq \emptyset$. Moreover, there is some $A'\in\mathcal{P}(\mathcal{A})$ with $A'(k,k)\neq 0$. Therefore, for every $k$, $1\leq k\leq n$, there is  some matrix $A'\in \mathcal{S}^*$, where $A'_{kk}\neq 0$, and $\Lambda(A'_{kk})\cap \Delta_2\neq \emptyset$. Accordingly, considering (\ref{f}), one can conclude that $f_{\Delta_i}^{\mathcal{L}}(k)=?$, for all $k=1,2,\ldots,n$ and $i=1,2,3$.\label{ex8}
\end{ex}

%%%%%%%%%%%%%%%%%%%%%%%%%%%%%%%%%%%%%%%%%%%%%%%%%%%%%%
\subsection{Modal Strong Structural Controllability}
Let $\Delta\subseteq \mathbb{C}$ be a nonempty discrete or continuous set. For a given pattern matrix $\mathcal{A}\in \{0,*,?\}^{N\times N}$, let $A\in\mathcal{P}(\mathcal{A})$, defined in (\ref{sm}), be a system matrix associated with an LTI network $\mathcal{N}$, and let $\mathcal{L}=(l_1,\ldots, l_n)$ be a sequence including the dimensions of submatrices $A_{ii}$s, $i=1,2,\ldots,n$. Assume that a $\Delta$-characteristic vector $f^{\mathcal{L}}_{\Delta}$, defined in (\ref{f}), is given, which represents some information about $\Lambda(A_{ii})\cap \Delta$, $i=1,2,\ldots, n$. Let $\mathcal{S}(f^{\mathcal{L}}_{\Delta},\mathcal{A})$, defined in (\ref{sp}), be a $\Delta$-specified pattern class of $f^{\mathcal{L}}_{\Delta}$ and $\mathcal{A}$.

\begin{deff}
An LTI network $\mathcal{N}$ with dynamics (\ref{e1}) is \emph{(modal) strongly structurally controllable with respect to $\Delta$} if for every $\lambda\in\Delta$ and for all    $A\in \mathcal{S}(f^{\mathcal{L}}_{\Delta},\mathcal{A})$ that $\lambda\in \Lambda(A)$, $\lambda$ is a controllable eigenvalue. We refer to this network as $\Delta$-SSC.
\label{ss}
\end{deff}

Based on the classical definition of strong strutural controllability in the literature, a network   is called strongly structurally controllable if it is $\mathbb{C}$-SSC. Now, let $\Delta'=\mathbb{C}\setminus \Delta$ (note that for an undirected network, we define $\Delta'=\mathbb{R}\setminus \Delta$). Then, based on Definition \ref{ss}, an LTI network $\mathcal{N}$ with dynamics (\ref{e1}) is \emph{strongly structurally controllable} if it is both $\Delta$-SSC and $\Delta'$-SSC. In a more general case, for some $k\geq 1$ and the disjoint sets $\Delta_1, \ldots, \Delta_k$, where $\mathbb{C}=\bigcup_{i=1}^k \Delta_k$, an LTI network $\mathcal{N}$ is \emph{strongly structurally controllable} if for every $1\leq i\leq k$, it is  $\Delta_i$-SSC. Notice that for a larger $k$, a less conservative condition can be obtained for strong structural controllability of networks. 

Now, consider Example \ref{stab}, where for a given pattern matrix $\mathcal{A}\in \{0,*,?\}^{N\times N}$,  $A\in \mathcal{P}(\mathcal{A})$. Moreover,  for all $k\in J$, the subsystem $k$ is stable. An LTI network $\mathcal{N}$ is called \emph{strongly structurally stabilizable} if for $\Delta=\{y\in\mathbb{C}\mid \Re(y)\geq 0\}$, the network is $\Delta$-SSC.  

We note that if there is no constraint on the system matrix $A$ other than $A\in \mathcal{P}(\mathcal{A})$,  then from Theorem 14 of \cite{jia2020unifying}, an LTI network is strongly structurally controllable if and only if it is strongly structurally stabilizable; however, since in Example \ref{stab}, we consider a smaller set of system matrices $\mathcal{S}(f^{\mathcal{L}}_{\Delta}, \mathcal{A})$, where $\mathcal{S}(f^{\mathcal{L}}_{\Delta}, \mathcal{A})\subseteq \mathcal{P}(\mathcal{A})$, and we assume that for every $k\in J$, all eigenvalues of $A_{kk}$ are in the open left-half plane, the strong structural controllability and stabilizabity will not necessarily be equivalent.  

\subsection{Problem Formulation}

Consider an LTI network with dynamics (\ref{e1}), which includes $n$ LTI structured subsystems, and its system matrix $A\in \mathcal{P}(\mathcal{A})$ is described in (\ref{sm}). %We consider an arbitrary set $\Delta\subseteq \mathbb{C}$, and we assume that a sequence $\mathcal{L}$, including the dimension of any subsystem, along with a $\Delta$-characteristic vector  $f^{\mathcal{L}}_{\Delta}$ are given. 
Given a nonempty set $\Delta\subseteq \mathbb{C}$, our focus in this work is on the combinatorial characterizations
of modal strong structural controllability of an LTI network with respect to $\Delta$. %Then, we assume that $A\in \mathcal{S}(f^{\mathcal{L}}_{\Delta},\mathcal{A})$, where $\mathcal{S}(f^{\mathcal{L}}_{\Delta},\mathcal{A})$ is a $\Delta$-specified pattern class of $f^{\mathcal{L}}_{\Delta}$ and $\mathcal{A}$. 
Therefore, the main problem that we aim to investigate in this paper is the following. 

\emph{Problem 1:} Given an arbitrary $\Delta\subseteq \mathbb{C}$, find graph-theoretic conditions under which an LTI network with dynamics (\ref{e1}) is strongly structurally controllable with respect to $\Delta$.

Let us illustrate this problem on an applied example and discuss how the existing results on strong structural controllability are too conservative in this case. In this example, we let $\Delta=\{y\in\mathbb{C}| \Re(y)\geq 0\}$.  
The LTI network, which can be, for example,  a network of  robots or a  platoon of heterogeneous vehicles moving along a ring road,  includes stable subsystems. Due to the physical constraints, it is known that all subsystems are stable and have no associated eigenvalues in the closed right-half plane.

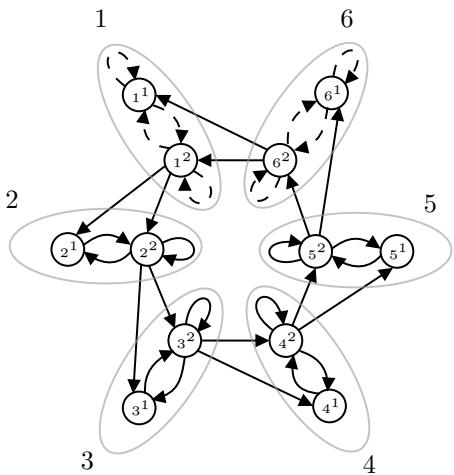
\begin{figure}[b]
\centering

\tikzset{every picture/.style={line width=0.75pt}} %set default line width to 0.75pt        

\begin{tikzpicture}[x=0.75pt,y=0.75pt,yscale=-1,xscale=1]
%uncomment if require: \path (0,11616); %set diagram left start at 0, and has height of 11616

%Shape: Circle [id:dp2148782782020089] 
\draw   (135.44,8742.35) .. controls (140,8742.34) and (143.7,8746.02) .. (143.7,8750.58) .. controls (143.71,8755.14) and (140.02,8758.84) .. (135.47,8758.85) .. controls (130.91,8758.85) and (127.21,8755.17) .. (127.2,8750.61) .. controls (127.2,8746.05) and (130.88,8742.35) .. (135.44,8742.35) -- cycle ;
%Shape: Circle [id:dp2808434799498767] 
\draw   (156.44,8775.35) .. controls (161,8775.34) and (164.7,8779.02) .. (164.7,8783.58) .. controls (164.71,8788.14) and (161.02,8791.84) .. (156.47,8791.85) .. controls (151.91,8791.85) and (148.21,8788.17) .. (148.2,8783.61) .. controls (148.2,8779.05) and (151.88,8775.35) .. (156.44,8775.35) -- cycle ;
%Shape: Circle [id:dp4630460625057502] 
\draw   (206.44,8775.35) .. controls (211,8775.34) and (214.7,8779.02) .. (214.7,8783.58) .. controls (214.71,8788.14) and (211.02,8791.84) .. (206.47,8791.85) .. controls (201.91,8791.85) and (198.21,8788.17) .. (198.2,8783.61) .. controls (198.2,8779.05) and (201.88,8775.35) .. (206.44,8775.35) -- cycle ;
%Shape: Circle [id:dp3435132640934846] 
\draw   (232.44,8741.35) .. controls (237,8741.34) and (240.7,8745.02) .. (240.7,8749.58) .. controls (240.71,8754.14) and (237.02,8757.84) .. (232.47,8757.85) .. controls (227.91,8757.85) and (224.21,8754.17) .. (224.2,8749.61) .. controls (224.2,8745.05) and (227.88,8741.35) .. (232.44,8741.35) -- cycle ;
%Shape: Circle [id:dp5480450224684223] 
\draw   (224.44,8821.35) .. controls (229,8821.34) and (232.7,8825.02) .. (232.7,8829.58) .. controls (232.71,8834.14) and (229.02,8837.84) .. (224.47,8837.85) .. controls (219.91,8837.85) and (216.21,8834.17) .. (216.2,8829.61) .. controls (216.2,8825.05) and (219.88,8821.35) .. (224.44,8821.35) -- cycle ;
%Shape: Circle [id:dp17656547746667095] 
\draw   (265.44,8821.35) .. controls (270,8821.34) and (273.7,8825.02) .. (273.7,8829.58) .. controls (273.71,8834.14) and (270.02,8837.84) .. (265.47,8837.85) .. controls (260.91,8837.85) and (257.21,8834.17) .. (257.2,8829.61) .. controls (257.2,8825.05) and (260.88,8821.35) .. (265.44,8821.35) -- cycle ;
%Shape: Circle [id:dp28056328551272935] 
\draw   (209.44,8866.35) .. controls (214,8866.34) and (217.7,8870.02) .. (217.7,8874.58) .. controls (217.71,8879.14) and (214.02,8882.84) .. (209.47,8882.85) .. controls (204.91,8882.85) and (201.21,8879.17) .. (201.2,8874.61) .. controls (201.2,8870.05) and (204.88,8866.35) .. (209.44,8866.35) -- cycle ;
%Shape: Circle [id:dp20399212070669015] 
\draw   (231.44,8899.35) .. controls (236,8899.34) and (239.7,8903.02) .. (239.7,8907.58) .. controls (239.71,8912.14) and (236.02,8915.84) .. (231.47,8915.85) .. controls (226.91,8915.85) and (223.21,8912.17) .. (223.2,8907.61) .. controls (223.2,8903.05) and (226.88,8899.35) .. (231.44,8899.35) -- cycle ;
%Shape: Circle [id:dp30150493666582046] 
\draw   (158.44,8866.35) .. controls (163,8866.34) and (166.7,8870.02) .. (166.7,8874.58) .. controls (166.71,8879.14) and (163.02,8882.84) .. (158.47,8882.85) .. controls (153.91,8882.85) and (150.21,8879.17) .. (150.2,8874.61) .. controls (150.2,8870.05) and (153.88,8866.35) .. (158.44,8866.35) -- cycle ;
%Shape: Circle [id:dp31656427905225826] 
\draw   (135.44,8900.35) .. controls (140,8900.34) and (143.7,8904.02) .. (143.7,8908.58) .. controls (143.71,8913.14) and (140.02,8916.84) .. (135.47,8916.85) .. controls (130.91,8916.85) and (127.21,8913.17) .. (127.2,8908.61) .. controls (127.2,8904.05) and (130.88,8900.35) .. (135.44,8900.35) -- cycle ;
%Shape: Circle [id:dp8555191737317314] 
\draw   (139.44,8820.35) .. controls (144,8820.34) and (147.7,8824.02) .. (147.7,8828.58) .. controls (147.71,8833.14) and (144.02,8836.84) .. (139.47,8836.85) .. controls (134.91,8836.85) and (131.21,8833.17) .. (131.2,8828.61) .. controls (131.2,8824.05) and (134.88,8820.35) .. (139.44,8820.35) -- cycle ;
%Shape: Circle [id:dp05436618074318522] 
\draw   (99.44,8820.35) .. controls (104,8820.34) and (107.7,8824.02) .. (107.7,8828.58) .. controls (107.71,8833.14) and (104.02,8836.84) .. (99.47,8836.85) .. controls (94.91,8836.85) and (91.21,8833.17) .. (91.2,8828.61) .. controls (91.2,8824.05) and (94.88,8820.35) .. (99.44,8820.35) -- cycle ;
%Curve Lines [id:da5581316611253468] 
\draw    (232.44,8828.35) .. controls (240.73,8822.08) and (246.46,8821.55) .. (254.83,8826.69) ;
\draw [shift={(257.24,8828.26)}, rotate = 214.79] [fill={rgb, 255:red, 0; green, 0; blue, 0 }  ][line width=0.08]  [draw opacity=0] (7.14,-3.43) -- (0,0) -- (7.14,3.43) -- cycle    ;
%Curve Lines [id:da8817715051524397] 
\draw    (257.5,8831.5) .. controls (251.98,8837.94) and (243.92,8839.3) .. (234.78,8833.31) ;
\draw [shift={(232.37,8831.58)}, rotate = 398.02] [fill={rgb, 255:red, 0; green, 0; blue, 0 }  ][line width=0.08]  [draw opacity=0] (7.14,-3.43) -- (0,0) -- (7.14,3.43) -- cycle    ;
%Curve Lines [id:da433926092322388] 
\draw    (107.44,8826.35) .. controls (115.73,8820.08) and (120.62,8819.55) .. (128.85,8824.69) ;
\draw [shift={(131.24,8826.26)}, rotate = 214.79] [fill={rgb, 255:red, 0; green, 0; blue, 0 }  ][line width=0.08]  [draw opacity=0] (7.14,-3.43) -- (0,0) -- (7.14,3.43) -- cycle    ;
%Curve Lines [id:da08409052860082356] 
\draw    (131.5,8829.5) .. controls (125.98,8835.94) and (118.77,8837.3) .. (109.76,8831.31) ;
\draw [shift={(107.37,8829.58)}, rotate = 398.02] [fill={rgb, 255:red, 0; green, 0; blue, 0 }  ][line width=0.08]  [draw opacity=0] (7.14,-3.43) -- (0,0) -- (7.14,3.43) -- cycle    ;
%Curve Lines [id:da12568404518151088] 
\draw [color={rgb, 255:red, 0; green, 0; blue, 0 }  ,draw opacity=1 ][line width=0.75]  [dash pattern={on 4.5pt off 4.5pt}]  (162.47,8788.85) .. controls (184.58,8800.03) and (160.58,8816.15) .. (158.6,8794.72) ;
\draw [shift={(158.47,8791.85)}, rotate = 449.93] [fill={rgb, 255:red, 0; green, 0; blue, 0 }  ,fill opacity=1 ][line width=0.08]  [draw opacity=0] (7.14,-3.43) -- (0,0) -- (7.14,3.43) -- cycle    ;
%Curve Lines [id:da9111234086939466] 
\draw [color={rgb, 255:red, 0; green, 0; blue, 0 }  ,draw opacity=1 ]   (146.7,8825.58) .. controls (164.84,8814.89) and (171.06,8839.66) .. (149.19,8832.49) ;
\draw [shift={(146.7,8831.58)}, rotate = 381.8] [fill={rgb, 255:red, 0; green, 0; blue, 0 }  ,fill opacity=1 ][line width=0.08]  [draw opacity=0] (7.14,-3.43) -- (0,0) -- (7.14,3.43) -- cycle    ;
%Curve Lines [id:da31123748683160346] 
\draw [color={rgb, 255:red, 0; green, 0; blue, 0 }  ,draw opacity=1 ] [dash pattern={on 4.5pt off 4.5pt}]  (206.47,8791.85) .. controls (198.86,8814.44) and (183.01,8800.86) .. (198.15,8790.01) ;
\draw [shift={(200.5,8788.5)}, rotate = 509.93] [fill={rgb, 255:red, 0; green, 0; blue, 0 }  ,fill opacity=1 ][line width=0.08]  [draw opacity=0] (7.14,-3.43) -- (0,0) -- (7.14,3.43) -- cycle    ;
%Curve Lines [id:da24328209634821696] 
\draw [color={rgb, 255:red, 0; green, 0; blue, 0 }  ,draw opacity=1 ]   (217.5,8833.5) .. controls (198.3,8840.22) and (193.85,8822.06) .. (215.38,8825.11) ;
\draw [shift={(218.2,8825.61)}, rotate = 191.69] [fill={rgb, 255:red, 0; green, 0; blue, 0 }  ,fill opacity=1 ][line width=0.08]  [draw opacity=0] (7.14,-3.43) -- (0,0) -- (7.14,3.43) -- cycle    ;
%Curve Lines [id:da5815770460664151] 
\draw [color={rgb, 255:red, 0; green, 0; blue, 0 }  ,draw opacity=1 ]   (159,8866.5) .. controls (161.4,8846.34) and (179.93,8851.06) .. (166.37,8868.28) ;
\draw [shift={(164.5,8870.5)}, rotate = 311.82] [fill={rgb, 255:red, 0; green, 0; blue, 0 }  ,fill opacity=1 ][line width=0.08]  [draw opacity=0] (7.14,-3.43) -- (0,0) -- (7.14,3.43) -- cycle    ;
%Curve Lines [id:da8951929085347534] 
\draw [color={rgb, 255:red, 0; green, 0; blue, 0 }  ,draw opacity=1 ]   (203.2,8869.61) .. controls (183.33,8855.1) and (202.81,8842.54) .. (207.79,8863.98) ;
\draw [shift={(208.33,8866.83)}, rotate = 261.4] [fill={rgb, 255:red, 0; green, 0; blue, 0 }  ,fill opacity=1 ][line width=0.08]  [draw opacity=0] (7.14,-3.43) -- (0,0) -- (7.14,3.43) -- cycle    ;
%Curve Lines [id:da20142595652687945] 
\draw [color={rgb, 255:red, 0; green, 0; blue, 0 }  ,draw opacity=1 ] [dash pattern={on 4.5pt off 4.5pt}]  (233.44,8741.35) .. controls (235.84,8721.19) and (254.37,8725.91) .. (240.8,8743.13) ;
\draw [shift={(238.94,8745.35)}, rotate = 311.82] [fill={rgb, 255:red, 0; green, 0; blue, 0 }  ,fill opacity=1 ][line width=0.08]  [draw opacity=0] (7.14,-3.43) -- (0,0) -- (7.14,3.43) -- cycle    ;
%Curve Lines [id:da783779044923236] 
\draw [color={rgb, 255:red, 0; green, 0; blue, 0 }  ,draw opacity=1 ] [dash pattern={on 4.5pt off 4.5pt}]  (128.31,8746.12) .. controls (108.43,8731.62) and (127.91,8719.06) .. (132.9,8740.49) ;
\draw [shift={(133.44,8743.35)}, rotate = 261.4] [fill={rgb, 255:red, 0; green, 0; blue, 0 }  ,fill opacity=1 ][line width=0.08]  [draw opacity=0] (7.14,-3.43) -- (0,0) -- (7.14,3.43) -- cycle    ;
%Curve Lines [id:da4772971581312033] 
\draw  [dash pattern={on 4.5pt off 4.5pt}]  (142.47,8755.85) .. controls (149.86,8758.29) and (156.37,8761.87) .. (156.54,8772.43) ;
\draw [shift={(156.44,8775.35)}, rotate = 274.72] [fill={rgb, 255:red, 0; green, 0; blue, 0 }  ][line width=0.08]  [draw opacity=0] (7.14,-3.43) -- (0,0) -- (7.14,3.43) -- cycle    ;
%Curve Lines [id:da4008079025329423] 
\draw  [dash pattern={on 4.5pt off 4.5pt}]  (230.47,8756.85) .. controls (228.65,8769.48) and (226.8,8774.7) .. (216.24,8777.79) ;
\draw [shift={(213.5,8778.5)}, rotate = 347.01] [fill={rgb, 255:red, 0; green, 0; blue, 0 }  ][line width=0.08]  [draw opacity=0] (7.14,-3.43) -- (0,0) -- (7.14,3.43) -- cycle    ;
%Curve Lines [id:da8214395576271778] 
\draw  [dash pattern={on 4.5pt off 4.5pt}]  (151.5,8777.5) .. controls (139.6,8776.58) and (136.92,8770.65) .. (138.06,8760.72) ;
\draw [shift={(138.47,8757.85)}, rotate = 459.58] [fill={rgb, 255:red, 0; green, 0; blue, 0 }  ][line width=0.08]  [draw opacity=0] (7.14,-3.43) -- (0,0) -- (7.14,3.43) -- cycle    ;
%Curve Lines [id:da8055029877181199] 
\draw  [dash pattern={on 4.5pt off 4.5pt}]  (210.44,8775.35) .. controls (209.6,8764.74) and (215.12,8758.82) .. (222.53,8754.91) ;
\draw [shift={(225.2,8753.61)}, rotate = 515.9200000000001] [fill={rgb, 255:red, 0; green, 0; blue, 0 }  ][line width=0.08]  [draw opacity=0] (7.14,-3.43) -- (0,0) -- (7.14,3.43) -- cycle    ;
%Curve Lines [id:da9385427876452483] 
\draw    (138.44,8901.35) .. controls (137.6,8890.74) and (143.12,8884.82) .. (150.53,8880.91) ;
\draw [shift={(153.2,8879.61)}, rotate = 515.9200000000001] [fill={rgb, 255:red, 0; green, 0; blue, 0 }  ][line width=0.08]  [draw opacity=0] (7.14,-3.43) -- (0,0) -- (7.14,3.43) -- cycle    ;
%Curve Lines [id:da6211873110144355] 
\draw    (158.47,8882.85) .. controls (156.65,8895.48) and (154.8,8900.7) .. (144.24,8903.79) ;
\draw [shift={(141.5,8904.5)}, rotate = 347.01] [fill={rgb, 255:red, 0; green, 0; blue, 0 }  ][line width=0.08]  [draw opacity=0] (7.14,-3.43) -- (0,0) -- (7.14,3.43) -- cycle    ;
%Curve Lines [id:da3524689152808309] 
\draw    (226.44,8901.35) .. controls (214.54,8900.43) and (211.86,8894.49) .. (213,8884.56) ;
\draw [shift={(213.41,8881.69)}, rotate = 459.58] [fill={rgb, 255:red, 0; green, 0; blue, 0 }  ][line width=0.08]  [draw opacity=0] (7.14,-3.43) -- (0,0) -- (7.14,3.43) -- cycle    ;
%Curve Lines [id:da29799205753697944] 
\draw    (216.47,8879.85) .. controls (223.86,8882.29) and (230.37,8886.72) .. (230.54,8897.41) ;
\draw [shift={(230.44,8900.35)}, rotate = 274.72] [fill={rgb, 255:red, 0; green, 0; blue, 0 }  ][line width=0.08]  [draw opacity=0] (7.14,-3.43) -- (0,0) -- (7.14,3.43) -- cycle    ;
%Straight Lines [id:da04285712935079422] 
\draw    (105.81,8819.51) -- (148.5,8786.5) ;
\draw [shift={(103.44,8821.35)}, rotate = 322.29] [fill={rgb, 255:red, 0; green, 0; blue, 0 }  ][line width=0.08]  [draw opacity=0] (7.14,-3.43) -- (0,0) -- (7.14,3.43) -- cycle    ;
%Straight Lines [id:da924839646984126] 
\draw    (140.53,8817.55) -- (151.5,8789.5) ;
\draw [shift={(139.44,8820.35)}, rotate = 291.36] [fill={rgb, 255:red, 0; green, 0; blue, 0 }  ][line width=0.08]  [draw opacity=0] (7.14,-3.43) -- (0,0) -- (7.14,3.43) -- cycle    ;
%Straight Lines [id:da5412015612070826] 
\draw    (152.36,8865.73) -- (140.47,8836.85) ;
\draw [shift={(153.5,8868.5)}, rotate = 247.62] [fill={rgb, 255:red, 0; green, 0; blue, 0 }  ][line width=0.08]  [draw opacity=0] (7.14,-3.43) -- (0,0) -- (7.14,3.43) -- cycle    ;
%Straight Lines [id:da5348933441098955] 
\draw    (132.62,8898.35) -- (136.47,8835.85) ;
\draw [shift={(132.44,8901.35)}, rotate = 273.52] [fill={rgb, 255:red, 0; green, 0; blue, 0 }  ][line width=0.08]  [draw opacity=0] (7.14,-3.43) -- (0,0) -- (7.14,3.43) -- cycle    ;
%Straight Lines [id:da415553950357612] 
\draw    (198.2,8874.61) -- (166.7,8874.58) ;
\draw [shift={(201.2,8874.61)}, rotate = 180.05] [fill={rgb, 255:red, 0; green, 0; blue, 0 }  ][line width=0.08]  [draw opacity=0] (7.14,-3.43) -- (0,0) -- (7.14,3.43) -- cycle    ;
%Straight Lines [id:da5442719629278565] 
\draw    (220.5,8906.31) -- (164.5,8879.5) ;
\draw [shift={(223.2,8907.61)}, rotate = 205.59] [fill={rgb, 255:red, 0; green, 0; blue, 0 }  ][line width=0.08]  [draw opacity=0] (7.14,-3.43) -- (0,0) -- (7.14,3.43) -- cycle    ;
%Straight Lines [id:da7123843609413225] 
\draw    (223.35,8840.63) -- (212.5,8867.5) ;
\draw [shift={(224.47,8837.85)}, rotate = 111.98] [fill={rgb, 255:red, 0; green, 0; blue, 0 }  ][line width=0.08]  [draw opacity=0] (7.14,-3.43) -- (0,0) -- (7.14,3.43) -- cycle    ;
%Straight Lines [id:da27095386646545716] 
\draw    (260.96,8839.5) -- (215.5,8869.5) ;
\draw [shift={(263.47,8837.85)}, rotate = 146.58] [fill={rgb, 255:red, 0; green, 0; blue, 0 }  ][line width=0.08]  [draw opacity=0] (7.14,-3.43) -- (0,0) -- (7.14,3.43) -- cycle    ;
%Straight Lines [id:da8989167374808824] 
\draw    (211.47,8793.34) -- (221.44,8822.35) ;
\draw [shift={(210.5,8790.5)}, rotate = 71.04] [fill={rgb, 255:red, 0; green, 0; blue, 0 }  ][line width=0.08]  [draw opacity=0] (7.14,-3.43) -- (0,0) -- (7.14,3.43) -- cycle    ;
%Straight Lines [id:da09402016882912934] 
\draw    (236.04,8759.46) -- (226.44,8821.35) ;
\draw [shift={(236.5,8756.5)}, rotate = 98.82] [fill={rgb, 255:red, 0; green, 0; blue, 0 }  ][line width=0.08]  [draw opacity=0] (7.14,-3.43) -- (0,0) -- (7.14,3.43) -- cycle    ;
%Straight Lines [id:da8592285190252522] 
\draw    (167.7,8783.58) -- (198.2,8783.61) ;
\draw [shift={(164.7,8783.58)}, rotate = 0.05] [fill={rgb, 255:red, 0; green, 0; blue, 0 }  ][line width=0.08]  [draw opacity=0] (7.14,-3.43) -- (0,0) -- (7.14,3.43) -- cycle    ;
%Straight Lines [id:da2125111068994756] 
\draw    (146.4,8751.9) -- (200.5,8778.5) ;
\draw [shift={(143.7,8750.58)}, rotate = 26.18] [fill={rgb, 255:red, 0; green, 0; blue, 0 }  ][line width=0.08]  [draw opacity=0] (7.14,-3.43) -- (0,0) -- (7.14,3.43) -- cycle    ;
%Shape: Ellipse [id:dp3493789699928824] 
\draw  [color={rgb, 255:red, 155; green, 155; blue, 155 }  ,draw opacity=0.6 ] (118.89,8726.56) .. controls (127.56,8720.83) and (146.54,8734.27) .. (161.26,8756.57) .. controls (175.99,8778.86) and (180.89,8801.59) .. (172.22,8807.32) .. controls (163.54,8813.04) and (144.57,8799.61) .. (129.84,8777.31) .. controls (115.12,8755.01) and (110.21,8732.29) .. (118.89,8726.56) -- cycle ;
%Shape: Ellipse [id:dp9681882621382776] 
\draw  [color={rgb, 255:red, 155; green, 155; blue, 155 }  ,draw opacity=0.6 ] (70.17,8826.89) .. controls (70.18,8816.49) and (91.85,8808.09) .. (118.57,8808.11) .. controls (145.29,8808.14) and (166.95,8816.59) .. (166.94,8826.99) .. controls (166.93,8837.39) and (145.25,8845.79) .. (118.53,8845.76) .. controls (91.81,8845.74) and (70.16,8837.29) .. (70.17,8826.89) -- cycle ;
%Shape: Ellipse [id:dp03614915665867824] 
\draw  [color={rgb, 255:red, 155; green, 155; blue, 155 }  ,draw opacity=0.6 ] (196.17,8828.89) .. controls (196.18,8818.49) and (217.85,8810.09) .. (244.57,8810.11) .. controls (271.29,8810.14) and (292.95,8818.59) .. (292.94,8828.99) .. controls (292.93,8839.39) and (271.25,8847.79) .. (244.53,8847.76) .. controls (217.81,8847.74) and (196.16,8839.29) .. (196.17,8828.89) -- cycle ;
%Shape: Ellipse [id:dp635375891390898] 
\draw  [color={rgb, 255:red, 155; green, 155; blue, 155 }  ,draw opacity=0.6 ] (192.8,8806.25) .. controls (184.13,8800.51) and (189.09,8777.8) .. (203.87,8755.53) .. controls (218.64,8733.27) and (237.64,8719.88) .. (246.31,8725.62) .. controls (254.97,8731.37) and (250.01,8754.08) .. (235.24,8776.35) .. controls (220.46,8798.61) and (201.46,8812) .. (192.8,8806.25) -- cycle ;
%Shape: Ellipse [id:dp9297173136377761] 
\draw  [color={rgb, 255:red, 155; green, 155; blue, 155 }  ,draw opacity=0.6 ] (193.89,8848.56) .. controls (202.56,8842.83) and (221.54,8856.27) .. (236.26,8878.57) .. controls (250.99,8900.86) and (255.89,8923.59) .. (247.22,8929.32) .. controls (238.54,8935.04) and (219.57,8921.61) .. (204.84,8899.31) .. controls (190.12,8877.01) and (185.21,8854.29) .. (193.89,8848.56) -- cycle ;
%Shape: Ellipse [id:dp23351023149204075] 
\draw  [color={rgb, 255:red, 155; green, 155; blue, 155 }  ,draw opacity=0.6 ] (119.8,8930.25) .. controls (111.13,8924.51) and (116.09,8901.8) .. (130.87,8879.53) .. controls (145.64,8857.27) and (164.64,8843.88) .. (173.31,8849.62) .. controls (181.97,8855.37) and (177.01,8878.08) .. (162.24,8900.35) .. controls (147.46,8922.61) and (128.46,8936) .. (119.8,8930.25) -- cycle ;

% Text Node
\draw (136.45,8750.6) node  [font=\tiny,rotate=-1.27] [align=left] {$\displaystyle 1^{1}$};
% Text Node
\draw (157.45,8783.6) node  [font=\tiny,rotate=-1.27] [align=left] {$\displaystyle 1^{2}$};
% Text Node
\draw (207.45,8783.6) node  [font=\tiny,rotate=-1.27] [align=left] {$\displaystyle 6^{2}$};
% Text Node
\draw (233.45,8749.6) node  [font=\tiny,rotate=-1.27] [align=left] {$\displaystyle 6^{1}$};
% Text Node
\draw (225.45,8829.6) node  [font=\tiny,rotate=-1.27] [align=left] {$\displaystyle 5^{2}$};
% Text Node
\draw (266.45,8829.6) node  [font=\tiny,rotate=-1.27] [align=left] {$\displaystyle 5^{1}$};
% Text Node
\draw (210.45,8874.6) node  [font=\tiny,rotate=-1.27] [align=left] {$\displaystyle 4^{2}$};
% Text Node
\draw (232.45,8907.6) node  [font=\tiny,rotate=-1.27] [align=left] {$\displaystyle 4^{1}$};
% Text Node
\draw (159.45,8874.6) node  [font=\tiny,rotate=-1.27] [align=left] {$\displaystyle 3^{2}$};
% Text Node
\draw (136.45,8908.6) node  [font=\tiny,rotate=-1.27] [align=left] {$\displaystyle 3^{1}$};
% Text Node
\draw (140.45,8828.6) node  [font=\tiny,rotate=-1.27] [align=left] {$\displaystyle 2^{2}$};
% Text Node
\draw (100.45,8828.6) node  [font=\tiny,rotate=-1.27] [align=left] {$\displaystyle 2^{1}$};
% Text Node
\draw (118.75,8712.33) node   [align=left] {\begin{minipage}[lt]{8.67pt}\setlength\topsep{0pt}
$\displaystyle {\displaystyle 1}$
\end{minipage}};
% Text Node
\draw (242.75,8712.33) node   [align=left] {\begin{minipage}[lt]{8.67pt}\setlength\topsep{0pt}
$\displaystyle 6$
\end{minipage}};
% Text Node
\draw (73.75,8803.33) node   [align=left] {\begin{minipage}[lt]{8.67pt}\setlength\topsep{0pt}
$\displaystyle 2$
\end{minipage}};
% Text Node
\draw (284.75,8805.33) node   [align=left] {\begin{minipage}[lt]{8.67pt}\setlength\topsep{0pt}
$\displaystyle 5$
\end{minipage}};
% Text Node
\draw (111.75,8935.33) node   [align=left] {\begin{minipage}[lt]{8.67pt}\setlength\topsep{0pt}
$\displaystyle 3$
\end{minipage}};
% Text Node
\draw (253.97,8937.15) node   [align=left] {\begin{minipage}[lt]{8.67pt}\setlength\topsep{0pt}
$\displaystyle {\displaystyle 4}$
\end{minipage}};

\end{tikzpicture}

\caption{ Global graph ${\mathcal{G}}$ for Example \ref{ex_stab0}.}
\label{global}
\end{figure}

\begin{ex}
Let $\gamma_{i1}$, and  $\gamma_{i2}$ be some nonzero real parameters, and assume that $\beta_{i1}, \beta_{i2}, \beta_{i3}>0$, $i=2,3,\ldots,n-1$. Now, consider an LTI network $\mathcal{N}$, whose system matrix $A$, described in (\ref{sm}),  is defined as:

\begin{equation}A=\begin{bmatrix}
A_{11} & 0 & \ldots & \ldots & 0 & A_{1n}\\
A_{21} & A_{22} & 0 & \ldots & \ldots & 0\\
0 & A_{32} & A_{33} & 0 & \ldots & 0\\
\vdots  & \ddots & \ddots & \ddots & \ddots & \vdots\\
0 & \ldots & \ldots & 0 & A_{n,n-1} & A_{nn} 
\end{bmatrix}, \label{A2}\end{equation}
where  we have  
\begin{equation}
\begin{aligned}A_{i,i-1}=\begin{bmatrix}
0 & \gamma_{i1}\\ 0 & \gamma_{i2}
\end{bmatrix}, \:\:\:\:\:\:\:\:\qquad\qquad1\leq i\leq n,\\  A_{ii}=\begin{bmatrix}
0 & -\beta_{i1} \\ \beta_{i2} & -\beta_{i3}
\end{bmatrix},\:\:\:\:\:\:\:\:\qquad\qquad1< i< n.\label{Ai}\end{aligned}\end{equation}
Then, besides the zero/nonzero/arbitrary pattern of any matrix $A_{ii}$, $i=2,3,\ldots, n-1$, we know the sign of its nonzero parameters. Now, one can see that for all numerical realizations of system matrix $A$, any subsystem $i$, $i=2,3,\ldots, n-1$, is stable (because the eigenvalues of any subsystem $i$ are the roots of the characteristic equation $\lambda^2+\beta_{i3}\lambda+\beta_{i1}\beta_{i2}=0$).  Now, let $\mathcal{A}'\in\{?\}^{2\times 2}$ be a $2\times 2$ pattern matrix with all entries equal to $?$. We let $A_{11},A_{nn}\in \mathcal{P}(\mathcal{A}')$. Thus, there is no available information about the zero/nonzero/arbitrary pattern of $A_{11}$ and $A_{nn}$. However, we know that these matrices are also associated with stable subsystems.  %Thus, we have  $f^{\mathcal{L}}_{\Delta}(i)=*$, for all $i=1,2, \ldots, n$. For $n=6$, the corresponding $\Delta$-network graph $G_{\mathcal{N}}^{\Delta}$ is depicted in Fig.~\ref{globgraph_ex}. 
Thus,  with $\Delta=\{y\in\mathbb{C}\mid \Re(y)\geq 0\}$, we have  $f^{\mathcal{L}}_{\Delta}(i)=*$, for all $i=1,2, \ldots, n$.
For $n=6$, the global graph $\mathcal{G}$ associated with this network is depicted in Fig. \ref{global}.  The $\Delta$-specified pattern class $\mathcal{S}(f^{\mathcal{L}}_{\Delta}, \mathcal{A})$ includes all system matrices of the same zero/nonzero/arbitrary structure, whose all subsystems  are stable. 
 In this example, we aim to examine the strong structural stabilizabity of the entire network, that is, the strong structural controllability with respect to $\Delta=\{y\in\mathbb{C}\mid \Re(y)\geq 0\}$.
 %The most precise and relevant results existing in the literature to consider this problem are those ones on strong structural controllability of nonzero eigenvalues (see e.g. \cite{trefois2015zero,jia2020unifying}). However, 
 Now, one can verify that by applying the existing results in the literature  on strong structural controllability (see e.g. \cite{trefois2015zero,jia2020unifying}),  a set of control nodes should include at least 7 vertices. However, we will show that one control subsystem, including two vertices, can render this network strongly structurally stabilizable. 
\label{ex_stab0}
\end{ex}

In this paper, we will also study the maximum geometric multiplicities of eigenvalues of $A\in \mathcal{S}(f^{\mathcal{L}}_{\Delta},\mathcal{A})$, for all eigenvalues that belong to the set $\Delta$. Then, the next problems of this work can be stated as follows. 

\emph{Problem 2:} Provide an upper bound on the maximum geometric multiplicity of eigenvalues of all $A\in \mathcal{S}(f^{\mathcal{L}}_{\Delta},\mathcal{A})$ that belong to a set $\Delta\in\mathbb{C}$. 

\emph{Problem 3:} Given a nonempty set $\Delta\in\mathbb{C}$, find a combinatorial condition under which no matrix $A\in \mathcal{S}(f^{\mathcal{L}}_{\Delta},\mathcal{A})$ has any eigenvalue in $\Delta$.

\section{$\Delta$-network Graphs and Coloring Process}\label{Delta_color}

In this section, we introduce some essential notions that are employed in providing the main results of this work. First, given a nonempty set $\Delta$, a $\Delta$-characteristic vector $f^{\mathcal{L}}_{\Delta}$, and the global graph of an LTI network, we discuss how the corresponding $\Delta$-network graph is formed. Next, a coloring process  and zero forcing sets are introduced.

\subsection{$\Delta$-network Graphs}\label{section}
Consider an LTI network $\mathcal{N}$, as described in (\ref{e1}), whose system matrix $A$ is defined in (\ref{sm}), and for some pattern matrix $\mathcal{A}\in\{0,*,?\}^{N\times N}$, we have $A\in\mathcal{P}(\mathcal{A})$. Let $\Delta\subseteq \mathbb{C}$ be a nonempty set. Now, assume that we are either given a vector $f^{\mathcal{L}}_{\Delta}\in \{0,*,?\}^{n\times n}$, representing the extra information about $\Lambda(A_{ii})\cap \Delta$, for $i=1,2,\ldots,n$ (similar to Examples \ref{stab} and \ref{int}), or we have the set of system matrices $\mathcal{S}^*$ associated with this network, and $f^{\mathcal{L}}_{\Delta}$ is defined according to (\ref{f}) (see Examples \ref{ex5}, \ref{ex6}, \ref{ex7}, and \ref{ex8}).   

Now, having $f^{\mathcal{L}}_{\Delta}$ and $\mathcal{A}$, we associate to the network $\mathcal{N}$ a $\Delta$-network graph $G_{\mathcal{N}}^{\Delta}=(V_\mathcal{N}, E_\mathcal{N})$, where $V_{\mathcal{N}}=\{1,2,\ldots, n\}$ and $E_{\mathcal{N}}=E_{\mathcal{N}}^*\cup E_{\mathcal{N}}^?$. The self-loops of this graph are defined based on the entries of the vector $f^{\mathcal{L}}_{\Delta}$ as follows. For $i=1,2,\ldots, n$, there exists a self-loop $(i,i)\in E_{\mathcal{N}}$ if and only if $f^{\mathcal{L}}_{\Delta}(i)\neq 0$. Moreover, $(i,i)\in E_{\mathcal{N}}^*$ if and only if $f^{\mathcal{L}}_{\Delta}(i)= *$, and $(i,i)\in E_{\mathcal{N}}^?$ if and only if $f^{\mathcal{L}}_{\Delta}(i)= ?$.

For $i\neq j$, we have  $(j,i)\notin E_\mathcal{N}$ if and only if $\mathcal{A}_{ij}=0$. Moreover,  $(j,i)\in E_{\mathcal{N}}^*$ if and only if the pattern matrix $\mathcal{A}_{ij}$ has full row rank, and $(j,i)\in E_{\mathcal{N}}^?$ otherwise. In the next section of the paper, a graph-theoretic condition is presented; by checking this condition in the associated global graph, we can find whether a pattern matrix $\mathcal{A}_{ij}$ has full row rank or not. We differentiate the edges in $E_{\mathcal{N}}^*$ and $E_{\mathcal{N}}^?$ through the solid and dotted arrows, respectively.

Notice that we may have no information about the pattern matrix of any subsystem $k$, $k=1,2,\ldots, n$. In other words, for $k=1,2,\ldots, n$ and  $1\leq i,j\leq l_i$, we have  $\mathcal{A}_{kk}(i,j)=?$. However, if $f^{\mathcal{L}}_{\Delta}$ is available, then one can describe the self-loops of the $\Delta$-network graph. For instance, for $\Delta=\{y\in\mathbb{C}\mid \Re(y)\geq 0\}$, if no information about the structure of the subsystems is available, but we know   which subsystems are stable, then we can determine the set of nodes of the $\Delta$-network graph that have self-loops, which are drawn by a solid arrows.

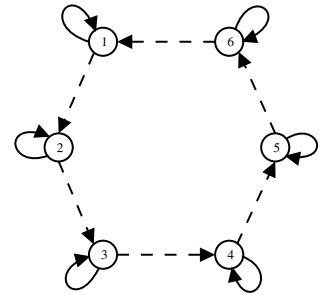
\begin{figure}[b]
\centering

\tikzset{every picture/.style={line width=0.75pt}} %set default line width to 0.75pt        

\begin{tikzpicture}[x=0.75pt,y=0.75pt,yscale=-.86,xscale=.86]
%uncomment if require: \path (0,11616); %set diagram left start at 0, and has height of 11616

%Shape: Circle [id:dp06900926076234803] 
\draw   (364.44,8758.35) .. controls (369,8758.34) and (372.7,8762.02) .. (372.7,8766.58) .. controls (372.71,8771.14) and (369.02,8774.84) .. (364.47,8774.85) .. controls (359.91,8774.85) and (356.21,8771.17) .. (356.2,8766.61) .. controls (356.2,8762.05) and (359.88,8758.35) .. (364.44,8758.35) -- cycle ;
%Shape: Circle [id:dp8906913651799573] 
\draw   (438.44,8758.35) .. controls (443,8758.34) and (446.7,8762.02) .. (446.7,8766.58) .. controls (446.71,8771.14) and (443.02,8774.84) .. (438.47,8774.85) .. controls (433.91,8774.85) and (430.21,8771.17) .. (430.2,8766.61) .. controls (430.2,8762.05) and (433.88,8758.35) .. (438.44,8758.35) -- cycle ;
%Shape: Circle [id:dp015687758480868785] 
\draw   (364.44,8883.35) .. controls (369,8883.34) and (372.7,8887.02) .. (372.7,8891.58) .. controls (372.71,8896.14) and (369.02,8899.84) .. (364.47,8899.85) .. controls (359.91,8899.85) and (356.21,8896.17) .. (356.2,8891.61) .. controls (356.2,8887.05) and (359.88,8883.35) .. (364.44,8883.35) -- cycle ;
%Shape: Circle [id:dp19203687551803528] 
\draw   (438.44,8883.35) .. controls (443,8883.34) and (446.7,8887.02) .. (446.7,8891.58) .. controls (446.71,8896.14) and (443.02,8899.84) .. (438.47,8899.85) .. controls (433.91,8899.85) and (430.21,8896.17) .. (430.2,8891.61) .. controls (430.2,8887.05) and (433.88,8883.35) .. (438.44,8883.35) -- cycle ;
%Shape: Circle [id:dp5288891501717434] 
\draw   (338.44,8820.35) .. controls (343,8820.34) and (346.7,8824.02) .. (346.7,8828.58) .. controls (346.71,8833.14) and (343.02,8836.84) .. (338.47,8836.85) .. controls (333.91,8836.85) and (330.21,8833.17) .. (330.2,8828.61) .. controls (330.2,8824.05) and (333.88,8820.35) .. (338.44,8820.35) -- cycle ;
%Shape: Circle [id:dp27909980354856767] 
\draw   (465.44,8820.35) .. controls (470,8820.34) and (473.7,8824.02) .. (473.7,8828.58) .. controls (473.71,8833.14) and (470.02,8836.84) .. (465.47,8836.85) .. controls (460.91,8836.85) and (457.21,8833.17) .. (457.2,8828.61) .. controls (457.2,8824.05) and (460.88,8820.35) .. (465.44,8820.35) -- cycle ;
%Straight Lines [id:da9007538708116749] 
\draw  [dash pattern={on 4.5pt off 4.5pt}]  (430.2,8766.61) -- (375.7,8766.58) ;
\draw [shift={(372.7,8766.58)}, rotate = 360.03] [fill={rgb, 255:red, 0; green, 0; blue, 0 }  ][line width=0.08]  [draw opacity=0] (8.93,-4.29) -- (0,0) -- (8.93,4.29) -- cycle    ;
%Straight Lines [id:da9249379181953263] 
\draw  [dash pattern={on 4.5pt off 4.5pt}]  (359,8773.5) -- (339.64,8817.6) ;
\draw [shift={(338.44,8820.35)}, rotate = 293.7] [fill={rgb, 255:red, 0; green, 0; blue, 0 }  ][line width=0.08]  [draw opacity=0] (8.93,-4.29) -- (0,0) -- (8.93,4.29) -- cycle    ;
%Straight Lines [id:da48188968192954906] 
\draw  [dash pattern={on 4.5pt off 4.5pt}]  (338.47,8836.85) -- (357.83,8882.74) ;
\draw [shift={(359,8885.5)}, rotate = 247.12] [fill={rgb, 255:red, 0; green, 0; blue, 0 }  ][line width=0.08]  [draw opacity=0] (8.93,-4.29) -- (0,0) -- (8.93,4.29) -- cycle    ;
%Straight Lines [id:da8453159939851713] 
\draw  [dash pattern={on 4.5pt off 4.5pt}]  (372.7,8891.58) -- (427.2,8891.61) ;
\draw [shift={(430.2,8891.61)}, rotate = 180.03] [fill={rgb, 255:red, 0; green, 0; blue, 0 }  ][line width=0.08]  [draw opacity=0] (8.93,-4.29) -- (0,0) -- (8.93,4.29) -- cycle    ;
%Straight Lines [id:da2646113292307688] 
\draw  [dash pattern={on 4.5pt off 4.5pt}]  (443,8885.5) -- (464.21,8839.57) ;
\draw [shift={(465.47,8836.85)}, rotate = 474.79] [fill={rgb, 255:red, 0; green, 0; blue, 0 }  ][line width=0.08]  [draw opacity=0] (8.93,-4.29) -- (0,0) -- (8.93,4.29) -- cycle    ;
%Straight Lines [id:da12976759809598581] 
\draw  [dash pattern={on 4.5pt off 4.5pt}]  (465.44,8820.35) -- (445.23,8775.24) ;
\draw [shift={(444,8772.5)}, rotate = 425.86] [fill={rgb, 255:red, 0; green, 0; blue, 0 }  ][line width=0.08]  [draw opacity=0] (8.93,-4.29) -- (0,0) -- (8.93,4.29) -- cycle    ;
%Curve Lines [id:da5461304256910975] 
\draw [color={rgb, 255:red, 0; green, 0; blue, 0 }  ,draw opacity=1 ]   (356.2,8766.61) .. controls (325.78,8756.75) and (346.88,8728.94) .. (359.39,8758.45) ;
\draw [shift={(360.33,8760.83)}, rotate = 249.7] [fill={rgb, 255:red, 0; green, 0; blue, 0 }  ,fill opacity=1 ][line width=0.08]  [draw opacity=0] (8.93,-4.29) -- (0,0) -- (8.93,4.29) -- cycle    ;
%Curve Lines [id:da8256992248075352] 
\draw [color={rgb, 255:red, 0; green, 0; blue, 0 }  ,draw opacity=1 ]   (332.44,8833.35) .. controls (309.16,8842.08) and (304.7,8812.23) .. (330.92,8821.38) ;
\draw [shift={(333.44,8822.35)}, rotate = 202.48] [fill={rgb, 255:red, 0; green, 0; blue, 0 }  ,fill opacity=1 ][line width=0.08]  [draw opacity=0] (8.93,-4.29) -- (0,0) -- (8.93,4.29) -- cycle    ;
%Curve Lines [id:da013351325333925912] 
\draw [color={rgb, 255:red, 0; green, 0; blue, 0 }  ,draw opacity=1 ]   (362,8899.5) .. controls (348.49,8926.52) and (330.32,8902.32) .. (354.23,8893.41) ;
\draw [shift={(357,8892.5)}, rotate = 524.05] [fill={rgb, 255:red, 0; green, 0; blue, 0 }  ,fill opacity=1 ][line width=0.08]  [draw opacity=0] (8.93,-4.29) -- (0,0) -- (8.93,4.29) -- cycle    ;
%Curve Lines [id:da9243232173545521] 
\draw [color={rgb, 255:red, 0; green, 0; blue, 0 }  ,draw opacity=1 ]   (472,8823.5) .. controls (490.43,8811.86) and (500.39,8837.86) .. (473.59,8833.95) ;
\draw [shift={(471,8833.5)}, rotate = 371.31] [fill={rgb, 255:red, 0; green, 0; blue, 0 }  ,fill opacity=1 ][line width=0.08]  [draw opacity=0] (8.93,-4.29) -- (0,0) -- (8.93,4.29) -- cycle    ;
%Curve Lines [id:da91476739246463] 
\draw [color={rgb, 255:red, 0; green, 0; blue, 0 }  ,draw opacity=1 ]   (446,8892.5) .. controls (470.25,8902.2) and (448.4,8928.84) .. (440.67,8901.22) ;
\draw [shift={(440,8898.5)}, rotate = 437.66] [fill={rgb, 255:red, 0; green, 0; blue, 0 }  ,fill opacity=1 ][line width=0.08]  [draw opacity=0] (8.93,-4.29) -- (0,0) -- (8.93,4.29) -- cycle    ;
%Curve Lines [id:da8626194828686671] 
\draw [color={rgb, 255:red, 0; green, 0; blue, 0 }  ,draw opacity=1 ]   (442,8759.5) .. controls (451.26,8731.71) and (477.16,8752.88) .. (448.05,8765.62) ;
\draw [shift={(445.7,8766.58)}, rotate = 339.01] [fill={rgb, 255:red, 0; green, 0; blue, 0 }  ,fill opacity=1 ][line width=0.08]  [draw opacity=0] (8.93,-4.29) -- (0,0) -- (8.93,4.29) -- cycle    ;

% Text Node
\draw (365.45,8766.6) node  [font=\tiny,rotate=-1.27] [align=left] {1};
% Text Node
\draw (439.45,8766.6) node  [font=\tiny,rotate=-1.27] [align=left] {6};
% Text Node
\draw (365.45,8891.6) node  [font=\tiny,rotate=-1.27] [align=left] {3};
% Text Node
\draw (439.45,8891.6) node  [font=\tiny,rotate=-1.27] [align=left] {4};
% Text Node
\draw (339.45,8828.6) node  [font=\tiny,rotate=-1.27] [align=left] {2};
% Text Node
\draw (466.45,8828.6) node  [font=\tiny,rotate=-1.27] [align=left] {5};

\end{tikzpicture}

\caption{$\Delta$-network graph $G_{\mathcal{N}}^{\Delta}$ associated with the network in Example \ref{ex_stab0}.}
\label{globgraph_ex}
\end{figure}

\begin{ex} Consider the LTI network in Example \ref{ex_stab0}. As mentioned before, all subsystem are stable. Now, let $\Delta=\{y\in\mathbb{C}\mid \Re(y)\geq 0\}$. Thus, we have  $f^{\mathcal{L}}_{\Delta}(i)=*$, for all $i=1, \ldots, n$. Therefore, all nodes in the $\Delta$-network graph have self-loops represented by solid arrows. 
   It is also obvious that $A_{i,i-1}$, for $1\leq i\leq n$, is rank-deficient. Thus, for $i=1, \ldots, n$,  there is a dotted edge from node $i-1$ to node $i$ in the $\Delta$-network graph.   For $n=6$, the corresponding $\Delta$-network graph $G_{\mathcal{N}}^{\Delta}$ is depicted in Fig.~\ref{globgraph_ex}.
\label{ex_stab}
\end{ex}

\begin{figure}[t]
\centering

%%%%%%%%%%%%%%%%%%%%%%%%%%

\tikzset{every picture/.style={line width=0.75pt}} %set default line width to 0.75pt        

\begin{tikzpicture}[x=0.75pt,y=0.75pt,yscale=-.89,xscale=.89]
%uncomment if require: \path (0,11411); %set diagram left start at 0, and has height of 11411

%Shape: Circle [id:dp06771290151330489] 
\draw   (107.44,10876.68) .. controls (112,10876.67) and (115.7,10880.36) .. (115.7,10884.91) .. controls (115.71,10889.47) and (112.02,10893.17) .. (107.47,10893.18) .. controls (102.91,10893.19) and (99.21,10889.5) .. (99.2,10884.94) .. controls (99.2,10880.39) and (102.88,10876.69) .. (107.44,10876.68) -- cycle ;
%Shape: Circle [id:dp4474229546172923] 
\draw   (161.11,10836.68) .. controls (165.66,10836.67) and (169.36,10840.36) .. (169.37,10844.91) .. controls (169.38,10849.47) and (165.69,10853.17) .. (161.13,10853.18) .. controls (156.58,10853.19) and (152.88,10849.5) .. (152.87,10844.94) .. controls (152.86,10840.39) and (156.55,10836.69) .. (161.11,10836.68) -- cycle ;
%Shape: Circle [id:dp4974451197909937] 
\draw   (126.44,10941.68) .. controls (131,10941.67) and (134.7,10945.36) .. (134.7,10949.91) .. controls (134.71,10954.47) and (131.02,10958.17) .. (126.47,10958.18) .. controls (121.91,10958.19) and (118.21,10954.5) .. (118.2,10949.94) .. controls (118.2,10945.39) and (121.88,10941.69) .. (126.44,10941.68) -- cycle ;
%Shape: Circle [id:dp9033828874066825] 
\draw   (195.77,10942.01) .. controls (200.33,10942) and (204.03,10945.69) .. (204.04,10950.25) .. controls (204.04,10954.8) and (200.36,10958.5) .. (195.8,10958.51) .. controls (191.24,10958.52) and (187.54,10954.83) .. (187.54,10950.28) .. controls (187.53,10945.72) and (191.22,10942.02) .. (195.77,10942.01) -- cycle ;
%Shape: Circle [id:dp055541942906038866] 
\draw   (215.77,10876.68) .. controls (220.33,10876.67) and (224.03,10880.36) .. (224.04,10884.91) .. controls (224.04,10889.47) and (220.36,10893.17) .. (215.8,10893.18) .. controls (211.24,10893.19) and (207.54,10889.5) .. (207.54,10884.94) .. controls (207.53,10880.39) and (211.22,10876.69) .. (215.77,10876.68) -- cycle ;
%Curve Lines [id:da4282879255912724] 
\draw    (109.44,10876.68) .. controls (112.86,10865.79) and (129.12,10847.38) .. (150.5,10845.14) ;
\draw [shift={(153.2,10844.94)}, rotate = 537.75] [fill={rgb, 255:red, 0; green, 0; blue, 0 }  ][line width=0.08]  [draw opacity=0] (8.04,-3.86) -- (0,0) -- (8.04,3.86) -- cycle    ;
%Curve Lines [id:da8615543540169777] 
\draw  [dash pattern={on 4.5pt off 4.5pt}]  (154.54,10850.28) .. controls (144.79,10869.11) and (135.23,10875.43) .. (116.72,10879.66) ;
\draw [shift={(114.04,10880.25)}, rotate = 348.18] [fill={rgb, 255:red, 0; green, 0; blue, 0 }  ][line width=0.08]  [draw opacity=0] (8.04,-3.86) -- (0,0) -- (8.04,3.86) -- cycle    ;
%Curve Lines [id:da717867616063592] 
\draw    (133.67,10945.33) .. controls (140.99,10939.18) and (166.82,10933.47) .. (186.56,10944.53) ;
\draw [shift={(189,10946)}, rotate = 213.02] [fill={rgb, 255:red, 0; green, 0; blue, 0 }  ][line width=0.08]  [draw opacity=0] (8.04,-3.86) -- (0,0) -- (8.04,3.86) -- cycle    ;
%Curve Lines [id:da6911170313355888] 
\draw  [dash pattern={on 4.5pt off 4.5pt}]  (187.87,10952.94) .. controls (177.65,10961.11) and (155.3,10964.22) .. (136.36,10955.55) ;
\draw [shift={(133.7,10954.25)}, rotate = 387.78] [fill={rgb, 255:red, 0; green, 0; blue, 0 }  ][line width=0.08]  [draw opacity=0] (8.04,-3.86) -- (0,0) -- (8.04,3.86) -- cycle    ;
%Straight Lines [id:da06336427809492418] 
\draw    (121.68,10939.78) -- (109.47,10893.18) ;
\draw [shift={(122.44,10942.68)}, rotate = 255.32] [fill={rgb, 255:red, 0; green, 0; blue, 0 }  ][line width=0.08]  [draw opacity=0] (8.04,-3.86) -- (0,0) -- (8.04,3.86) -- cycle    ;
%Straight Lines [id:da35258626685919126] 
\draw    (212.92,10896.05) -- (198.83,10942.17) ;
\draw [shift={(213.8,10893.18)}, rotate = 106.99] [fill={rgb, 255:red, 0; green, 0; blue, 0 }  ][line width=0.08]  [draw opacity=0] (8.04,-3.86) -- (0,0) -- (8.04,3.86) -- cycle    ;
%Straight Lines [id:da31631661695739943] 
\draw    (171.75,10845.74) -- (212.33,10877) ;
\draw [shift={(169.37,10843.91)}, rotate = 37.6] [fill={rgb, 255:red, 0; green, 0; blue, 0 }  ][line width=0.08]  [draw opacity=0] (8.04,-3.86) -- (0,0) -- (8.04,3.86) -- cycle    ;
%Curve Lines [id:da24784531714244107] 
\draw [color={rgb, 255:red, 0; green, 0; blue, 0 }  ,draw opacity=1 ]   (99.31,10884.12) .. controls (73.3,10874.79) and (92.28,10849.1) .. (101.6,10875.72) ;
\draw [shift={(102.44,10878.35)}, rotate = 253.84] [fill={rgb, 255:red, 0; green, 0; blue, 0 }  ,fill opacity=1 ][line width=0.08]  [draw opacity=0] (8.04,-3.86) -- (0,0) -- (8.04,3.86) -- cycle    ;
%Curve Lines [id:da7044706176110518] 
\draw [color={rgb, 255:red, 0; green, 0; blue, 0 }  ,draw opacity=1 ]   (125.47,10958.18) .. controls (116.38,10985.21) and (97.89,10966.32) .. (116.04,10955.18) ;
\draw [shift={(118.5,10953.83)}, rotate = 513.95] [fill={rgb, 255:red, 0; green, 0; blue, 0 }  ,fill opacity=1 ][line width=0.08]  [draw opacity=0] (8.04,-3.86) -- (0,0) -- (8.04,3.86) -- cycle    ;
%Shape: Circle [id:dp6743463461959192] 
\draw   (291.44,10876.68) .. controls (296,10876.67) and (299.7,10880.36) .. (299.7,10884.91) .. controls (299.71,10889.47) and (296.02,10893.17) .. (291.47,10893.18) .. controls (286.91,10893.19) and (283.21,10889.5) .. (283.2,10884.94) .. controls (283.2,10880.39) and (286.88,10876.69) .. (291.44,10876.68) -- cycle ;
%Shape: Circle [id:dp2577956137378201] 
\draw   (345.11,10836.68) .. controls (349.66,10836.67) and (353.36,10840.36) .. (353.37,10844.91) .. controls (353.38,10849.47) and (349.69,10853.17) .. (345.13,10853.18) .. controls (340.58,10853.19) and (336.88,10849.5) .. (336.87,10844.94) .. controls (336.86,10840.39) and (340.55,10836.69) .. (345.11,10836.68) -- cycle ;
%Shape: Circle [id:dp27614303488766057] 
\draw   (310.44,10941.68) .. controls (315,10941.67) and (318.7,10945.36) .. (318.7,10949.91) .. controls (318.71,10954.47) and (315.02,10958.17) .. (310.47,10958.18) .. controls (305.91,10958.19) and (302.21,10954.5) .. (302.2,10949.94) .. controls (302.2,10945.39) and (305.88,10941.69) .. (310.44,10941.68) -- cycle ;
%Shape: Circle [id:dp6663789366582977] 
\draw   (379.77,10942.01) .. controls (384.33,10942) and (388.03,10945.69) .. (388.04,10950.25) .. controls (388.04,10954.8) and (384.36,10958.5) .. (379.8,10958.51) .. controls (375.24,10958.52) and (371.54,10954.83) .. (371.54,10950.28) .. controls (371.53,10945.72) and (375.22,10942.02) .. (379.77,10942.01) -- cycle ;
%Shape: Circle [id:dp732468526646646] 
\draw   (399.77,10876.68) .. controls (404.33,10876.67) and (408.03,10880.36) .. (408.04,10884.91) .. controls (408.04,10889.47) and (404.36,10893.17) .. (399.8,10893.18) .. controls (395.24,10893.19) and (391.54,10889.5) .. (391.54,10884.94) .. controls (391.53,10880.39) and (395.22,10876.69) .. (399.77,10876.68) -- cycle ;
%Curve Lines [id:da5830823799190792] 
\draw    (293.44,10876.68) .. controls (296.86,10865.79) and (313.12,10847.38) .. (334.5,10845.14) ;
\draw [shift={(337.2,10844.94)}, rotate = 537.75] [fill={rgb, 255:red, 0; green, 0; blue, 0 }  ][line width=0.08]  [draw opacity=0] (8.04,-3.86) -- (0,0) -- (8.04,3.86) -- cycle    ;
%Curve Lines [id:da152355913900855] 
\draw  [dash pattern={on 4.5pt off 4.5pt}]  (338.54,10850.28) .. controls (328.79,10869.11) and (319.23,10875.43) .. (300.72,10879.66) ;
\draw [shift={(298.04,10880.25)}, rotate = 348.18] [fill={rgb, 255:red, 0; green, 0; blue, 0 }  ][line width=0.08]  [draw opacity=0] (8.04,-3.86) -- (0,0) -- (8.04,3.86) -- cycle    ;
%Curve Lines [id:da14354319450519237] 
\draw    (317.67,10945.33) .. controls (324.99,10939.18) and (350.82,10933.47) .. (370.56,10944.53) ;
\draw [shift={(373,10946)}, rotate = 213.02] [fill={rgb, 255:red, 0; green, 0; blue, 0 }  ][line width=0.08]  [draw opacity=0] (8.04,-3.86) -- (0,0) -- (8.04,3.86) -- cycle    ;
%Curve Lines [id:da40249733846191393] 
\draw  [dash pattern={on 4.5pt off 4.5pt}]  (371.87,10952.94) .. controls (361.65,10961.11) and (339.3,10964.22) .. (320.36,10955.55) ;
\draw [shift={(317.7,10954.25)}, rotate = 387.78] [fill={rgb, 255:red, 0; green, 0; blue, 0 }  ][line width=0.08]  [draw opacity=0] (8.04,-3.86) -- (0,0) -- (8.04,3.86) -- cycle    ;
%Straight Lines [id:da2251729341549038] 
\draw    (305.68,10939.78) -- (293.47,10893.18) ;
\draw [shift={(306.44,10942.68)}, rotate = 255.32] [fill={rgb, 255:red, 0; green, 0; blue, 0 }  ][line width=0.08]  [draw opacity=0] (8.04,-3.86) -- (0,0) -- (8.04,3.86) -- cycle    ;
%Straight Lines [id:da9527411075687711] 
\draw    (396.92,10896.05) -- (382.83,10942.17) ;
\draw [shift={(397.8,10893.18)}, rotate = 106.99] [fill={rgb, 255:red, 0; green, 0; blue, 0 }  ][line width=0.08]  [draw opacity=0] (8.04,-3.86) -- (0,0) -- (8.04,3.86) -- cycle    ;
%Straight Lines [id:da9855354083962371] 
\draw    (355.75,10845.74) -- (396.33,10877) ;
\draw [shift={(353.37,10843.91)}, rotate = 37.6] [fill={rgb, 255:red, 0; green, 0; blue, 0 }  ][line width=0.08]  [draw opacity=0] (8.04,-3.86) -- (0,0) -- (8.04,3.86) -- cycle    ;
%Curve Lines [id:da288993629775407] 
\draw [color={rgb, 255:red, 0; green, 0; blue, 0 }  ,draw opacity=1 ]   (340.11,10837.68) .. controls (335.92,10809.05) and (352.54,10809.91) .. (349.21,10834.57) ;
\draw [shift={(348.77,10837.35)}, rotate = 280.37] [fill={rgb, 255:red, 0; green, 0; blue, 0 }  ,fill opacity=1 ][line width=0.08]  [draw opacity=0] (8.04,-3.86) -- (0,0) -- (8.04,3.86) -- cycle    ;
%Curve Lines [id:da018669091835059826] 
\draw [color={rgb, 255:red, 0; green, 0; blue, 0 }  ,draw opacity=1 ] [dash pattern={on 4.5pt off 4.5pt}]  (283.31,10884.12) .. controls (256.33,10871.88) and (277.17,10848.93) .. (285.69,10875.72) ;
\draw [shift={(286.44,10878.35)}, rotate = 255.57] [fill={rgb, 255:red, 0; green, 0; blue, 0 }  ,fill opacity=1 ][line width=0.08]  [draw opacity=0] (8.04,-3.86) -- (0,0) -- (8.04,3.86) -- cycle    ;
%Curve Lines [id:da6338026446697331] 
\draw [color={rgb, 255:red, 0; green, 0; blue, 0 }  ,draw opacity=1 ]   (402.77,10876.68) .. controls (406.55,10854.03) and (430.81,10863.39) .. (410.31,10881.02) ;
\draw [shift={(408.27,10882.68)}, rotate = 322.18] [fill={rgb, 255:red, 0; green, 0; blue, 0 }  ,fill opacity=1 ][line width=0.08]  [draw opacity=0] (8.04,-3.86) -- (0,0) -- (8.04,3.86) -- cycle    ;
%Curve Lines [id:da43395171014840805] 
\draw [color={rgb, 255:red, 0; green, 0; blue, 0 }  ,draw opacity=1 ] [dash pattern={on 4.5pt off 4.5pt}]  (309.47,10958.18) .. controls (300.38,10985.21) and (281.89,10966.32) .. (300.04,10955.18) ;
\draw [shift={(302.5,10953.83)}, rotate = 513.95] [fill={rgb, 255:red, 0; green, 0; blue, 0 }  ,fill opacity=1 ][line width=0.08]  [draw opacity=0] (8.04,-3.86) -- (0,0) -- (8.04,3.86) -- cycle    ;
%Curve Lines [id:da14000498433720399] 
\draw [color={rgb, 255:red, 0; green, 0; blue, 0 }  ,draw opacity=1 ][line width=0.75]    (386.8,10954.51) .. controls (408.22,10966.24) and (386.62,10985.28) .. (382.21,10961.29) ;
\draw [shift={(381.8,10958.51)}, rotate = 443.52] [fill={rgb, 255:red, 0; green, 0; blue, 0 }  ,fill opacity=1 ][line width=0.08]  [draw opacity=0] (8.04,-3.86) -- (0,0) -- (8.04,3.86) -- cycle    ;
%Shape: Circle [id:dp9158017219487349] 
\draw   (196.44,11066.68) .. controls (201,11066.67) and (204.7,11070.36) .. (204.7,11074.91) .. controls (204.71,11079.47) and (201.02,11083.17) .. (196.47,11083.18) .. controls (191.91,11083.19) and (188.21,11079.5) .. (188.2,11074.94) .. controls (188.2,11070.39) and (191.88,11066.69) .. (196.44,11066.68) -- cycle ;
%Shape: Circle [id:dp3434445403058337] 
\draw   (250.11,11026.68) .. controls (254.66,11026.67) and (258.36,11030.36) .. (258.37,11034.91) .. controls (258.38,11039.47) and (254.69,11043.17) .. (250.13,11043.18) .. controls (245.58,11043.19) and (241.88,11039.5) .. (241.87,11034.94) .. controls (241.86,11030.39) and (245.55,11026.69) .. (250.11,11026.68) -- cycle ;
%Shape: Circle [id:dp6032993812559229] 
\draw   (215.44,11131.68) .. controls (220,11131.67) and (223.7,11135.36) .. (223.7,11139.91) .. controls (223.71,11144.47) and (220.02,11148.17) .. (215.47,11148.18) .. controls (210.91,11148.19) and (207.21,11144.5) .. (207.2,11139.94) .. controls (207.2,11135.39) and (210.88,11131.69) .. (215.44,11131.68) -- cycle ;
%Shape: Circle [id:dp7757208444541868] 
\draw   (284.77,11132.01) .. controls (289.33,11132) and (293.03,11135.69) .. (293.04,11140.25) .. controls (293.04,11144.8) and (289.36,11148.5) .. (284.8,11148.51) .. controls (280.24,11148.52) and (276.54,11144.83) .. (276.54,11140.28) .. controls (276.53,11135.72) and (280.22,11132.02) .. (284.77,11132.01) -- cycle ;
%Shape: Circle [id:dp2993509845220661] 
\draw   (304.77,11066.68) .. controls (309.33,11066.67) and (313.03,11070.36) .. (313.04,11074.91) .. controls (313.04,11079.47) and (309.36,11083.17) .. (304.8,11083.18) .. controls (300.24,11083.19) and (296.54,11079.5) .. (296.54,11074.94) .. controls (296.53,11070.39) and (300.22,11066.69) .. (304.77,11066.68) -- cycle ;
%Curve Lines [id:da46152090843246873] 
\draw    (198.44,11066.68) .. controls (201.86,11055.79) and (218.12,11037.38) .. (239.5,11035.14) ;
\draw [shift={(242.2,11034.94)}, rotate = 537.75] [fill={rgb, 255:red, 0; green, 0; blue, 0 }  ][line width=0.08]  [draw opacity=0] (8.04,-3.86) -- (0,0) -- (8.04,3.86) -- cycle    ;
%Curve Lines [id:da4756584480649195] 
\draw  [dash pattern={on 4.5pt off 4.5pt}]  (243.54,11040.28) .. controls (233.79,11059.11) and (224.23,11065.43) .. (205.72,11069.66) ;
\draw [shift={(203.04,11070.25)}, rotate = 348.18] [fill={rgb, 255:red, 0; green, 0; blue, 0 }  ][line width=0.08]  [draw opacity=0] (8.04,-3.86) -- (0,0) -- (8.04,3.86) -- cycle    ;
%Curve Lines [id:da027305368857017953] 
\draw    (222.67,11135.33) .. controls (229.99,11129.18) and (255.82,11123.47) .. (275.56,11134.53) ;
\draw [shift={(278,11136)}, rotate = 213.02] [fill={rgb, 255:red, 0; green, 0; blue, 0 }  ][line width=0.08]  [draw opacity=0] (8.04,-3.86) -- (0,0) -- (8.04,3.86) -- cycle    ;
%Curve Lines [id:da26531293840173076] 
\draw  [dash pattern={on 4.5pt off 4.5pt}]  (276.87,11142.94) .. controls (266.65,11151.11) and (244.3,11154.22) .. (225.36,11145.55) ;
\draw [shift={(222.7,11144.25)}, rotate = 387.78] [fill={rgb, 255:red, 0; green, 0; blue, 0 }  ][line width=0.08]  [draw opacity=0] (8.04,-3.86) -- (0,0) -- (8.04,3.86) -- cycle    ;
%Straight Lines [id:da009324107067658849] 
\draw    (210.68,11129.78) -- (198.47,11083.18) ;
\draw [shift={(211.44,11132.68)}, rotate = 255.32] [fill={rgb, 255:red, 0; green, 0; blue, 0 }  ][line width=0.08]  [draw opacity=0] (8.04,-3.86) -- (0,0) -- (8.04,3.86) -- cycle    ;
%Straight Lines [id:da6515140418216843] 
\draw    (301.92,11086.05) -- (287.83,11132.17) ;
\draw [shift={(302.8,11083.18)}, rotate = 106.99] [fill={rgb, 255:red, 0; green, 0; blue, 0 }  ][line width=0.08]  [draw opacity=0] (8.04,-3.86) -- (0,0) -- (8.04,3.86) -- cycle    ;
%Straight Lines [id:da5675506842705809] 
\draw    (260.75,11035.74) -- (301.33,11067) ;
\draw [shift={(258.37,11033.91)}, rotate = 37.6] [fill={rgb, 255:red, 0; green, 0; blue, 0 }  ][line width=0.08]  [draw opacity=0] (8.04,-3.86) -- (0,0) -- (8.04,3.86) -- cycle    ;
%Curve Lines [id:da3048598140190868] 
\draw [color={rgb, 255:red, 0; green, 0; blue, 0 }  ,draw opacity=1 ] [dash pattern={on 4.5pt off 4.5pt}]  (245.11,11027.68) .. controls (240.92,10999.05) and (257.54,10999.91) .. (254.21,11024.57) ;
\draw [shift={(253.77,11027.35)}, rotate = 280.37] [fill={rgb, 255:red, 0; green, 0; blue, 0 }  ,fill opacity=1 ][line width=0.08]  [draw opacity=0] (8.04,-3.86) -- (0,0) -- (8.04,3.86) -- cycle    ;
%Curve Lines [id:da4991315158680465] 
\draw [color={rgb, 255:red, 0; green, 0; blue, 0 }  ,draw opacity=1 ] [dash pattern={on 4.5pt off 4.5pt}]  (188.31,11074.12) .. controls (162.3,11061.88) and (180.34,11040.81) .. (190.52,11065.88) ;
\draw [shift={(191.44,11068.35)}, rotate = 250.99] [fill={rgb, 255:red, 0; green, 0; blue, 0 }  ,fill opacity=1 ][line width=0.08]  [draw opacity=0] (8.04,-3.86) -- (0,0) -- (8.04,3.86) -- cycle    ;
%Curve Lines [id:da8332296481872739] 
\draw [color={rgb, 255:red, 0; green, 0; blue, 0 }  ,draw opacity=1 ] [dash pattern={on 4.5pt off 4.5pt}]  (307.77,11066.68) .. controls (311.55,11044.03) and (335.81,11053.39) .. (315.31,11071.02) ;
\draw [shift={(313.27,11072.68)}, rotate = 322.18] [fill={rgb, 255:red, 0; green, 0; blue, 0 }  ,fill opacity=1 ][line width=0.08]  [draw opacity=0] (8.04,-3.86) -- (0,0) -- (8.04,3.86) -- cycle    ;
%Curve Lines [id:da5354247258952805] 
\draw [color={rgb, 255:red, 0; green, 0; blue, 0 }  ,draw opacity=1 ] [dash pattern={on 4.5pt off 4.5pt}]  (214.47,11148.18) .. controls (205.38,11175.21) and (186.89,11156.32) .. (205.04,11145.18) ;
\draw [shift={(207.5,11143.83)}, rotate = 513.95] [fill={rgb, 255:red, 0; green, 0; blue, 0 }  ,fill opacity=1 ][line width=0.08]  [draw opacity=0] (8.04,-3.86) -- (0,0) -- (8.04,3.86) -- cycle    ;
%Curve Lines [id:da9772255003705139] 
\draw [color={rgb, 255:red, 0; green, 0; blue, 0 }  ,draw opacity=1 ][line width=0.75]  [dash pattern={on 4.5pt off 4.5pt}]  (291.8,11144.51) .. controls (313.22,11156.24) and (291.62,11175.28) .. (287.21,11151.29) ;
\draw [shift={(286.8,11148.51)}, rotate = 443.52] [fill={rgb, 255:red, 0; green, 0; blue, 0 }  ,fill opacity=1 ][line width=0.08]  [draw opacity=0] (8.04,-3.86) -- (0,0) -- (8.04,3.86) -- cycle    ;

% Text Node
\draw (108.45,10884.93) node  [font=\tiny,rotate=-1.27] [align=left] {$\displaystyle 1$};
% Text Node
\draw (162.12,10844.93) node  [font=\tiny,rotate=-1.27] [align=left] {$\displaystyle 5$};
% Text Node
\draw (127.45,10949.93) node  [font=\tiny,rotate=-1.27] [align=left] {2};
% Text Node
\draw (196.79,10950.26) node  [font=\tiny,rotate=-1.27] [align=left] {3};
% Text Node
\draw (216.79,10884.93) node  [font=\tiny,rotate=-1.27] [align=left] {$\displaystyle 4$};
% Text Node
\draw (292.45,10884.93) node  [font=\tiny,rotate=-1.27] [align=left] {$\displaystyle 1$};
% Text Node
\draw (346.12,10844.93) node  [font=\tiny,rotate=-1.27] [align=left] {$\displaystyle 5$};
% Text Node
\draw (311.45,10949.93) node  [font=\tiny,rotate=-1.27] [align=left] {$\displaystyle 2$};
% Text Node
\draw (380.79,10950.26) node  [font=\tiny,rotate=-1.27] [align=left] {$\displaystyle 3$};
% Text Node
\draw (400.79,10884.93) node  [font=\tiny,rotate=-1.27] [align=left] {$\displaystyle 4$};
% Text Node
\draw (197.45,11074.93) node  [font=\tiny,rotate=-1.27] [align=left] {$\displaystyle 1$};
% Text Node
\draw (251.12,11034.93) node  [font=\tiny,rotate=-1.27] [align=left] {$\displaystyle 5$};
% Text Node
\draw (216.45,11139.93) node  [font=\tiny,rotate=-1.27] [align=left] {2};
% Text Node
\draw (285.79,11140.26) node  [font=\tiny,rotate=-1.27] [align=left] {3};
% Text Node
\draw (305.79,11074.93) node  [font=\tiny,rotate=-1.27] [align=left] {$\displaystyle 4$};
% Text Node
\draw (161,11002.33) node   [align=left] {\begin{minipage}[lt]{15.64pt}\setlength\topsep{0pt}
(a)
\end{minipage}};
% Text Node
\draw (346,11003) node   [align=left] {\begin{minipage}[lt]{15.64pt}\setlength\topsep{0pt}
(b)
\end{minipage}};
% Text Node
\draw (251,11193) node   [align=left] {\begin{minipage}[lt]{15.64pt}\setlength\topsep{0pt}
(c)
\end{minipage}};

\end{tikzpicture}

%%%%%%%%%%5
\caption{(a) Global graph ${\mathcal{G}}$ for an N1DS, which is the same as its corresponding $\Delta$-network graph for $\Delta=\{0\}$; (b) associated $\Delta$-network graph for $\Delta=\mathbb{C}\setminus\{0\}$; (c) associated $\Delta$-network graph for $\Delta=\mathbb{C}$.}
\label{stable}
\end{figure}
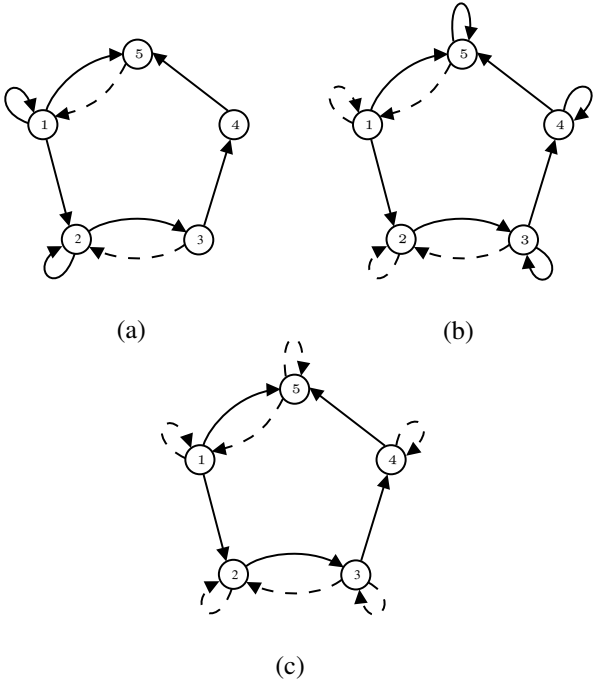

\begin{ex}
Consider an N1DS with dynamics (\ref{e1}), whose corresponding global graph is depicted in Fig.~\ref{stable}(a), and assume that the set of system matrices $\mathcal{S}^*$ is the same as $\mathcal{P}(\mathcal{A})$.  First, let $\Delta=\{0\}$. As discussed in Example \ref{ex5}, one can see that the $\Delta$-network graph is the same as the global graph. Now, let  $\Delta=\mathbb{C}\setminus\{0\}$. In this case, as explained in Example \ref{ex6}, for $k=1,2,\ldots, n$, $f^{\mathcal{L}}_{\Delta}(k)=*$ if $\mathcal{A}(k,k)=0$, and $f^{\mathcal{L}}_{\Delta}(k)=?$ otherwise. Thus, the corresponding $\Delta$-network graph is the one shown in Fig.~\ref{stable}(b). Finally, one can see that, for $\Delta=\mathbb{C}$ (as in Example \ref{ex7}), the graph shown in Fig.~\ref{stable}(c) is the associated $\Delta$-network graph.
\label{ex10}\end{ex}

\subsection{Coloring Process and Zero Forcing Sets}
 
Our combinatorial characterization of modal strong structural controllability
involves a  coloring process on the network. Let $G=(V,E)$ be a graph with nodes $Z\subset V$ colored black and other nodes as white. Let $E=E^*\cup E^?$. Now, consider the following ``coloring rule'':
\begin{itemize}
\item[] $u$ is colored black if it is the only white strong out-neighbor of some node $v$ \big($(v,u)\in E^*$\big). Then, we say that $v$ forces $u$ and represent it by $v\rightarrow u$ ($u$ and $v$ may be the same node; moreover, the forcing node $v$ can be black or white).
\end{itemize}

The coloring process refers to the sequential coloring of nodes through employing the coloring rule as many times as possible (until no further coloring rule can be applied). For more information on  this coloring process, the reader is referred to \cite{shima_ECC, barioli2009minimum,trefois2015zero}.

The set of black nodes after the termination of the coloring process is referred to as the {\textit{derived set}} of $Z$ and is denoted by $\mathcal{D}(Z)$. Furthermore, if for the initial set of black nodes $Z$, $\mathcal{D}(Z)=V$, then $Z$ is called a \textit{zero forcing set (ZFS)}.

{We also note that checking for whether an initial set of black nodes is a zero forcing set or not can be done in a polynomial time; however, finding a zero forcing set with the smallest cardinality is NP-hard \cite{ trefois2015zero}.} 

\begin{ex}
Consider the graph in Fig.~\ref{stable}(a), and let $Z=\{1\}$ be the set of  initial black nodes.  
The steps of the coloring
process are: (a) $4 \rightarrow 5$, (b) $1\rightarrow 2$, (c) $2 \rightarrow3$, and $3 \rightarrow 4$. Moreover, the coloring process and its successive steps in Fig.~\ref{stable}(b) are as follows: 
(a) $5 \rightarrow 5$, (b) $4\rightarrow 4$, and $1 \rightarrow 2$, (c) $2 \rightarrow 3$.
Thus, $Z=\{1\}$ is a  zero forcing set for both graphs in Figs.~\ref{stable}(a) and (b), while one can see that the coloring process cannot be initiated in Fig.~\ref{stable}(c). 
\label{ex11}
\end{ex}

\section{Main Results}\label{results}

In this section, we discuss the main results of this paper. 

\subsection{Modal Strong Structural Controllability}

Let $\Delta\subseteq \mathbb{C}$ be a given nonempty set. Now, consider an LTI network $\mathcal{N}$, with dynamics (\ref{e1}) and the pattern matrix $\mathcal{A}\in\{0,*,?\}^{N\times N}$. Let $f^{\mathcal{L}}_{\Delta}\in\{0,*,?\}^{n}$ be a given $\Delta$-characteristic vector. Accordingly, the corresponding $\Delta$-network graph $G_{\mathcal{N}}^{\Delta}=(V_{\mathcal{N}}, E_{\mathcal{N}})$, with $V_{\mathcal{N}}=\{1,2,\ldots, n\}$, and $E_{\mathcal{N}}=E_{\mathcal{N}}^?\cup E_{\mathcal{N}}^*$, can be obtained as discussed in Section \ref{section}. We recall that every node $i\in V_{\mathcal{N}}$ is a super node, which represents subsystem $i$, and $G_i=(V_i,E_i)$ is the $i$-th node graph, where $V_i=\{i^1,\ldots, i^{l_i}\}$ is called the set of vertices of node $i$. For a subset of nodes $Z\subseteq V_{\mathcal{N}}$, $\mbox{Ver}(Z)$ is the set of all vertices associated with the nodes in $Z$. 

We first study the controllability of any eigenvalue $\lambda\in \Delta$ for a family of LTI networks with
$A\in\mathcal{S}(f^{\mathcal{L}}_{\Delta}, \mathcal{A})$. 
Let $Z\subseteq V_{\mathcal{N}} $ be a subset of nodes of the $\Delta$-network graph. Moreover, $\mathcal{D}(Z)$ is the derived set of $Z$ in the $\Delta$-network graph.      In this direction, we first show that if a left eigenvector $\nu\in\mathbb{R}^N$ associated with  an eigenvalue $\lambda\in\Delta$ for some $A\in\mathcal{S}(f^{\mathcal{L}}_{\Delta}, \mathcal{A})$ vanishes at entries indexed by all the vertices of $Z$, this eigenvector has to vanish at every entry indexed by the vertices of  $\mathcal{D}(Z)$.

\begin{lem}
Let $Z\subset V_{\mathcal{N}}$, and let  
$A\in\mathcal{S}(f^{\mathcal{L}}_{\Delta}, \mathcal{A})$. Assume that $\nu\in \mathbb{C}^N$ be a left eigenvector of $A$ associated with some $\lambda\in\Delta$. If $\nu(k)=0$ for all $k\in \mbox{Ver}(Z)$, then $\nu(k)=0$ for all $k\in \mbox{Ver}(\mathcal{D}(Z))$.
\label{l2}
\end{lem}
  
\textbf{Proof.} If $Z=\mathcal{D}(Z)$, the statement of the lemma follows. Otherwise, the coloring process can be applied to the $\Delta$-network graph. Thus, we have one of the following cases: (a) There exists a black node $i\in V_{\mathcal{N}}$  with exactly one white out-neighbor $j$, where $(i,j)\in E_{\mathcal{N}}^*$; (b) there exists a white node $i$ which has no white out-neighbor except itself, and $(i,i)\in E_{\mathcal{N}}^*$; and (c) there exists a white node $i$ with no self-loop (i.e., $(i,i)\notin E_{\mathcal{N}}$) and exactly one white out-neighbor $j$, where $(i,j)\in E_{\mathcal{N}}^*$.

Consider matrix $A$ as a block matrix described in (\ref{sm}), and let $\nu^T=\begin{bmatrix}\nu_1^T & \ldots &\nu_n^T\end{bmatrix}$, where, for $i=1,2,\ldots, n$,  $\nu_i=\begin{bmatrix}\nu_i^1 & \ldots &\nu_i^{l_i}\end{bmatrix}^T\in\mathbb{C}^{l_i}$. Now, one can write the matrix equation $\nu^TA=\lambda\nu^T$ as:
\begin{equation}\begin{aligned} &\colvec[1]{\nu_1^T & \nu_2^T & \ldots & \nu_n^T
  }\colvec[.9]{\lambda I_{l_1}-A_{11} & -A_{12} & \ldots & -A_{1n}\\
  -A_{21} & \lambda I_{l_2}-A_{22} & \ldots & -A_{2n}\\
  \ddots & \ddots & \ddots & \ddots\\
  -A_{n1} & -A_{n2} & \ldots & \lambda I_{l_n}-A_{nn}
  }\\&=0.\end{aligned}
\label{A2_2}\end{equation}
  
For the first case, the $i$-th block column of equation  (\ref{A2_2}) assumes the form, 
   \begin{equation}
      {\nu_i^T(\lambda I_{l_i}-A_{ii})+\nu_j^T A_{ji}+\sum_{k\in N_{\rm out}(i), k\neq i,j} \nu_k^T A_{ki}=0.}
 \label{strong}  
   \end{equation}
Since node $i$ is black, $i$ belongs to $Z$, and $\nu(k)=0$, for all $ k \in \mbox{Ver}(\{i\})$. Thus, we have $\nu_i=0$. Likewise, all the out-neighbors of node $i$ except  $j$ are black, implying that $\nu(k)=0$, for all $ k \in \mbox{Ver}(N_{\rm out}(i)\setminus \{j\})$. Thus, (\ref{strong}) reduces to $\nu_j^TA_{ji}=0$.
Since $(i,j)\in E_{\mathcal{N}}^*$, the submatrix $A_{ji}$ has full row rank, and hence we have that $\nu_j=0$.

Now, consider case (b). Note that in the $\Delta$-network graph $G_{\mathcal{N}}^{\Delta}$, $(i,i)\in E_{\mathcal{N}}^*$ implies that $f^{\mathcal{L}}_{\Delta}(i)=*$. Thus, from (\ref{f}), $\Lambda(A_{ii})\cap\Delta=\emptyset$, and since $\lambda\in \Delta$, $\lambda\notin \Lambda(A_{ii})$. Thus,  one can conclude that  the matrix $\lambda I_{l_i}-A_{ii}$ is nonsingular. Moreover, we have,
 \begin{equation}
      {\nu_i^T(\lambda I_{l_i}-A_{ii})+\sum_{k\in N_{\rm out}(i), k\neq i} \nu_k^T A_{ki}=0.}
 \label{strong2}  
   \end{equation}
Since all out-neighbors of node $i$ except itself are black, $\nu_k=0$, for all $ k \in N_{\rm out}(i), k\neq i$. Thus, (\ref{strong2}) reduces to $\nu_i^T(\lambda I_{l_i}-A_{ii})=0$. Now, since  $\lambda I_{l_i}-A_{ii}$ is nonsingular, we have $\nu_i=0$.

Finally, assume that case (c) is valid. Since $(i,i)\notin E_{\mathcal{N}}$, we have $f^{\mathcal{L}}_{\Delta}(i)=0$. Thus, from (\ref{f}), one can see that $\Delta=\{\lambda\}$ and $A_{ii}=\lambda I_{l_i}$. Accordingly, (\ref{strong}) reduces to, 
   \begin{equation}
      {\nu_j^T A_{ji}+\sum_{k\in N_{\rm out}(i), k\neq j} \nu_k^T A_{ki}=0.}
 \label{sss}  
   \end{equation}
Again, since all out-neighbors of $i$ except $j$ are black, $\nu_k=0$, for all $ k \in N_{\rm out}(i), k\neq j$. Hence, we have  $\nu_j^T A_{ji}=0$. In addition, based on the defnition of a $\Delta$-network graph,  $(i,j)\in E_{\mathcal{N}}^*$ implies that $A_{ji}$ has full row rank.  Thereby, one has $\nu_j=0$. Therefore, if node $j$  becomes black during the coloring process, the subvector of $\nu$ corresponding to the vertices of $j$, that is,   $\nu_j$ should be zero. Hence, by the termination of the  coloring process, we have $\nu_k=0$, for all $ k \in \mbox{Ver}(\mathcal{D}(Z))$.   %
\carre

The next theorem, regarding the modal strong structural controllability, is one of the main results of this work.

\begin{thm}
Given some $\Delta\subseteq \mathbb{C}$, an LTI network with dynamics (\ref{e1}) and  the $\Delta$-network graph $G_{\mathcal{N}}^{\Delta}=(V_{\mathcal{N}}, E_{\mathcal{N}})$ is $\Delta$-SSC if $V_C$ is a ZFS of $G_{\mathcal{N}}^{\Delta}$.
\label{th1}
\end{thm}

\textbf{Proof.} Suppose that $V_C$ is a ZFS of $G_{\mathcal{N}}^{\Delta}$, but there exists  $A\in\mathcal{S}(f^{\mathcal{L}}_{\Delta}, \mathcal{A})$ 
with some uncontrollable  eigenvalue  in $\Delta$. 
Then, $A$ has a nonzero left eigenvector $\nu=\begin{bmatrix}\nu_1^T & \ldots &\nu_n^T\end{bmatrix}^T\in\mathbb{C}^N$ associated with $\lambda\in \Delta$ such that $\nu^TB=0$, or equivalently, $\nu(i)=0$, for all $ i \in \mbox{Ver}(V_C)$. However, since $\mathcal{D}(V_C)=V_{\mathcal{N}}$, it follows from Lemma \ref{l2} that $\nu=0$; which is a contradiction. 
\carre  
 
\begin{ex}
Consider Example \ref{ex_stab0}, where the associated $\Delta$-network graph $G_{\mathcal{N}}^{\Delta}$ is shown in Fig.~\ref{globgraph_ex}. In  this example, one can see that  every set $V_C=\{i\}$, $i=1,2,\ldots,n$, is a ZFS of $G_{\mathcal{N}}^{\Delta}$, and thus renders the network $\Delta$-SSC.  
\end{ex}

Note that Theorem \ref{th1} provides a sufficient condition for modal strong structural controllability of an LTI network that includes $n$ subsystems with (probably) different dimensions. In the following, given a set $\Delta\subseteq \mathbb{C}$ that includes at least one real number, we show that if the LTI network is an N1DS (i.e.,  all its subsystems are single-state), then the controllability condition is necessary as well.

\begin{thm}
Consider some $\Delta\subseteq \mathbb{C}$, where $\Delta\cap \mathbb{R}\neq \emptyset$. Then, an N1DS with dynamics (\ref{e1}) and  the $\Delta$-network graph $G_{\mathcal{N}}^{\Delta}=(V_{\mathcal{N}}, E_{\mathcal{N}})$ is $\Delta$-SSC if and only if $V_C$ is a ZFS of $G_{\mathcal{N}}^{\Delta}$.
\label{th5}
\end{thm}

\textbf{Proof.} The sufficiency is proved by  considering Theorem \ref{th1}.  
To prove the necessity, assume that every eigenvalue $\lambda\in\Delta$ of all $A\in\mathcal{S}(f^{\mathcal{L}}_{\Delta}, \mathcal{A})$ is controllable, but $V_C$ is not a ZFS of $G_{\mathcal{N}}^{\Delta}$, that is, $\mathcal{D}(V_C)\neq V_{\mathcal{N}}$. In the $\Delta$-network graph, let $\mathcal{B}= \mathcal{D}(V_C)$, and $\mathcal{W}=V_{\mathcal{N}}\setminus \mathcal{D}(V_C)$. Without loss of generality, index the nodes of $\mathcal{B}$ first. Let $r=\mid\mathcal{B}\mid$.
Then, the pattern matrix $\mathcal{A}$ can be partitioned as,
\begin{equation*}\mathcal{A}=\begin{bmatrix}
\mathcal{A}_{1} & \mathcal{A}_2\\
\mathcal{A}_3 & \mathcal{A}_4
\end{bmatrix},\label{A22}
\end{equation*}
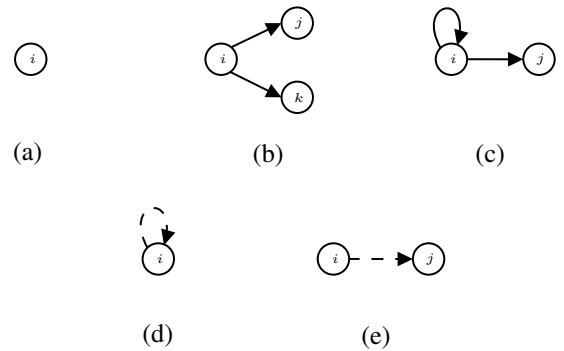
\begin{figure}[b]
\centering

%%%%%%%%%%%%%%%%%%%%%%%%%%

\tikzset{every picture/.style={line width=0.75pt}} %set default line width to 0.75pt        

\begin{tikzpicture}[x=0.75pt,y=0.75pt,yscale=-.96,xscale=.96]
%uncomment if require: \path (0,11616); %set diagram left start at 0, and has height of 11616

%Shape: Circle [id:dp5062683336554414] 
\draw   (72.44,11268.68) .. controls (77,11268.67) and (80.7,11272.36) .. (80.7,11276.91) .. controls (80.71,11281.47) and (77.02,11285.17) .. (72.47,11285.18) .. controls (67.91,11285.19) and (64.21,11281.5) .. (64.2,11276.94) .. controls (64.2,11272.39) and (67.88,11268.69) .. (72.44,11268.68) -- cycle ;
%Shape: Circle [id:dp20779620509951546] 
\draw   (172.44,11268.68) .. controls (177,11268.67) and (180.7,11272.36) .. (180.7,11276.91) .. controls (180.71,11281.47) and (177.02,11285.17) .. (172.47,11285.18) .. controls (167.91,11285.19) and (164.21,11281.5) .. (164.2,11276.94) .. controls (164.2,11272.39) and (167.88,11268.69) .. (172.44,11268.68) -- cycle ;
%Shape: Circle [id:dp775596075759952] 
\draw   (212.44,11288.68) .. controls (217,11288.67) and (220.7,11292.36) .. (220.7,11296.91) .. controls (220.71,11301.47) and (217.02,11305.17) .. (212.47,11305.18) .. controls (207.91,11305.19) and (204.21,11301.5) .. (204.2,11296.94) .. controls (204.2,11292.39) and (207.88,11288.69) .. (212.44,11288.68) -- cycle ;
%Shape: Circle [id:dp5127928783063724] 
\draw   (212.44,11249.68) .. controls (217,11249.67) and (220.7,11253.36) .. (220.7,11257.91) .. controls (220.71,11262.47) and (217.02,11266.17) .. (212.47,11266.18) .. controls (207.91,11266.19) and (204.21,11262.5) .. (204.2,11257.94) .. controls (204.2,11253.39) and (207.88,11249.69) .. (212.44,11249.68) -- cycle ;
%Shape: Circle [id:dp2768405690225588] 
\draw   (293.44,11268.68) .. controls (298,11268.67) and (301.7,11272.36) .. (301.7,11276.91) .. controls (301.71,11281.47) and (298.02,11285.17) .. (293.47,11285.18) .. controls (288.91,11285.19) and (285.21,11281.5) .. (285.2,11276.94) .. controls (285.2,11272.39) and (288.88,11268.69) .. (293.44,11268.68) -- cycle ;
%Shape: Circle [id:dp6057445251118629] 
\draw   (339.44,11268.68) .. controls (344,11268.67) and (347.7,11272.36) .. (347.7,11276.91) .. controls (347.71,11281.47) and (344.02,11285.17) .. (339.47,11285.18) .. controls (334.91,11285.19) and (331.21,11281.5) .. (331.2,11276.94) .. controls (331.2,11272.39) and (334.88,11268.69) .. (339.44,11268.68) -- cycle ;
%Shape: Circle [id:dp5148927411600592] 
\draw   (139.44,11373.68) .. controls (144,11373.67) and (147.7,11377.36) .. (147.7,11381.91) .. controls (147.71,11386.47) and (144.02,11390.17) .. (139.47,11390.18) .. controls (134.91,11390.19) and (131.21,11386.5) .. (131.2,11381.94) .. controls (131.2,11377.39) and (134.88,11373.69) .. (139.44,11373.68) -- cycle ;
%Shape: Circle [id:dp16393127348688186] 
\draw   (231.44,11373.68) .. controls (236,11373.67) and (239.7,11377.36) .. (239.7,11381.91) .. controls (239.71,11386.47) and (236.02,11390.17) .. (231.47,11390.18) .. controls (226.91,11390.19) and (223.21,11386.5) .. (223.2,11381.94) .. controls (223.2,11377.39) and (226.88,11373.69) .. (231.44,11373.68) -- cycle ;
%Shape: Circle [id:dp6703556784111679] 
\draw   (281.44,11373.68) .. controls (286,11373.67) and (289.7,11377.36) .. (289.7,11381.91) .. controls (289.71,11386.47) and (286.02,11390.17) .. (281.47,11390.18) .. controls (276.91,11390.19) and (273.21,11386.5) .. (273.2,11381.94) .. controls (273.2,11377.39) and (276.88,11373.69) .. (281.44,11373.68) -- cycle ;
%Straight Lines [id:da8972472436047929] 
\draw    (201.5,11259.24) -- (178,11270.5) ;
\draw [shift={(204.2,11257.94)}, rotate = 154.4] [fill={rgb, 255:red, 0; green, 0; blue, 0 }  ][line width=0.08]  [draw opacity=0] (8.04,-3.86) -- (0,0) -- (8.04,3.86) -- cycle    ;
%Straight Lines [id:da0291076048127894] 
\draw    (201.51,11295.61) -- (177,11283.5) ;
\draw [shift={(204.2,11296.94)}, rotate = 206.3] [fill={rgb, 255:red, 0; green, 0; blue, 0 }  ][line width=0.08]  [draw opacity=0] (8.04,-3.86) -- (0,0) -- (8.04,3.86) -- cycle    ;
%Straight Lines [id:da18308565584564906] 
\draw    (328.2,11276.94) -- (301.7,11276.91) ;
\draw [shift={(331.2,11276.94)}, rotate = 180.06] [fill={rgb, 255:red, 0; green, 0; blue, 0 }  ][line width=0.08]  [draw opacity=0] (8.04,-3.86) -- (0,0) -- (8.04,3.86) -- cycle    ;
%Straight Lines [id:da891686893025794] 
\draw  [dash pattern={on 4.5pt off 4.5pt}]  (270.2,11381.94) -- (252,11381.93) -- (239.7,11381.91) ;
\draw [shift={(273.2,11381.94)}, rotate = 180.05] [fill={rgb, 255:red, 0; green, 0; blue, 0 }  ][line width=0.08]  [draw opacity=0] (8.04,-3.86) -- (0,0) -- (8.04,3.86) -- cycle    ;
%Curve Lines [id:da5917951305830498] 
\draw [color={rgb, 255:red, 0; green, 0; blue, 0 }  ,draw opacity=1 ]   (287.77,11271.01) .. controls (274.41,11247.24) and (303.18,11240.88) .. (297.11,11267.15) ;
\draw [shift={(296.44,11269.68)}, rotate = 286.35] [fill={rgb, 255:red, 0; green, 0; blue, 0 }  ,fill opacity=1 ][line width=0.08]  [draw opacity=0] (8.04,-3.86) -- (0,0) -- (8.04,3.86) -- cycle    ;
%Curve Lines [id:da7839983382988429] 
\draw [color={rgb, 255:red, 0; green, 0; blue, 0 }  ,draw opacity=1 ] [dash pattern={on 4.5pt off 4.5pt}]  (133.77,11376.01) .. controls (120.41,11352.24) and (149.18,11345.88) .. (143.11,11372.15) ;
\draw [shift={(142.44,11374.68)}, rotate = 286.35] [fill={rgb, 255:red, 0; green, 0; blue, 0 }  ,fill opacity=1 ][line width=0.08]  [draw opacity=0] (8.04,-3.86) -- (0,0) -- (8.04,3.86) -- cycle    ;

% Text Node
\draw (73.45,11276.93) node  [font=\tiny,rotate=-1.27] [align=left] {$\displaystyle i$};
% Text Node
\draw (173.45,11276.93) node  [font=\tiny,rotate=-1.27] [align=left] {$\displaystyle i$};
% Text Node
\draw (213.45,11296.93) node  [font=\tiny,rotate=-1.27] [align=left] {$\displaystyle k$};
% Text Node
\draw (213.45,11257.93) node  [font=\tiny,rotate=-1.27] [align=left] {$\displaystyle j$};
% Text Node
\draw (294.45,11276.93) node  [font=\tiny,rotate=-1.27] [align=left] {$\displaystyle i$};
% Text Node
\draw (340.45,11276.93) node  [font=\tiny,rotate=-1.27] [align=left] {$\displaystyle j$};
% Text Node
\draw (140.45,11381.93) node  [font=\tiny,rotate=-1.27] [align=left] {$\displaystyle i$};
% Text Node
\draw (232.45,11381.93) node  [font=\tiny,rotate=-1.27] [align=left] {$\displaystyle i$};
% Text Node
\draw (282.45,11381.93) node  [font=\tiny,rotate=-1.27] [align=left] {$\displaystyle j$};
% Text Node
\draw (74,11326.33) node   [align=left] {\begin{minipage}[lt]{15.64pt}\setlength\topsep{0pt}
(a)
\end{minipage}};
% Text Node
\draw (200,11327.33) node   [align=left] {\begin{minipage}[lt]{15.64pt}\setlength\topsep{0pt}
(b)
\end{minipage}};
% Text Node
\draw (317,11327.33) node   [align=left] {\begin{minipage}[lt]{15.64pt}\setlength\topsep{0pt}
(c)
\end{minipage}};
% Text Node
\draw (142,11422.33) node   [align=left] {\begin{minipage}[lt]{15.64pt}\setlength\topsep{0pt}
(d)
\end{minipage}};
% Text Node
\draw (257,11423.33) node   [align=left] {\begin{minipage}[lt]{15.64pt}\setlength\topsep{0pt}
(e)
\end{minipage}};

\end{tikzpicture}

%%%%%%%%%%5
\caption{(a)-(e) Different cases for the white out-neighbors of a white node $i$.}
\label{white}
\end{figure}
where $\mathcal{A}_{1}\in\{*,0,?\}^{r\times r}$, $\mathcal{A}_2\in \mathbb{R}^{r\times (n-r)}$, $\mathcal{A}_3\in\mathbb{R}^{(n-r)\times r}$, and     
$\mathcal{A}_{4}\in\{*,0,?\}^{(n-r)\times (n-r)}$. Let $\Scale[1]{\nu^T=\begin{bmatrix}0_{r}^T & \mathbf{1}^T_{n-r}\end{bmatrix}}$, and assume that $\mu\in\Delta\cap \mathbb{R}$. Note that no node in $V_{\mathcal{N}}$ has exactly one  strong out-neighbor in $\mathcal{W}$. 
{Otherwise,} the  coloring process can be applied, contradicting the definition of the derived set of $V_C$. 
Therefore, either every node in $\mathcal{B}$ has no white out-neighbor, or it has at least one white weak or two white strong out-neighbors. Then, in every column of $\mathcal{A}_{3}$, if all entries are not zero, then either there is at least one question mark ? (i.e., an arbitrary entry, which can be zero or nonzero) or there are at least two stars * (i.e., nonzero entries). Thus, the  nonzero and arbitrary entries of $\mathcal{A}_{3}$ can be chosen in a way that $\mathbf{1}_{n-r}^T A_{3}=0$. Now, consider a node $i$ in $\mathcal{W}$. Then, we have one of the following cases, which are shown in Figs.~\ref{white}(a)-(e), respectively: (a) $i$  has no white out-neighbor; (b) $i$ has at least two white strong out-neighbors $j$ and $k$ ($j,k\neq i$); (c) $i$ has at least two white strong out-neighbors $i$ and $j$; (d) $i$ has at least one white weak out-neighbor $i$; (e) $i$ has at least one white weak out-neighbor $j$ ($j\neq i$). Since node $i$ in $G_{\mathcal{N}}$ has no self-loop in cases (a), (b), and (e), then $f_{\Delta}^{\mathcal{L}}(i)=0$, and one can conclude that $\Delta=\{\mu\}$, and $A(i,i)=\mu$, for all $ A \in \mathcal{S}(f^{\mathcal{L}}_{\Delta}, \mathcal{A})$. Moreover, in case (b), for all $ A \in \mathcal{S}(f^{\mathcal{L}}_{\Delta}, \mathcal{A})$, we have $A(j,i)\neq 0$ and $A(k,i)\neq 0$. Also, in case (e), $A(j,i)$ can be zero or nonzero. Consider now case (c), where $i$ is strong out-neighbor of itself. Then, since $f_{\Delta}^{\mathcal{L}}(i)=*$, we have $A(i,i)\cap \Delta=\emptyset$. Then, because $\mu\in \Delta$, one can see that $A(i,i)\neq \mu$,  for all $ A \in \mathcal{S}(f^{\mathcal{L}}_{\Delta}, \mathcal{A})$. Moreover, one has $A(j,i)\neq 0$. Finally, in case (d), one can see that $f_{\Delta}^{\mathcal{L}}(i)=?$, and thus, either $A(i,i)$ can have any real value, or for the set $\Delta$ with $|\Delta|>1$, we have $A(i,i)\in \Delta$, for all $ A \in \mathcal{S}(f^{\mathcal{L}}_{\Delta}, \mathcal{A})$. Overall, by considering all of these cases, one can conclude that the arbitrary and nonzero entries of any column in $\mathcal{A}_4$ can be chosen in such a way that $\mathbf{1}_{n-r}^T A_{4}=\mu$. Consequently, we have some $A\in\mathcal{S}(f^{\mathcal{L}}_{\Delta}, \mathcal{A})$ for which $\nu^T A=\mu \nu^T$. Moreover, $\nu^T B=0$ as $V_C\subseteq \mathcal{B}$. Hence, there is some $A\in\mathcal{S}(f^{\mathcal{L}}_{\Delta}, \mathcal{A})$ with uncontrollable eigenvalue $\mu\in\Delta$, establishing a contradiction.
\carre

\begin{figure}[t]
\centering
%%%%%%%%%%%%%%%%%%%%%%%%

\tikzset{every picture/.style={line width=0.75pt}} %set default line width to 0.75pt        

\begin{tikzpicture}[x=0.75pt,y=0.75pt,yscale=-.96,xscale=.96]
%uncomment if require: \path (0,11616); %set diagram left start at 0, and has height of 11616

%Shape: Circle [id:dp3422310031412428] 
\draw   (166.44,10607.68) .. controls (171,10607.67) and (174.7,10611.36) .. (174.7,10615.91) .. controls (174.71,10620.47) and (171.02,10624.17) .. (166.47,10624.18) .. controls (161.91,10624.19) and (158.21,10620.5) .. (158.2,10615.94) .. controls (158.2,10611.39) and (161.88,10607.69) .. (166.44,10607.68) -- cycle ;
%Shape: Circle [id:dp6464158855122366] 
\draw  [fill={rgb, 255:red, 255; green, 255; blue, 255 }  ,fill opacity=1 ] (220.11,10567.68) .. controls (224.66,10567.67) and (228.36,10571.36) .. (228.37,10575.91) .. controls (228.38,10580.47) and (224.69,10584.17) .. (220.13,10584.18) .. controls (215.58,10584.19) and (211.88,10580.5) .. (211.87,10575.94) .. controls (211.86,10571.39) and (215.55,10567.69) .. (220.11,10567.68) -- cycle ;
%Shape: Circle [id:dp11853407917026781] 
\draw   (185.44,10672.68) .. controls (190,10672.67) and (193.7,10676.36) .. (193.7,10680.91) .. controls (193.71,10685.47) and (190.02,10689.17) .. (185.47,10689.18) .. controls (180.91,10689.19) and (177.21,10685.5) .. (177.2,10680.94) .. controls (177.2,10676.39) and (180.88,10672.69) .. (185.44,10672.68) -- cycle ;
%Shape: Circle [id:dp23901468798420678] 
\draw   (254.77,10673.01) .. controls (259.33,10673) and (263.03,10676.69) .. (263.04,10681.25) .. controls (263.04,10685.8) and (259.36,10689.5) .. (254.8,10689.51) .. controls (250.24,10689.52) and (246.54,10685.83) .. (246.54,10681.28) .. controls (246.53,10676.72) and (250.22,10673.02) .. (254.77,10673.01) -- cycle ;
%Shape: Circle [id:dp11987498302209376] 
\draw   (274.77,10607.68) .. controls (279.33,10607.67) and (283.03,10611.36) .. (283.04,10615.91) .. controls (283.04,10620.47) and (279.36,10624.17) .. (274.8,10624.18) .. controls (270.24,10624.19) and (266.54,10620.5) .. (266.54,10615.94) .. controls (266.53,10611.39) and (270.22,10607.69) .. (274.77,10607.68) -- cycle ;
%Straight Lines [id:da12227164074625074] 
\draw  [dash pattern={on 4.5pt off 4.5pt}]  (180.68,10670.78) -- (168.47,10624.18) ;
\draw [shift={(181.44,10673.68)}, rotate = 255.32] [fill={rgb, 255:red, 0; green, 0; blue, 0 }  ][line width=0.08]  [draw opacity=0] (8.04,-3.86) -- (0,0) -- (8.04,3.86) -- cycle    ;
%Straight Lines [id:da4261470127459488] 
\draw    (272.8,10624.18) -- (258.71,10670.3) ;
\draw [shift={(257.83,10673.17)}, rotate = 286.99] [fill={rgb, 255:red, 0; green, 0; blue, 0 }  ][line width=0.08]  [draw opacity=0] (8.04,-3.86) -- (0,0) -- (8.04,3.86) -- cycle    ;
%Straight Lines [id:da1268504048415786] 
\draw    (230.75,10576.74) -- (271.33,10608) ;
\draw [shift={(228.37,10574.91)}, rotate = 37.6] [fill={rgb, 255:red, 0; green, 0; blue, 0 }  ][line width=0.08]  [draw opacity=0] (8.04,-3.86) -- (0,0) -- (8.04,3.86) -- cycle    ;
%Curve Lines [id:da46968510249975304] 
\draw [color={rgb, 255:red, 0; green, 0; blue, 0 }  ,draw opacity=1 ]   (215.11,10568.68) .. controls (210.92,10540.05) and (227.54,10540.91) .. (224.21,10565.57) ;
\draw [shift={(223.77,10568.35)}, rotate = 280.37] [fill={rgb, 255:red, 0; green, 0; blue, 0 }  ,fill opacity=1 ][line width=0.08]  [draw opacity=0] (8.04,-3.86) -- (0,0) -- (8.04,3.86) -- cycle    ;
%Curve Lines [id:da68791068936826] 
\draw [color={rgb, 255:red, 0; green, 0; blue, 0 }  ,draw opacity=1 ]   (277.77,10607.68) .. controls (281.55,10585.03) and (305.81,10594.39) .. (285.31,10612.02) ;
\draw [shift={(283.27,10613.68)}, rotate = 322.18] [fill={rgb, 255:red, 0; green, 0; blue, 0 }  ,fill opacity=1 ][line width=0.08]  [draw opacity=0] (8.04,-3.86) -- (0,0) -- (8.04,3.86) -- cycle    ;
%Curve Lines [id:da4612131205446659] 
\draw [color={rgb, 255:red, 0; green, 0; blue, 0 }  ,draw opacity=1 ]   (184.47,10689.18) .. controls (175.38,10716.21) and (156.89,10697.32) .. (175.04,10686.18) ;
\draw [shift={(177.5,10684.83)}, rotate = 513.95] [fill={rgb, 255:red, 0; green, 0; blue, 0 }  ,fill opacity=1 ][line width=0.08]  [draw opacity=0] (8.04,-3.86) -- (0,0) -- (8.04,3.86) -- cycle    ;
%Curve Lines [id:da8396095116481712] 
\draw [color={rgb, 255:red, 0; green, 0; blue, 0 }  ,draw opacity=1 ][line width=0.75]  [dash pattern={on 4.5pt off 4.5pt}]  (260.8,10686.51) .. controls (282.22,10698.24) and (259.69,10716.35) .. (255.21,10692.3) ;
\draw [shift={(254.8,10689.51)}, rotate = 443.52] [fill={rgb, 255:red, 0; green, 0; blue, 0 }  ,fill opacity=1 ][line width=0.08]  [draw opacity=0] (8.04,-3.86) -- (0,0) -- (8.04,3.86) -- cycle    ;
%Shape: Circle [id:dp7787327057337772] 
\draw  [fill={rgb, 255:red, 0; green, 0; blue, 0 }  ,fill opacity=1 ] (357.44,10607.68) .. controls (362,10607.67) and (365.7,10611.36) .. (365.7,10615.91) .. controls (365.71,10620.47) and (362.02,10624.17) .. (357.47,10624.18) .. controls (352.91,10624.19) and (349.21,10620.5) .. (349.2,10615.94) .. controls (349.2,10611.39) and (352.88,10607.69) .. (357.44,10607.68) -- cycle ;
%Shape: Circle [id:dp991132130878198] 
\draw   (411.11,10567.68) .. controls (415.66,10567.67) and (419.36,10571.36) .. (419.37,10575.91) .. controls (419.38,10580.47) and (415.69,10584.17) .. (411.13,10584.18) .. controls (406.58,10584.19) and (402.88,10580.5) .. (402.87,10575.94) .. controls (402.86,10571.39) and (406.55,10567.69) .. (411.11,10567.68) -- cycle ;
%Shape: Circle [id:dp995859118209991] 
\draw   (376.44,10672.68) .. controls (381,10672.67) and (384.7,10676.36) .. (384.7,10680.91) .. controls (384.71,10685.47) and (381.02,10689.17) .. (376.47,10689.18) .. controls (371.91,10689.19) and (368.21,10685.5) .. (368.2,10680.94) .. controls (368.2,10676.39) and (371.88,10672.69) .. (376.44,10672.68) -- cycle ;
%Shape: Circle [id:dp299240507100516] 
\draw   (445.77,10673.01) .. controls (450.33,10673) and (454.03,10676.69) .. (454.04,10681.25) .. controls (454.04,10685.8) and (450.36,10689.5) .. (445.8,10689.51) .. controls (441.24,10689.52) and (437.54,10685.83) .. (437.54,10681.28) .. controls (437.53,10676.72) and (441.22,10673.02) .. (445.77,10673.01) -- cycle ;
%Shape: Circle [id:dp6603727468281022] 
\draw   (465.77,10607.68) .. controls (470.33,10607.67) and (474.03,10611.36) .. (474.04,10615.91) .. controls (474.04,10620.47) and (470.36,10624.17) .. (465.8,10624.18) .. controls (461.24,10624.19) and (457.54,10620.5) .. (457.54,10615.94) .. controls (457.53,10611.39) and (461.22,10607.69) .. (465.77,10607.68) -- cycle ;
%Straight Lines [id:da4927992135157788] 
\draw  [dash pattern={on 4.5pt off 4.5pt}]  (371.68,10670.78) -- (359.47,10624.18) ;
\draw [shift={(372.44,10673.68)}, rotate = 255.32] [fill={rgb, 255:red, 0; green, 0; blue, 0 }  ][line width=0.08]  [draw opacity=0] (8.04,-3.86) -- (0,0) -- (8.04,3.86) -- cycle    ;
%Straight Lines [id:da22969342642302015] 
\draw    (463.8,10624.18) -- (449.71,10670.3) ;
\draw [shift={(448.83,10673.17)}, rotate = 286.99] [fill={rgb, 255:red, 0; green, 0; blue, 0 }  ][line width=0.08]  [draw opacity=0] (8.04,-3.86) -- (0,0) -- (8.04,3.86) -- cycle    ;
%Straight Lines [id:da6661162294396039] 
\draw    (421.75,10576.74) -- (462.33,10608) ;
\draw [shift={(419.37,10574.91)}, rotate = 37.6] [fill={rgb, 255:red, 0; green, 0; blue, 0 }  ][line width=0.08]  [draw opacity=0] (8.04,-3.86) -- (0,0) -- (8.04,3.86) -- cycle    ;
%Curve Lines [id:da002029247958623026] 
\draw [color={rgb, 255:red, 0; green, 0; blue, 0 }  ,draw opacity=1 ]   (406.11,10568.68) .. controls (401.92,10540.05) and (418.54,10540.91) .. (415.21,10565.57) ;
\draw [shift={(414.77,10568.35)}, rotate = 280.37] [fill={rgb, 255:red, 0; green, 0; blue, 0 }  ,fill opacity=1 ][line width=0.08]  [draw opacity=0] (8.04,-3.86) -- (0,0) -- (8.04,3.86) -- cycle    ;
%Curve Lines [id:da8102015729152838] 
\draw [color={rgb, 255:red, 0; green, 0; blue, 0 }  ,draw opacity=1 ] [dash pattern={on 4.5pt off 4.5pt}]  (349.31,10615.12) .. controls (323.3,10602.88) and (341.34,10581.81) .. (351.52,10606.88) ;
\draw [shift={(352.44,10609.35)}, rotate = 250.99] [fill={rgb, 255:red, 0; green, 0; blue, 0 }  ,fill opacity=1 ][line width=0.08]  [draw opacity=0] (8.04,-3.86) -- (0,0) -- (8.04,3.86) -- cycle    ;
%Curve Lines [id:da0213751544162315] 
\draw [color={rgb, 255:red, 0; green, 0; blue, 0 }  ,draw opacity=1 ] [dash pattern={on 4.5pt off 4.5pt}]  (468.77,10607.68) .. controls (472.55,10585.03) and (496.81,10594.39) .. (476.31,10612.02) ;
\draw [shift={(474.27,10613.68)}, rotate = 322.18] [fill={rgb, 255:red, 0; green, 0; blue, 0 }  ,fill opacity=1 ][line width=0.08]  [draw opacity=0] (8.04,-3.86) -- (0,0) -- (8.04,3.86) -- cycle    ;
%Curve Lines [id:da18216832871430944] 
\draw [color={rgb, 255:red, 0; green, 0; blue, 0 }  ,draw opacity=1 ]   (375.47,10689.18) .. controls (366.38,10716.21) and (347.89,10697.32) .. (366.04,10686.18) ;
\draw [shift={(368.5,10684.83)}, rotate = 513.95] [fill={rgb, 255:red, 0; green, 0; blue, 0 }  ,fill opacity=1 ][line width=0.08]  [draw opacity=0] (8.04,-3.86) -- (0,0) -- (8.04,3.86) -- cycle    ;
%Curve Lines [id:da823259751529257] 
\draw [color={rgb, 255:red, 0; green, 0; blue, 0 }  ,draw opacity=1 ][line width=0.75]  [dash pattern={on 4.5pt off 4.5pt}]  (452.8,10685.51) .. controls (474.22,10697.24) and (452.62,10716.28) .. (448.21,10692.29) ;
\draw [shift={(447.8,10689.51)}, rotate = 443.52] [fill={rgb, 255:red, 0; green, 0; blue, 0 }  ,fill opacity=1 ][line width=0.08]  [draw opacity=0] (8.04,-3.86) -- (0,0) -- (8.04,3.86) -- cycle    ;
%Straight Lines [id:da4412839555872452] 
\draw    (209.5,10577.78) -- (170,10608.5) ;
\draw [shift={(211.87,10575.94)}, rotate = 142.13] [fill={rgb, 255:red, 0; green, 0; blue, 0 }  ][line width=0.08]  [draw opacity=0] (8.04,-3.86) -- (0,0) -- (8.04,3.86) -- cycle    ;
%Straight Lines [id:da36450420737440314] 
\draw    (400.5,10577.78) -- (361,10608.5) ;
\draw [shift={(402.87,10575.94)}, rotate = 142.13] [fill={rgb, 255:red, 0; green, 0; blue, 0 }  ][line width=0.08]  [draw opacity=0] (8.04,-3.86) -- (0,0) -- (8.04,3.86) -- cycle    ;
%Straight Lines [id:da5037526544438888] 
\draw    (454.54,10615.94) -- (365.7,10615.91) ;
\draw [shift={(457.54,10615.94)}, rotate = 180.02] [fill={rgb, 255:red, 0; green, 0; blue, 0 }  ][line width=0.08]  [draw opacity=0] (8.04,-3.86) -- (0,0) -- (8.04,3.86) -- cycle    ;
%Straight Lines [id:da953071480582091] 
\draw    (263.54,10615.94) -- (174.7,10615.91) ;
\draw [shift={(266.54,10615.94)}, rotate = 180.02] [fill={rgb, 255:red, 0; green, 0; blue, 0 }  ][line width=0.08]  [draw opacity=0] (8.04,-3.86) -- (0,0) -- (8.04,3.86) -- cycle    ;
%Straight Lines [id:da009667005739230738] 
\draw    (196.7,10680.93) -- (246.54,10681.28) ;
\draw [shift={(193.7,10680.91)}, rotate = 0.39] [fill={rgb, 255:red, 0; green, 0; blue, 0 }  ][line width=0.08]  [draw opacity=0] (8.04,-3.86) -- (0,0) -- (8.04,3.86) -- cycle    ;
%Straight Lines [id:da887039505252087] 
\draw    (387.7,10680.93) -- (437.54,10681.28) ;
\draw [shift={(384.7,10680.91)}, rotate = 0.39] [fill={rgb, 255:red, 0; green, 0; blue, 0 }  ][line width=0.08]  [draw opacity=0] (8.04,-3.86) -- (0,0) -- (8.04,3.86) -- cycle    ;
%Straight Lines [id:da5338937115372155] 
\draw    (405.19,10585.39) -- (380.44,10673.68) ;
\draw [shift={(406,10582.5)}, rotate = 105.66] [fill={rgb, 255:red, 0; green, 0; blue, 0 }  ][line width=0.08]  [draw opacity=0] (8.04,-3.86) -- (0,0) -- (8.04,3.86) -- cycle    ;
%Straight Lines [id:da041951628216437165] 
\draw    (214.19,10585.39) -- (189.44,10673.68) ;
\draw [shift={(215,10582.5)}, rotate = 105.66] [fill={rgb, 255:red, 0; green, 0; blue, 0 }  ][line width=0.08]  [draw opacity=0] (8.04,-3.86) -- (0,0) -- (8.04,3.86) -- cycle    ;
%Straight Lines [id:da6683753368645651] 
\draw    (225.81,10585.39) -- (250.77,10674.01) ;
\draw [shift={(225,10582.5)}, rotate = 74.27] [fill={rgb, 255:red, 0; green, 0; blue, 0 }  ][line width=0.08]  [draw opacity=0] (8.04,-3.86) -- (0,0) -- (8.04,3.86) -- cycle    ;
%Straight Lines [id:da849878307527882] 
\draw    (416.81,10585.39) -- (441.77,10674.01) ;
\draw [shift={(416,10582.5)}, rotate = 74.27] [fill={rgb, 255:red, 0; green, 0; blue, 0 }  ][line width=0.08]  [draw opacity=0] (8.04,-3.86) -- (0,0) -- (8.04,3.86) -- cycle    ;

% Text Node
\draw (167.45,10615.93) node  [font=\tiny,rotate=-1.27] [align=left] {1};
% Text Node
\draw (186.45,10680.93) node  [font=\tiny,rotate=-1.27] [align=left] {2};
% Text Node
\draw (255.79,10681.26) node  [font=\tiny,rotate=-1.27] [align=left] {3};
% Text Node
\draw (275.79,10615.93) node  [font=\tiny,rotate=-1.27] [align=left] {4};
% Text Node
\draw (358.45,10615.93) node  [font=\tiny,rotate=-1.27] [align=left] {\textcolor[rgb]{1,1,1}{1}};
% Text Node
\draw (412.12,10575.93) node  [font=\tiny,rotate=-1.27] [align=left] {5};
% Text Node
\draw (377.45,10680.93) node  [font=\tiny,rotate=-1.27] [align=left] {2};
% Text Node
\draw (446.79,10681.26) node  [font=\tiny,rotate=-1.27] [align=left] {3};
% Text Node
\draw (466.79,10615.93) node  [font=\tiny,rotate=-1.27] [align=left] {4};
% Text Node
\draw (221,10734) node   [align=left] {\begin{minipage}[lt]{15.64pt}\setlength\topsep{0pt}
(a)
\end{minipage}};
% Text Node
\draw (412,10734) node   [align=left] {\begin{minipage}[lt]{15.64pt}\setlength\topsep{0pt}
(b)
\end{minipage}};
% Text Node
\draw (297.25,10588.83) node   [align=left] {\begin{minipage}[lt]{9.180000000000001pt}\setlength\topsep{0pt}
\textbf{{\Large +}}
\end{minipage}};
% Text Node
\draw (171.25,10709.83) node   [align=left] {\begin{minipage}[lt]{9.180000000000001pt}\setlength\topsep{0pt}
\textbf{{\Large -}}
\end{minipage}};
% Text Node
\draw (223.25,10539.83) node   [align=left] {\begin{minipage}[lt]{9.180000000000001pt}\setlength\topsep{0pt}
\textbf{{\Large -}}
\end{minipage}};
% Text Node
\draw (221.12,10575.93) node  [font=\tiny,rotate=-1.27] [align=left] {5};

\end{tikzpicture}

%%%%%%%%%%5
\caption{(a) Global graph $\mathcal{G}$; (b) corresponding $\Delta$-network graph $G_{\mathcal{N}}^{\Delta}$.}
\label{stabexample}
\end{figure}
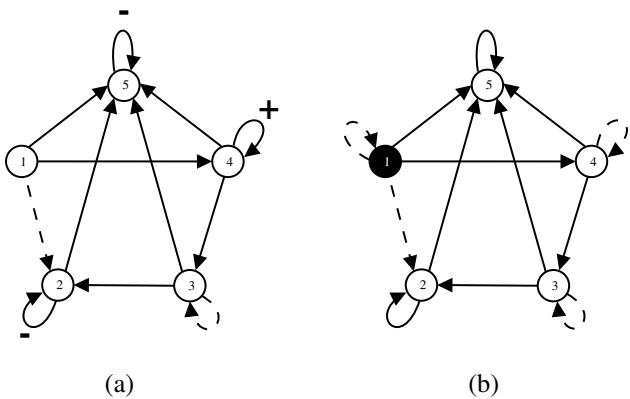

\begin{ex}
Consider an LTI N1DS with the global graph $\mathcal{G}$ in Fig.~\ref{stabexample}(a). Let $\mathcal{A}$ be the associated pattern matrix, and let $\Delta=\{y\in \mathbb{C}\mid\Re(y)\geq 0\}$. Furthermore, $\mathcal{S}(f_{\Delta}^\mathcal{L}, \mathcal{A})$ is the set of all $A\in\mathcal{P}(\mathcal{A})$, where $A(1,1)=0$, $A(2,2)< 0$, $A(4,4)>0$, and $A(5,5)<0$.  As discussed in Example \ref{stab}, if $i=2$ or $i=5$, then $f_{\Delta}^{\mathcal{L}}(i)=*$, and otherwise, $f_{\Delta}^{\mathcal{L}}(i)=?$. Thus, the corresponding $\Delta$-network graph $G_{\mathcal{N}}^{\Delta}$, depicted in Fig.~\ref{stabexample}(b), is obtained. From Theorem \ref{th5}, the network is strongly structurally stabilizable if and only if $V_C$ is a ZFS of $G_{\mathcal{N}}^{\Delta}$. For example,    one can see that $V_C=\{1\}$ is a ZFS of $G_{\mathcal{N}}^{\Delta}$ and renders the network strongly structurally stabilizable. 
\end{ex}
%%%%%%%%%%%%%%%%%%%%%%%%%%%%%%%%%%%%%%%%%%%%%%%%%%%%%%%%%%%%%%%%%%%%

%%%%%%%%%%%%%%%%%%%%%%%%%%%%%%%%%%%%%%%%%%%%%%%%%%%%%%%%%%%%%%%%%%%%5

\begin{ex}
Consider the LTI network of Example \ref{ex10} (Fig. \ref{stable}). It was shown in Example \ref{ex11} that $V_C=\{1\}$ is a ZFS of the $\Delta$-network graph for both $\Delta=\{0\}$ and $\Delta=\mathbb{C}\setminus \{0\}$. Thus,  the network is strongly structurally controllable with respect to both $\Delta=\{0\}$ and $\Delta=\mathbb{C}\setminus \{0\}$.  Then, the network is strongly structurally controllable.  
\end{ex}

\subsection{Comparison to Existing Results}

In this part, we compare Theorem \ref{th5} of this paper with some of the existing results on strong structural controllability of LTI networks in the literature. In fact, the modal strong structural controllability results presented in the current work are based on forming the $\Delta$-network graph $G_{\mathcal{N}}^{\Delta}$ associated with a set of eigenvalues $\Delta$, and then finding a ZFS in $G_{\mathcal{N}}^{\Delta}$. In \cite{jia2019strong}, it has been shown that the algebraic and graph-theoretic conditions provided in the works \cite{mayeda1979strong, Olesky1993, Reinschke1992, jarczyk2011strong, chapman2013strong,  monshizadeh2014zero, trefois2015zero } are special cases of the main result of \cite{jia2019strong}. Thus, it suffices to make a comparison between the results of this paper and the work in \cite{jia2019strong}, whose main result is presented in the following.

\begin{thm}[see Theorem 11 of \cite{jia2019strong}]
Consider an N1DS with dynamics (\ref{e1}) and the pattern matrix $\mathcal{A}\in \{*,0,?\}^{n\times n}$. Let $\bar{\mathcal{A}}\in \{*,0,?\}^{n\times n}$ be another pattern matrix, where for $i\neq j$, $\mathcal{A}(i,j)=\bar{\mathcal{A}}(i,j)$, and for $1\leq i\leq n$,
\begin{equation}
   \bar{\mathcal{A}}(i,i)=
    \begin{cases}
      *, & \text{if}\  \mathcal{A}(i,i)=0,  
      \\
      ?, & \text{otherwise}.
    \end{cases}
    \label{Abar}
\end{equation}
Then, the network is strongly structurally controllable if and only if $V_C$ is a ZFS of both the graphs $G_1=G(\mathcal{A})$ and $G_2=G(\bar{\mathcal{A}})$.
\label{th_11}
\end{thm}

In order to compare Theorems \ref{th5} and \ref{th_11}, we note that a very special example of our result is an N1DS with  $\Delta=\{0\}$ and  $\Delta'=\mathbb{C}\setminus \{0\}$. From Example \ref{ex5}, one can conclude that the $\Delta$-network graph associated to an N1DS with the pattern matrix $\mathcal{A}$ is the same as the graph $G_1=G(\mathcal{A})$. Now, consider the pattern matrix $\bar{\mathcal{A}}$ in (\ref{Abar}). As inferred from Example \ref{ex6}, the $\Delta'$-network graph associated with an N1DS with the pattern matrix $\mathcal{A}$ can be defined as $G_2=G(\bar{\mathcal{A}})$.  Based on the definition, the network is strongly structurally controllable if it is both $\Delta$-SSC and $\Delta'$-SSC. Moreover, in view of  Theorem \ref{th5}, an N1DS is $\Delta$-SSC (respectively, $\Delta'$-SSC) if and only if $V_C$ is a ZFS of the $\Delta$-network graph  $G_1$ (respectively, $\Delta'$-network graph $G_2$),  verifying Theorem \ref{th_11}.

Now, let us discuss the result of the work \cite{monshizadeh2014zero}
as well. In this work, an N1DS with an arbitrary-diagonal pattern matrix $\mathcal{A}\in \{*,0,?\}^{n\times n}$ has been investigated, where we have $\mathcal{A}(i,i)=?$, for all $1\leq i\leq n$. Then, a loop-free graph\footnote{A loop-free graph does not contain any self-loop, while a loop graph can have self-loops on some of the nodes.} $G_{\mathrm{Lf}}$ associated with $\mathcal{A}$ has been  defined as $G_{\mathrm{Lf}}=G_{\text{loop-free}}(\mathcal{A})=(V,E)$, where for $i\neq j$, $(v_j,v_i)\notin E$ if and only if $\mathcal{A}(i,j)=0$. However, based on our definition, one can see that $G_{\mathrm{L}}=G(\mathcal{A})$ is a loop graph whose every node is a weak out-neighbor of itself. 

In \cite{monshizadeh2014zero}, an ``ordinary coloring rule" is introduced, which states that a white node $u$ can become black if it is the only white out-neighbor of a ``black node'' $v$. Now, let $Z$ be an initial set of black nodes in $G_{\mathrm{Lf}}$. Set $Z$ is called an ``ordinary zero forcing set'' if by the repeated application of the coloring rule in $G_{\mathrm{Lf}}$, all nodes become black.  The next result from \cite{monshizadeh2014zero} provides a graph-theoretic condition for strong structural controllability. 

\begin{thm}
An N1DS with dynamics (\ref{e1}) and an arbitrary-diagonal pattern matrix $\mathcal{A}$ is strongly structurally controllable if and only if $V_C$ is an ordinary ZFS of $G_{\mathrm{Lf}}$.
\end{thm}

Considering Example \ref{ex8}, we note that in our framework, the $\Delta$-network graph for any nonempty $\Delta\subseteq \mathbb{C}$ is the same as $G_{\mathrm{L}}=G(A)$, where there is a self-loop on every node of the graph, and all self-loops are represented  by dotted arrows. Therefore, based on Theorem \ref{th5}, given a $\Delta\subseteq \mathbb{C}$, where $\Delta\cap\mathbb{R}\neq \emptyset$, this network is $\Delta$-SSC if and only if $V_C$ is a ZFS of $G_{\mathrm{L}}$. Moreover, in \cite{jia2019strong, trefois2015zero}, it has been demonstrated that $V_C$ is an ordinary ZFS of $G_{\mathrm{Lf}}$ if and only if it is a ZFS of $G_{\mathrm{L}}$. Thus, one can see that our result is in line with the result of \cite{monshizadeh2014zero}.

\begin{corr}
Consider two arbitrary nonempty sets of numbers $\Delta_1$ and $\Delta_2$, where $\Delta_i\subseteq \mathbb{C}$ and $\Delta_i\cap\mathbb{R}=\emptyset$, $i=1,2$. An N1DS with dynamics (\ref{e1}) and an arbitrary-diagonal pattern matrix is $\Delta_1$-SSC if and only if it is $\Delta_2$-SSC.
\end{corr}

%\begin{rem}
 %   The modal strong structural controllability results in this work   are based on forming the $\Delta$-network graph $G_{\mathcal{N}}^{\Delta}$ associated with a set of eigenvalues $\Delta$ and then finding a ZFS in $G_{\mathcal{N}}^{\Delta}$. A very special example that we can mention is an N1DS with  $\Delta=\{0\}$ or $\Delta=\mathbb{R}\setminus \{0\}$, for which our results verifies the results of \cite{trefois2015zero, jia2020unifying}. Moreover,  we note that for all cases discussed in Examples \ref{ex7} and \ref{ex8}, every node of the $\Delta$-network graph is  a weak out-neighbor of itself. Thus, in the coloring process applied in these networks,  a white node $u$ can be black if it is the only white strong out-neighbor of a \emph{black} node, which is in line with the  results of \cite{monshizadeh2014zero}. Considering the network in Example \ref{ex8}, one can find  that if the network is $\Delta$-SSC for $\Delta=\{0\}$, it is also $\Delta'$-SSC, for every nonempty $\Delta'\subseteq \mathbb{R}$.   
%\end{rem}

\subsection{Maximum Geometric Multiplicity and Zero Forcing Sets}

Consider a block matrix $A$ in (\ref{sm}) with $A_{ii}\in\mathbb{R}^{l_i\times l_i}$. For a pattern matrix $\mathcal{A}$, a characteristic vector $f^{\mathcal{L}}_{\Delta}$,  and a nonempty set $\Delta\subseteq\mathbb{C}$, let $A\in\mathcal{S}(f^{\mathcal{L}}_{\Delta}, \mathcal{A})$. Now, consider the corresponding $\Delta$-network graph $G_{\mathcal{N}}^{\Delta}$. A minimal ZFS, denoted by $\mathrm{ZFS}_{\mathrm{min}}$, is a ZFS of $G_{\mathcal{N}}^{\Delta}$, where $\mbox{Ver}(\mathrm{ZFS}_{\mathrm{min}})$ has the minimum cardinality. The cardinality of the set of vertices of a minimal ZFS (i.e.,  $|\mbox{Ver}(\mathrm{ZFS}_{\mathrm{min}})|$)  is called the  zero forcing number and denoted by $\mathcal{Z}(G_{\mathcal{N}}^{\Delta})$.

{There is a relation between $\mathcal{Z}(G_{\mathcal{N}}^{\Delta})$ and the maximum nullity of $\lambda I-A$, for all $A\in\mathcal{S}(f^{\mathcal{L}}_{\Delta}, \mathcal{A})$ and $\lambda\in\Delta$.}  By this relation, we can provide an upper bound on the maximum geometric multiplicity of  eigenvalues of all $A\in\mathcal{S}(f^{\mathcal{L}}_{\Delta}, \mathcal{A})$ that belong to the set $\Delta$.

\begin{pro}
Given the $\Delta$-network graph $G_{\mathcal{N}}^{\Delta}=(V_{\mathcal{N}}, E_{\mathcal{N}})$ with the  zero forcing number $\mathcal{Z}(G_{\mathcal{N}}^{\Delta})$, we have $\Psi_{{\Delta}}(A)\leq~\mathcal{Z}(G_{\mathcal{N}})$,
for all $A\in\mathcal{S}(f^{\mathcal{L}}_{\Delta}, \mathcal{A})$.
\label{p8}
\end{pro}

\textbf{Proof.} 
Suppose that for some $\lambda\in\Delta$, $\psi_A(\lambda)=k$. Then, for every $X\subset V_{\mathcal{N}}$ with $|\mbox{Ver}(X)|=k-1$, there is a nonzero $\nu\in \mathcal{V}_A(\lambda)$ such that $\nu_i=0$, for all $i\in \mbox{Ver}(X)$ (see proof of Proposition 2.2 in \cite{work2008zero}). Now, let $\beta\in \Delta$ be an eigenvalue of some $A\in\mathcal{S}(f^{\mathcal{L}}_{\Delta}, \mathcal{A})$, with $\psi_A(\beta)> \mathcal{Z}(G_{\mathcal{N}}^{\Delta})$. Then, for a ZFS denoted by $Z$ with $|\mbox{Ver}(Z)|=\mathcal{Z}   (G_{\mathcal{N}}^{\Delta})$, there is a nonzero $\nu\in \mathcal{V}_A(\beta)$ such that $\nu_i=0$, for all $ i \in \mbox{Ver}(Z)$. Moreover, from Lemma \ref{l2}, since $\mathcal{D}(Z)=V_{\mathcal{N}}$, we have $\nu=0$, contradicting the assumption.  \carre

For an N1DS, where $l_1=\ldots=l_n=1$, by using the notion of  ZFS in the corresponding $\Delta$-network graph, one can further examine 
whether there is any $A\in\mathcal{S}(f^{\mathcal{L}}_{\Delta}, \mathcal{A})$ %for a given $G$
that has any eigenvalue in $\Delta$ or not.
\begin{thm}
Consider $\Delta\subseteq \mathbb{C}$, where $\Delta\cap\mathbb{R}\neq\emptyset$, and let   $G_{\mathcal{N}}^{\Delta}=(V_{\mathcal{N}}, E_{\mathcal{N}})$ be a $\Delta$-network graph associated with an N1DS. {No matrix $A\in\mathcal{S}(f^{\mathcal{L}}_{\Delta}, \mathcal{A})$ has any eigenvalue in $\Delta$} {if} and only if $\mathcal{Z}(G_{\mathcal{N}}^{\Delta})=0$.
\label{th6}
\end{thm}
  
\textbf{Proof.}  
The sufficiency part of the Theorem follows immediately from Proposition \ref{p8}.
To prove the necessity,  assume that no $A\in\mathcal{S}(f^{\mathcal{L}}_{\Delta}, \mathcal{A})$ has any eigenvalue in $\Delta$, but $\mathcal{Z}(G_{\mathcal{N}}^{\Delta})>0$. Then, for every subset $\mathcal{X}\subset V_{\mathcal{N}}$ such that $|\mathcal{X}|=\mathcal{Z}(G_{\mathcal{N}}^{\Delta})-1$, one should have $\mathcal{D}(\mathcal{X})\neq V_{\mathcal{N}}$. Since we have $\Lambda(A)\cap\Delta=\emptyset$, for all $ A \in \mathcal{S}(f^{\mathcal{L}}_{\Delta}, \mathcal{A})$, then based on Definition \ref{ss},  an N1DS with dynamics (\ref{e1}) and  any choice of the set $V_C $ is $\Delta$-SSC. Then, from Theorem \ref{th5}, $V_C=\mathcal{X}$ should be a ZFS,  contradicting the assumption. \carre

\subsection{Full Rank Condition for Pattern Matrices}
 
Consider an LTI network $\mathcal{N}$ with the system matrix $A$ described in (\ref{sm}), where $A\in\mathcal{S}(f^{\mathcal{L}}_{\Delta}, \mathcal{A})$. Let $\mathcal{G}=(V_{\mathcal{G}}, E_{\mathcal{G}})$ and $G_{\mathcal{N}}^{\Delta}=(V_{\mathcal{N}}, E_{\mathcal{N}})$ be the global graph and the corresponding $\Delta$-network graph, respectively. As discussed in Section \ref{section}, for $j\neq i$, in order to find whether  $(j,i)\in E_{\mathcal{N}}^*$, one should check if the pattern matrix $\mathcal{A}_{ij}\in \{0,*,?\}^{l_i\times l_j}$ has full row rank or not, which is discussed in the following.  
  
Consider a pattern matrix $\mathcal{A}'\in \{0,*,?\}^{q\times p}$, where $l=\max (q,p)$. Then, recall that one can associate to $\mathcal{A}'$ a graph $G(\mathcal{A}')=(V,E)$, where $V=\{v_1,\ldots,v_l\}$, and $(v_j, v_i) \in E$ if and only if $\mathcal{A}'(i,j)\neq 0$. In addition, corresponding to   $\mathcal{A}'$, one can define  a  bipartite graph  as $G_b(\mathcal{A}')=(V_r,V_c,E_b)$ with $V_r=\{v^r_1,\ldots,v^r_q\}$ and $V_c=\{v^c_1,\ldots,v^c_p\}$, and $(v^c_j, v^r_i) \in E_b$ if and only if $\mathcal{A}'(i,j)\neq 0$. Thus, $G_b$ can be considered as a bipartite graph equivalent to  the graph $G(\mathcal{A}')$.

Now, let $G_i=(V_i,E_i)$ be the node graph associated with subsystem $i$, $1\leq i\leq n$, where $V_i=\{i^1,\ldots, i^{l_i}\}$.  For $i\neq j$, consider the pattern matrix $\mathcal{A}_{ij}\in \{0,*,?\}^{l_i\times l_j}$, that describes the edges connecting the vertices of $V_j$ to the vertices of $V_i$. Let us define the bipartite graph $G_{j\rightarrow i}=G_b(\mathcal{A}_{ij})=(V_i, V_j, E_{j\rightarrow i})$. Similarly, let $G_{i\rightarrow j}=G_b(\mathcal{A}_{ji})=(V_j, V_i, E_{i\rightarrow j})$. Now, as a subgraph of $\mathcal{G}$, we define the graph $G_{\{i,j\}}=(V_i\cup V_j, E_{i\rightarrow j}\cup E_{j\rightarrow i})$, whose edge set is the set of all edges in $E_{\mathcal{G}}$ with one end node in $V_i$ and the other end node in $V_j$.

Now, assume that all nodes of $G_{\{i,j\}}$ are initially white, and apply the coloring rule in $G_{\{i,j\}}$ as many times as possible. In the following, we show that pattern matrix $\mathcal{A}_{ij}$ has full row rank if and only if all vertices in $V_i$ become black.

\begin{figure}[b]
\centering

%%%%%%%%%%%%%%%%%%%%%%%%%%%%%

\tikzset{every picture/.style={line width=0.75pt}} %set default line width to 0.75pt        

\begin{tikzpicture}[x=0.75pt,y=0.75pt,yscale=-.95,xscale=.95]
%uncomment if require: \path (0,10563); %set diagram left start at 0, and has height of 10563

%Shape: Circle [id:dp7318972946450595] 
\draw  [fill={rgb, 255:red, 255; green, 255; blue, 255 }  ,fill opacity=1 ] (120.11,10423.35) .. controls (124.66,10423.34) and (128.36,10427.02) .. (128.37,10431.58) .. controls (128.38,10436.14) and (124.69,10439.84) .. (120.13,10439.85) .. controls (115.58,10439.85) and (111.88,10436.17) .. (111.87,10431.61) .. controls (111.86,10427.05) and (115.55,10423.35) .. (120.11,10423.35) -- cycle ;
%Shape: Circle [id:dp5200637243696518] 
\draw  [fill={rgb, 255:red, 255; green, 255; blue, 255 }  ,fill opacity=1 ] (120.11,10363.01) .. controls (124.66,10363) and (128.36,10366.69) .. (128.37,10371.25) .. controls (128.38,10375.8) and (124.69,10379.5) .. (120.13,10379.51) .. controls (115.58,10379.52) and (111.88,10375.83) .. (111.87,10371.28) .. controls (111.86,10366.72) and (115.55,10363.02) .. (120.11,10363.01) -- cycle ;
%Shape: Circle [id:dp6876421638120607] 
\draw  [fill={rgb, 255:red, 255; green, 255; blue, 255 }  ,fill opacity=1 ] (247.77,10363.01) .. controls (252.33,10363) and (256.03,10366.69) .. (256.04,10371.25) .. controls (256.04,10375.8) and (252.36,10379.5) .. (247.8,10379.51) .. controls (243.24,10379.52) and (239.54,10375.83) .. (239.54,10371.28) .. controls (239.53,10366.72) and (243.22,10363.02) .. (247.77,10363.01) -- cycle ;
%Shape: Circle [id:dp8779785437084049] 
\draw  [fill={rgb, 255:red, 255; green, 255; blue, 255 }  ,fill opacity=1 ] (248.44,10422.68) .. controls (253,10422.67) and (256.7,10426.36) .. (256.7,10430.91) .. controls (256.71,10435.47) and (253.02,10439.17) .. (248.47,10439.18) .. controls (243.91,10439.19) and (240.21,10435.5) .. (240.2,10430.94) .. controls (240.2,10426.39) and (243.88,10422.69) .. (248.44,10422.68) -- cycle ;
%Straight Lines [id:da350380584582642] 
\draw    (239.54,10371.28) -- (131.37,10371.25) ;
\draw [shift={(128.37,10371.25)}, rotate = 360.01] [fill={rgb, 255:red, 0; green, 0; blue, 0 }  ][line width=0.08]  [draw opacity=0] (8.93,-4.29) -- (0,0) -- (8.93,4.29) -- cycle    ;
%Straight Lines [id:da2877540997893018] 
\draw    (240.2,10430.94) -- (131.37,10431.56) ;
\draw [shift={(128.37,10431.58)}, rotate = 359.66999999999996] [fill={rgb, 255:red, 0; green, 0; blue, 0 }  ][line width=0.08]  [draw opacity=0] (8.93,-4.29) -- (0,0) -- (8.93,4.29) -- cycle    ;
%Straight Lines [id:da9450001576683873] 
\draw  [dash pattern={on 4.5pt off 4.5pt}]  (242,10376.5) -- (129.42,10424.65) ;
\draw [shift={(126.67,10425.83)}, rotate = 336.84000000000003] [fill={rgb, 255:red, 0; green, 0; blue, 0 }  ][line width=0.08]  [draw opacity=0] (8.93,-4.29) -- (0,0) -- (8.93,4.29) -- cycle    ;
%Curve Lines [id:da92648629352277] 
\draw  [dash pattern={on 4.5pt off 4.5pt}]  (127.33,10367.17) .. controls (151.14,10353.5) and (210.64,10353.02) .. (234.84,10366.93) ;
\draw [shift={(237.33,10368.5)}, rotate = 214.76] [fill={rgb, 255:red, 0; green, 0; blue, 0 }  ][line width=0.08]  [draw opacity=0] (8.93,-4.29) -- (0,0) -- (8.93,4.29) -- cycle    ;
%Curve Lines [id:da8139666079762311] 
\draw  [dash pattern={on 4.5pt off 4.5pt}]  (238.04,10435.51) .. controls (209.71,10448.85) and (151.13,10446.65) .. (127,10437) ;
\draw [shift={(241,10434)}, rotate = 150.95] [fill={rgb, 255:red, 0; green, 0; blue, 0 }  ][line width=0.08]  [draw opacity=0] (8.93,-4.29) -- (0,0) -- (8.93,4.29) -- cycle    ;
%Shape: Circle [id:dp4981368344543615] 
\draw  [fill={rgb, 255:red, 0; green, 0; blue, 0 }  ,fill opacity=1 ] (343.11,10364.01) .. controls (347.66,10364) and (351.36,10367.69) .. (351.37,10372.25) .. controls (351.38,10376.8) and (347.69,10380.5) .. (343.13,10380.51) .. controls (338.58,10380.52) and (334.88,10376.83) .. (334.87,10372.28) .. controls (334.86,10367.72) and (338.55,10364.02) .. (343.11,10364.01) -- cycle ;
%Shape: Circle [id:dp29196407050900675] 
\draw  [fill={rgb, 255:red, 255; green, 255; blue, 255 }  ,fill opacity=1 ] (343.11,10424.01) .. controls (347.66,10424) and (351.36,10427.69) .. (351.37,10432.25) .. controls (351.38,10436.8) and (347.69,10440.5) .. (343.13,10440.51) .. controls (338.58,10440.52) and (334.88,10436.83) .. (334.87,10432.28) .. controls (334.86,10427.72) and (338.55,10424.02) .. (343.11,10424.01) -- cycle ;
%Shape: Circle [id:dp9904684275572575] 
\draw  [fill={rgb, 255:red, 255; green, 255; blue, 255 }  ,fill opacity=1 ] (420.11,10364.01) .. controls (424.66,10364) and (428.36,10367.69) .. (428.37,10372.25) .. controls (428.38,10376.8) and (424.69,10380.5) .. (420.13,10380.51) .. controls (415.58,10380.52) and (411.88,10376.83) .. (411.87,10372.28) .. controls (411.86,10367.72) and (415.55,10364.02) .. (420.11,10364.01) -- cycle ;
%Shape: Circle [id:dp9702970509560862] 
\draw  [fill={rgb, 255:red, 255; green, 255; blue, 255 }  ,fill opacity=1 ] (420.11,10424.01) .. controls (424.66,10424) and (428.36,10427.69) .. (428.37,10432.25) .. controls (428.38,10436.8) and (424.69,10440.5) .. (420.13,10440.51) .. controls (415.58,10440.52) and (411.88,10436.83) .. (411.87,10432.28) .. controls (411.86,10427.72) and (415.55,10424.02) .. (420.11,10424.01) -- cycle ;
%Straight Lines [id:da8059243278344246] 
\draw    (411.87,10432.28) -- (354.37,10432.25) ;
\draw [shift={(351.37,10432.25)}, rotate = 360.03] [fill={rgb, 255:red, 0; green, 0; blue, 0 }  ][line width=0.08]  [draw opacity=0] (8.93,-4.29) -- (0,0) -- (8.93,4.29) -- cycle    ;
%Straight Lines [id:da03198237289412997] 
\draw  [dash pattern={on 4.5pt off 4.5pt}]  (343.11,10424.01) -- (343.13,10383.51) ;
\draw [shift={(343.13,10380.51)}, rotate = 450.04] [fill={rgb, 255:red, 0; green, 0; blue, 0 }  ][line width=0.08]  [draw opacity=0] (8.93,-4.29) -- (0,0) -- (8.93,4.29) -- cycle    ;
%Straight Lines [id:da7061202196915766] 
\draw  [dash pattern={on 4.5pt off 4.5pt}]  (420.11,10424.01) -- (420.13,10383.51) ;
\draw [shift={(420.13,10380.51)}, rotate = 450.04] [fill={rgb, 255:red, 0; green, 0; blue, 0 }  ][line width=0.08]  [draw opacity=0] (8.93,-4.29) -- (0,0) -- (8.93,4.29) -- cycle    ;
%Straight Lines [id:da7192285301948298] 
\draw  [dash pattern={on 4.5pt off 4.5pt}]  (351.37,10372.25) -- (408.87,10372.27) ;
\draw [shift={(411.87,10372.28)}, rotate = 180.03] [fill={rgb, 255:red, 0; green, 0; blue, 0 }  ][line width=0.08]  [draw opacity=0] (8.93,-4.29) -- (0,0) -- (8.93,4.29) -- cycle    ;
%Curve Lines [id:da8418086739179014] 
\draw [color={rgb, 255:red, 0; green, 0; blue, 0 }  ,draw opacity=1 ]   (343.13,10440.51) .. controls (340.11,10463.18) and (316.72,10455.59) .. (333.36,10436.91) ;
\draw [shift={(335.33,10434.83)}, rotate = 495.24] [fill={rgb, 255:red, 0; green, 0; blue, 0 }  ,fill opacity=1 ][line width=0.08]  [draw opacity=0] (8.93,-4.29) -- (0,0) -- (8.93,4.29) -- cycle    ;
%Curve Lines [id:da11275053308591598] 
\draw [color={rgb, 255:red, 0; green, 0; blue, 0 }  ,draw opacity=1 ]   (420.11,10364.01) .. controls (420,10338.91) and (447.01,10348.28) .. (428.77,10366.05) ;
\draw [shift={(426.61,10368.01)}, rotate = 319.66999999999996] [fill={rgb, 255:red, 0; green, 0; blue, 0 }  ,fill opacity=1 ][line width=0.08]  [draw opacity=0] (8.04,-3.86) -- (0,0) -- (8.04,3.86) -- cycle    ;
%Curve Lines [id:da6713737309914762] 
\draw [color={rgb, 255:red, 0; green, 0; blue, 0 }  ,draw opacity=1 ] [dash pattern={on 4.5pt off 4.5pt}]  (335.98,10368.79) .. controls (316.6,10357.35) and (332.01,10334.43) .. (342.18,10361.37) ;
\draw [shift={(343.11,10364.01)}, rotate = 251.95] [fill={rgb, 255:red, 0; green, 0; blue, 0 }  ,fill opacity=1 ][line width=0.08]  [draw opacity=0] (8.04,-3.86) -- (0,0) -- (8.04,3.86) -- cycle    ;
%Curve Lines [id:da5006113473695575] 
\draw [color={rgb, 255:red, 0; green, 0; blue, 0 }  ,draw opacity=1 ][line width=0.75]  [dash pattern={on 4.5pt off 4.5pt}]  (427.37,10435.25) .. controls (444.47,10445.68) and (430.87,10470.45) .. (421.04,10443.19) ;
\draw [shift={(420.13,10440.51)}, rotate = 432.6] [fill={rgb, 255:red, 0; green, 0; blue, 0 }  ,fill opacity=1 ][line width=0.08]  [draw opacity=0] (8.04,-3.86) -- (0,0) -- (8.04,3.86) -- cycle    ;
%Curve Lines [id:da9756853575700541] 
\draw  [dash pattern={on 4.5pt off 4.5pt}]  (409.94,10437.28) .. controls (390.14,10444.89) and (363.35,10443.65) .. (351,10437) ;
\draw [shift={(413,10436)}, rotate = 155.77] [fill={rgb, 255:red, 0; green, 0; blue, 0 }  ][line width=0.08]  [draw opacity=0] (8.93,-4.29) -- (0,0) -- (8.93,4.29) -- cycle    ;
%Curve Lines [id:da3764791564947594] 
\draw  [dash pattern={on 4.5pt off 4.5pt}]  (352.32,10376.78) .. controls (382.46,10384.27) and (402.53,10401.84) .. (417,10425) ;
\draw [shift={(349,10376)}, rotate = 12.34] [fill={rgb, 255:red, 0; green, 0; blue, 0 }  ][line width=0.08]  [draw opacity=0] (8.93,-4.29) -- (0,0) -- (8.93,4.29) -- cycle    ;
%Curve Lines [id:da8452515390568869] 
\draw    (411.07,10426.02) .. controls (384.33,10416.6) and (361.51,10397.33) .. (348,10379) ;
\draw [shift={(414,10427)}, rotate = 197.82] [fill={rgb, 255:red, 0; green, 0; blue, 0 }  ][line width=0.08]  [draw opacity=0] (8.93,-4.29) -- (0,0) -- (8.93,4.29) -- cycle    ;

% Text Node
\draw (120.12,10430.6) node  [font=\tiny,rotate=-1.27] [align=left] {$\displaystyle 4^{2}$};
% Text Node
\draw (120.12,10370.26) node  [font=\tiny,rotate=-1.27] [align=left] {$\displaystyle 4^{1}$};
% Text Node
\draw (247.79,10370.26) node  [font=\tiny,rotate=-1.27] [align=left] {$\displaystyle 3^{1}$};
% Text Node
\draw (248.45,10429.93) node  [font=\tiny,rotate=-1.27] [align=left] {$\displaystyle 3^{2}$};
% Text Node
\draw (183,10479.33) node   [align=left] {\begin{minipage}[lt]{15.64pt}\setlength\topsep{0pt}
(a)
\end{minipage}};
% Text Node
\draw (343.12,10371.26) node  [font=\tiny,rotate=-1.27] [align=left] {$\displaystyle \textcolor[rgb]{1,1,1}{1}$};
% Text Node
\draw (343.12,10431.26) node  [font=\tiny,rotate=-1.27] [align=left] {$\displaystyle 4$};
% Text Node
\draw (420.12,10371.26) node  [font=\tiny,rotate=-1.27] [align=left] {$\displaystyle 2$};
% Text Node
\draw (420.12,10432.26) node  [font=\tiny,rotate=-1.27] [align=left] {$\displaystyle 3$};
% Text Node
\draw (383,10479.33) node   [align=left] {\begin{minipage}[lt]{15.64pt}\setlength\topsep{0pt}
(b)
\end{minipage}};

\end{tikzpicture}

%%%%%%%%%%5
\caption{(a) Graph $G_{\{3,4\}}$; (b) $\Delta$-network graph associated with the global graph in Fig.~\ref{globgraph}.}
\label{final}
\end{figure}
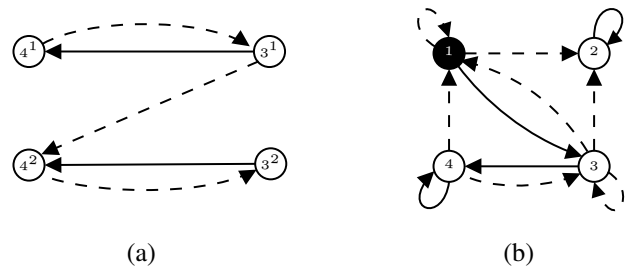

\begin{pro}
Consider the graph $G_{\{i,j\}}$, and let $Z=\emptyset$. Then, for $1\leq i,j\leq n$, where $i\neq j$,  $\mathcal{A}_{ij}$ has full row rank if and only if $V_i\subseteq \mathcal{D}(Z)$.
\label{thmf}
\end{pro}
 
\textbf{Proof.} For $l_i\leq l_j$, assume that all nodes of the graph $G(\mathcal{A}_{ij})=(V,E)$ with $V=\{v_1,\ldots, v_{l_i}\}$ are initially white, and perform a coloring process in this graph. Considering Theorem 10 in \cite{jia2020unifying}, one can conclude that $\mathcal{A}_{ij}$ has full row rank if and only if all nodes $v_1, \ldots, v_{{l_i}}$ in $G(\mathcal{A}_{ij})$ become finally black.  Now, consider the equivalent bipartite graph $G_{j\rightarrow i}=G_b(\mathcal{A}_{ij})$, and suppose that all vertices in $V_i\cup V_j$ are initially white. Assume that in some step  of the coloring process in the graph $G(\mathcal{A}_{ij})$, a white node $v_k$ has only one white strong out-neighbor $v_r$, and thus $v_k\rightarrow v_r$. One can see that equivalently, in the bipartite graph $G_{j\rightarrow i}$, there is a vertex $j^k\in V_j$ that has only a white strong out-neighbor, that is, $i^r\in V_i$, and we have $j^k\rightarrow i^r$. Thus, the nodes $v_1, \ldots, v_{l_i}$ in $G(\mathcal{A}_{ij})$ become finally black if and only if the vertices $i^1, \ldots, i^{l_i}$ in $G_{j\rightarrow i}$ are eventually black. Moreover, since the in-neighbors of any node in $V_i$ should  belong to $V_j$, one can see that during the coloring process, any vertex  $i^k$ in $G_{j\rightarrow i}$ becomes black if and only if it becomes black in $G_{\{i,j\}}$. Thus, $\mathcal{A}_{ij}$ has full row rank if and only if by the termination of the coloring process in $G_{\{i,j\}}$, all vertices in $V_i$ become black.   
\carre

\begin{ex}
Consider the global graph $\mathcal{G}$ in Fig.~\ref{globgraph}(b). As an example, $G_{b}(\mathcal{A}_{34})\cup G_{b}(\mathcal{A}_{43})$ is shown in Fig.~\ref{final}(a). One can see that by applying the coloring process in this graph, $V_4=\{4^1,4^2\}$ become black, while $V_3=\{3^1, 3^2\}$ cannot be black. Thus,  based on Proposition  \ref{thmf},  $\mathcal{A}_{43}$ has full row rank. Similarly, one can see that $\mathcal{A}_{12}=0$, $\mathcal{A}_{41}=0$, $\mathcal{A}_{32}=0$, $\mathcal{A}_{24}=0$, and $\mathcal{A}_{42}=0$. Moreover,  $\mathcal{A}_{31}$ and $\mathcal{A}_{43}$ have full row rank, while the nonzero pattern matrices $\mathcal{A}_{21}$, $\mathcal{A}_{13}$, $\mathcal{A}_{14}$,  $\mathcal{A}_{23}$, and $\mathcal{A}_{34}$ are rank-deficient. Now, let $\Delta=\{y\in \mathbb{C}\mid\Re(y)\geq 0\}$, and assume that subsystems 2 and 4 are stable. Then, one can provide the corresponding $\Delta$-network graph $G_{\mathcal{N}}^{\Delta}$ as shown in Fig.~\ref{final}(b). Note that since $V_C=\{1\}$ is a ZFS of $G_{\mathcal{N}}^{\Delta}$, the network is strongly structurally stabilizable.   
 \end{ex}

\section{Conclusions}\label{concl}

In this work, for a given nonempty set  $\Delta\subseteq \mathbb{C}$, we derived a graph-theoretic condition for  strong structural controllability of LTI networks, including dynamical subsystems, with respect to $\Delta$. The family of LTI networks consists of systems with the same zero/nonzero/arbitrary structure. However, in real applications, there might be some more information about the subsystems other than the zero/nonzero/arbitrary pattern of system matrices. For instance, one might know that some of the subsystems have no eigenvalue in $\Delta$. To deal with this more general case, we have defined a more restrictive family of LTI networks, for which the existing results on strong structural controllability seem conservative. Then, it was shown how one can provide a corresponding $\Delta$-network graph for a given set $\Delta$. In this setup, a correspondence between a set of control subsystems and a zero forcing set of a $\Delta$-network graph has been established. Moreover, given a set $\Delta$, that includes at least one real number, we have demonstrated that in N1DSs, the controllability condition is necessary and sufficient. Along the way, we have also derived structural conditions under which an N1DS admits no eigenvalue in $\Delta$ over all its pattern matrices (with extra spectrum features). In addition, the strong structural stabilizability of LTI networks was investigated as a special case. Finally, it has been shown how this work can generalize the previous works on strong structural controllability.

\bibliographystyle{IEEEtran}
\bibliography{library}

\end{document}